\documentclass[times,10pt]{amsart}
\usepackage{breqn}

\newcommand{\N}{\ensuremath{\mathbb{N}}}
\newcommand{\Z}{\ensuremath{\mathbb{Z}}}

\usepackage{hyperref}

\newcommand{\dd}{{\rm d}}
\newcommand{\F}{\ensuremath{\mathcal{F}}}
\newcommand{\I}{\ensuremath{\mathcal{I}}}
\newcommand{\R}{\ensuremath{\mathcal{R}}}
\newcommand{\n}{\ensuremath{\infty}}

\newcommand{\resp}{{\it resp. }}
\renewcommand{\phi}{\varphi}
\renewcommand{\cong}{\simeq}
\renewcommand{\resp}{{\it resp. }}
\renewcommand{\ell}{{\rm Pf}}
\usepackage{amssymb}
\usepackage{geometry}
\usepackage{amsthm}
\usepackage{color}
\usepackage{xcolor}
\usepackage{tikz}
\usetikzlibrary{positioning}
 \newtheorem{thm}{Theorem}[section]
 \newtheorem{defn}[thm]{Definition}
 \newtheorem{cor}[thm]{Corollary}
 \newtheorem{prop}[thm]{Proposition}
 \newtheorem{lemma}[thm]{Lemma}
 \newtheorem{rem}[thm]{Remark}

\date{}
\title{Kato-Milne cohomology over rational function fields in characteristic $p>0$}
\author{Ahmed Laghribi}
\author{Trisha Maiti}
\begin{document}
\maketitle

\begin{abstract}
     Let $\F$ be a field 
 of characteristic $p>0$. In this note, we describe the Kato-Milne cohomology of  the  rational function field $\F(x)$ in one variable   $x$ over $\F$.  In particular, we extend  the  Milnor-Scharlau exact sequence \cite[Theorem 5.2]{Laghribi-Maiti}  to  characteristic $p$. As an application of this exact sequence, we further obtain  the norm theorem \cite[Theorem A.3]{Appendix} in characteristic $p$.
\end{abstract}

\section{Introduction}\label{s1}
   { The well-known  result \cite[Theorem 5.3]{M} due to Milnor describes the Witt group of non-degenerate bilinear forms $W(F(x))$, where  $F$ is  a field of  characteristic  different from  $2$.  By proving a similar exact sequence in \cite[Theorem 6.2]{AJ}, Aravire and Jacob  described the Witt group of nonsingular quadratic forms $W_q(F(x))$, where  $F$ has characteristic  $2$.  Since $W_q(F(x))$ is   related to  the Kato-Milne cohomology $ H_2^{m+1}(F(x))$ by Kato's isomorphisms \cite[equation (1.2)]{ex},  we  have already  extended the exact sequence of  Aravire-Jacob in the setting of Kato-Milne cohomology in characteristic $2$, \cite[Theorem 5.2]{Laghribi-Maiti}.  Here we extend  the works of  
      \cite{ex} and \cite{Laghribi-Maiti} to characteristic $p>0$.}  {More precisely, we will prove} Theorem \ref{myth}, which is a kind of local-global principle for the Kato-Milne cohomology groups of rational function fields in one variable in  positive characteristic.  {The difference compared to the characteristic $2$ case lies in the definition of the transfer maps $ s_q^*$, where $q\in F[x]$ is a monic irreducible polynomial and  $F$ is a field of characteristic $p$.}  Since Kato's isomorphism \cite[equation (1.2)]{ex} does not exist in characteristic  $p>2$, we can not define $s_q^*$ using Scharlau's transfer as done {in characteristic $2$ case}. {To overcome this situation, we will consider the trace map ${\rm Tr}: F(q)\xrightarrow{} F$ to construct $s_q^*$ when $q$ is separable, and for inseparable polynomial the transfer  $s_{F(q)/F}: \Omega_{F(q)}^{m}\xrightarrow{} \Omega_{F}^{m-1}$ constructed in \cite[Proposition 3.4]{ALO} will be used}. {Other crucial computations to establish the transfer maps were necessary like lemmas \ref{l81} and \ref{l82}.} We start by giving some preliminary  definitions regarding to differential forms and  Kato-Milne cohomology.

      {Let $F$ be a field of characteristic $p>0$ and $\Omega_F^1$ be the $F$-vector space generated by the symbols $\dd a$ for $a\in F$ subject to the relations $\dd (a+b)= \dd a+ \dd b$ and $\dd(ab)= a\dd b+b\dd a$ for all  $a$, $b\in F$. For any positive integer $m$, let $\Omega_F^m:= \wedge^m\Omega_F^1$ be the $m$-th exterior product of $\Omega_F^1$. {We take $\Omega_F^0=F$ and $\Omega_F^{m}=0$ for $m<0$}. The  set $\Omega_F^m$ is called the space of $m$-differential forms over $F$  for all $m\in \Z$.  
Now  the differential operator $\dd: \Omega_F^{m-1}\xrightarrow{} \Omega_F^m$ is defined by  $$  \dd ( a\dd a_1\wedge\ldots\wedge \dd a_{m-1})= \dd a\wedge \dd a_1\wedge \ldots\wedge \dd a_{m-1}.$$ 
{A} $m$-differential form which lies in $\dd(\Omega_F^{m-1})$ for some $m\in \Z$ is called an exact form. Let $ \mathcal{B}:=\{e_i\mid\hspace{1mm} i\in I\}$ be a {$p$-basis} of $F$, where  $I$ is an ordered set. {This means that the set $\big\{\prod_{i\in I}e_i^{\epsilon_i}\mid \epsilon_i\in {\{0,\ldots,p-1\}}\; \text{and}\;\epsilon_i=0\;\text{ for almost all }\; i\in I\big\}$ is an {$F^p$-basis} of $F$}.  Let us consider the set of maps  $\Sigma_{m,F}=\{ \sigma: \{1,\ldots,m\}\xrightarrow{}I \mid \sigma(j)<\sigma(k)\; \text{when}\;j<k \}.$ It is well-known  that $\left\{\frac{\dd e_{\sigma}}{e_{\sigma}}:= \frac{\dd e_{\sigma(1)}}{e_{\sigma(1)}}\wedge \frac{\mathrm{d}e_{\sigma(2)}}{e_{\sigma(2)}}\wedge \ldots \wedge \frac{\mathrm{d}e_{\sigma(m)}}{e_{\sigma(m)}} \mid  \sigma\in \Sigma_{m,F}\right\}$ is an $F$-basis of $\Omega_F^m$.   For a fixed {$p$-basis} $\mathcal{B}$, we define the Artin-Schreier operator  $\wp: \Omega_F^m\xrightarrow{} \Omega_F^m$ by $$ \wp\left(  \sum_{\sigma} c_{\sigma} \frac{\dd e_{\sigma}}{e_{\sigma}}\right)= \sum_{\sigma} {(c^p_{\sigma}-c_{\sigma})} \frac{\dd e_{\sigma}}{e_{\sigma}}.$$ 
For $m=0$,  {this operator coincides with the usual} Artin-Schreier map $ \wp: F\xrightarrow{}F$,  {$x\mapsto x^p-x$} for all $x\in F$.
By a result of Cartier, if we change the {$p$-basis} of $F$, the image of $\wp$ will be changed by an exact form.  Therefore,  the  map $\wp: \Omega_F^m\xrightarrow{} \Omega_F^m/{\dd \Omega_F^{m-1}}$ given by  $$ \wp\left(  \sum_{\sigma} c_{\sigma} \frac{\dd e_{\sigma}}{e_{\sigma}}\right)= \sum_{\sigma} {(c^p_{\sigma}-c_{\sigma})} \frac{\dd e_{\sigma}}{e_{\sigma}}+\dd \Omega_F^{m-1} $$does not depend on the {$p$-basis} of $F$.  The cokernel of this map  is  called the  Kato-Milne cohomology group (of degree $m$) over $F$ and  denoted by ${H_p^{m+1}(F)}$. {Moreover, let $\nu_F(m)$ denote the kernel of the map $\wp$. It is clear that we have$${H_p^{m+1}(F)}= {\Omega^m_F}/\left({\dd\Omega_F^{m-1}+\wp(\Omega_F^m)}\right)$$ where $\wp$ is the  Artin-Shreier operator with respect to a fixed {$p$-basis} of $F$.} }

    Let    $v: F\xrightarrow{} \Z\cup\{\n\}$ be a discrete valuation of $F$ with $q$ as a uniformizer. We use $\overline{F}$ and $K$ to denote the  residue field  and the completion  of $F$, respectively. Let $ \{ e_i\mid i\in \I\}$ be a $p$-basis of $\overline{F}$ and $ \{ f_i \mid i\in \I\}$ a subset of $ K$, where  each $f_i$ is a lifting of $e_i$.  By Theorem  \ref{imp},   there exists a  Teichm\"uller lifting $\alpha: \overline{F}\xrightarrow{} K$ such that $ \alpha(e_i)=f_i$, for all $i\in \I$.  Now we define the group $$ W_1(H_p^{m+1}(K)):= {\rm coker} ( \alpha': H_p^{m+1}( \overline{F})\xrightarrow{} H_p^{m+1}(K))$$ where $\alpha'$ is  induced from the Teichm\"uller lifting $\alpha$.  Moreover, by  Corollary \ref{T2c}  the group $ W_1(H_p^{m+1}(K))$ is independent of the choice of a Teichm\"uller lifting. Let $ \beta_1: H_p^{m+1}( {F})\xrightarrow{} H_p^{m+1}(K)$ be the map induced by the inclusion $ {F}\xrightarrow{} K$, and $ \beta_2: H_p^{m+1}( K)\xrightarrow{}  W_1(H_p^{m+1}(K))$ be the projection map. Then, the {(second)} residue map $\partial_v: H_p^{m+1}({F})\xrightarrow{} W_1(H_p^{m+1}(K))$ is defined by $$\partial_v := \beta_2\circ \beta_1.$$ 
    
Now let {$F=\F(x)$ be the  rational function field in one variable $x$} over a field $\F$ of characteristic $p>0$.   Let $\overline{F}_q$ and $F_q$ be the residue field and the completion of $F$ with respect to the $q$-adic valuation, where $q=1/x$ or $q\in \F[x]$ a monic irreducible polynomial. Now it follows from Theorem \ref{fsec} that    the group $H_p^m(\overline{F}_q)$ is isomorphic to   the quotient group of $W_1(H_p^{m+1}(F_q))$. Let $t_q^*: H_p^{m}(\overline{F}_q)\xrightarrow{} H_p^m(\F)$ be the map  {induced by}  ${\rm Tr_*}: \Omega^{m-1}_{\overline{F}_q}\xrightarrow{}\Omega^{m-1}_{\F}$  as defined in \cite[Corollary 2.9]{AAB} for $q$ separable, and induced by the transfert map $s_{\overline{F}_q/\F}: \Omega_{\overline{F}_q}^{m-1}\xrightarrow{} \Omega_{\F}^{m-1}$  defined in \cite{ALO} for $q$ inseparable. Now we define the transfer map $ s_q^*: W_1(H_p^{m+1}(F_q))\xrightarrow{} H_p^{m}(\F)$   as the composition  map $ W_1(H_p^{m+1}(F_q))\xrightarrow{} H_p^m(\overline{F}_q)\xrightarrow{t_q^*}  H_p^m(\F)$, (see Definition \ref{dfntr} for details).
     Finally our main result 
is the following:
\medskip

\begin{thm}\label{myth} Let $\mathcal{F}$ be a field of characteristic {$p>0$}, and $F=\mathcal{F}(x)$ the rational function field in one variable $x$ over $\mathcal{F}$. Then, \textcolor{black}{the following   sequence is    exact}:\[ 0\xrightarrow{} H_p^{m+1}(\mathcal{F})\xrightarrow{\textcolor{black}{i}}  H_p^{m+1}(F)\xrightarrow{{\oplus_{q}\partial_q}}\displaystyle{ \oplus_{q}W_1(H_p^{m+1}(F_q))} \xrightarrow{\oplus s_q^*} H_p^m(\F)\xrightarrow{}0\] where $\textcolor{black}{i}$ is induced by 
 the inclusion $ \mathcal{F}\xrightarrow{}F$, and $q$ varies  over ${1}/{x}$ and  monic irreducible polynomials  of $\mathcal{F}[x]$. Moreover, $\partial_q$ is the second residue map with respect to the $q$-adic valuation.
\end{thm}
To prove this theorem, we follow the same approach  {as in  \cite{ex} and \cite{Laghribi-Maiti}}.
Now  the rest of this {note is} outlined as follows.       In Section \ref{5s2} we decompose the group $W_1(H_p^{m+1}(F))$ into the required form, where $F$ is a complete valued field. To this end we extend \cite[Lemma 2.3]{ex}  and \cite[Proposition 1.5]{AJ} to characteristic $p$ by proving Lemma  \ref{lde} and Lemma \ref{lu}, respectively.
Then, in Section \ref{5sts} we  define a filtration 
 $L_0\subset L_1\subset \ldots\subset H_p^{m+1}(\F(x))$ and some subgroups $S_q$ of $H_p^{m+1}(\F(x))$ (see Definition \ref{definse}) such that  $ H_p^{m+1}(\F(x))=\sum_{q} (S_q+S_q\wedge \frac{\dd x}{x})$, where $q$ varies over monic  irreducible polynomials over $\F$. To establish the required relations between $ S_q\wedge \frac{\dd q}{q}$ and $S_q\wedge \frac{\dd x}{x}$,   we prove Lemma \ref{linv} which extends \cite[Lemma 2.6]{AJ} to characteristic $p$.   These relations of Lemma \ref{linv} and  the decomposition of  $W_1(H_p^{m+1}(F_q))$ {given in} Theorem \ref{thd2}   induce the  following exact sequence \begin{equation}\label{ei1}
      0\xrightarrow{}L_0\xrightarrow{i_1}H_p^{m+1}(\F(x))\xrightarrow{\oplus \partial_q}\oplus_{q} W_1(H_p^{m+1}(F_q))\xrightarrow{}0\end{equation} where $i_1$ is the inclusion map  and   $q$ varies over monic irreducible polynomials of $\F[x]$ (Theorem \ref{fsec}).  Then, using the decomposition of $W_1(H_p^{m+1}(F_{1/x}))$, we get the following exact sequence  in Theorem \ref{ssec} \begin{equation}\label{ei2} 0\xrightarrow{}H_p^{m+1}(\F)\xrightarrow{i_2} L_0\xrightarrow{\partial_{1/x}} W_1(H_p^{m+1}(F_{1/x}))\xrightarrow{\eta} H_p^m(\F)\xrightarrow{}0
           \end{equation} where $i_2$ is induced by inclusion and $\eta(\phi)=\phi'$, which comes from Theorem \ref{thd1}. \\
           By (\ref{ei1}) and  (\ref{ei2}),  we get  Theorem \ref{thmil} which  extends  \cite[Theorem 5.3]{Milnor} in the setting of Kato-Milne cohomology. \\
           To introduce the transfer maps $s_q^*$,   we first establish an alternative  decomposition of $W_1(H_p^{m+1}(F_q))$  in Theorem \ref{thd3} which extends  \cite[Theorem 2.3]{Laghribi-Maiti} to characteristic $p$.       In Section \ref{5srl} we prove  the reciprocity law for $H_p^{m+1}(\F(x))$  i.e.,  $\oplus_{q} (s_q^*\circ \partial_q)(\phi)=0$ for all $\phi\in H_p^{m+1}(\F(x))$. Consequently, ${\rm Im}(\oplus \partial_q^*)\subset \ker(\oplus_q s^*_q)$. The  converse follows from a straightforward computation, which proves Theorem \ref{myth}.
             Finally as an application of Theorem \ref{myth}, we prove the norm theorem for Kato-Milne cohomology in characteristic $p$ (Theorem \ref{n2th}). Since the arguments are analogous to those presented in papers \cite{ex} and \cite{Laghribi-Maiti}, we omit most of the proofs.

    \section{The group $W_1(H_p^{m+1}(F))$}\label{5s2}
    
    Let  $F$ be a field of characteristic $p$ with   $B:=\{ t_i\mid i\in \I\}$  as  a $p$-basis, where $\I$ is an   ordered set.  Let   $T:= \oplus_{|\I|} \Z/p\Z$  \label{b} {which is ordered by the lexicographic ordering as follows: For} any $I= ( I_i)_{i\in \I}$ and $J=(J_i)_{i\in \I}\in T$, we say $I<J$ iff  there exists some $j\in \I$ such that $ I_j<J_j$ and $I_i= J_i$ for all $i>j$. For each  $I\in T$, the notation $t^I$ denotes  the element $\prod_{i\in \I} t_i^{I_i}$. Thus, each element $f\in F$ can be written uniquely as $$  f= \sum_{I\in T} t^I f_I^p,$$ where $f_I\in F$ and $f_I=0$ for {almost all} $I\in T$.
    
    For the {rest} of this section, we suppose   $ F$ a complete  valued field  of characteristic $p$ with respect to the valuation $ v: F\xrightarrow{} \Z\cup\{\n\}$ with $\pi$ as a uniformizer.  We use the notations $A_F$ and $\overline{F}$ to denote the valuation ring and the residue field of $F$, respectively. Let $ \{ t_i\mid i\in \mathcal{J}\} \subset A_F$ be such that $ \{ \overline{t_i}\mid i\in \mathcal{J}\}$ {is a } $p$-basis of $\overline{F}$. {By \cite[Remark 1.1.1.12.]{LX},} $\{ t_i\mid i\in \I\}$ will be a $p$-basis of $F$, where $ \I=\mathcal{J}\cup \{\infty\}$ and $t_{\infty}=\pi $.   The ordering of $\mathcal{J}$ is extended to $\I$ with $ i<\infty$ for all $i\in \mathcal{J}$.
     By Theorem \ref{imp}, there is a unique lifting monomorphism $\alpha: \overline{F}\xrightarrow{} F$ such that $\alpha( \overline{t_i})= t_i$ for all $i\in \mathcal{J}$.
    For all ${I= ( I_i)_{i\in \I}}\in T:= \oplus_{|\I|} \Z/p\Z$, we have  $v( t^I)= I_{\infty}$. For any $I\in T$ and $m\in \N_0$, we define $ \tilde{I}:=\{ i\in \I\mid I_i\neq 0\}$, $l(I):= {\rm Max} \tilde{I}$   and $(T)_m:=\{ I\in T\mid  I_i\in \{\overline{0},\overline{1}\}\text{  for all $i\in \I$ and } |\tilde{I}|=m\}$. Let $\overline{T}:= \oplus_{|\mathcal{J}|} \Z/p\Z\subset T$. Then for any $m\in \N_0$, we have \label{p1}$(\overline{T})_m\subset (T)_m$.  {Moreover, if} $\overline{F}\subset F$, then we take $ \overline{t}_i=t_i$ for all $i\in \mathcal{J}$.
    
    \begin{defn}\label{dr} Let $F$ be a complete valued field as described above such that $\overline{F}\subset F$.
        Then, we define the subgroup  $\R$  of $H_p^{m+1}(F)$ generated by the following elements $$   \sum_{I\in (T)_m}\sum_{J\in T\mid l(J)\notin \tilde{I}} t^J r^p_{I,J} \frac{\dd t_I}{t_I}$$   where $r_{I,J}\in \pi^{-1}\overline{F}[\pi^{-1}]$   and  $\frac{\dd t_I}{t_I}$ denotes the element $ \frac{\dd t_{i_1}}{t_{i_1}}\wedge \ldots\wedge \frac{\dd t_{i_m}}{t_{i_m}}\in \Omega_F^m$ such that $\tilde{I}= \{ i_1,\ldots,i_m\}$  and $i_1<\ldots<i_m$.  
    \end{defn}

    \begin{rem}\label{rv}(1) If $a\in F$ and $v(a)>0$, then by Hensel's lemma  $a\in \wp(F):= \{ x^p-x\mid x\in F\}$.\\
\noindent{(2)} If $v(a)<0$, then $a\notin \wp(F)$ if one of  the following  conditions {holds}:\\(i) $ v(a)$ is not divisible by $p$.\\(ii) The leading coefficient of $a$ is not {in $F^p$}.\\ 
Therefore, the separable parts of the generators of   $\R$ are not   in $\wp(F)$.
\end{rem}
\begin{lemma}\label{lde}
    Let $\{ t_i\mid i\in \I\}$ be a $p$-basis of $F$. Then, each element $\phi\in H_p^{m+1}(F)$ can be expressed as follows:$$\phi= \sum_{I\in (T)_m} A_I \frac{\dd t_I}{t_I}$$ where $ A_I= \sum\limits_{J\in T\atop J=0 \text{ or } l(J)\notin \tilde{I}} t^J r^p_{I,J}$ with  $r_{I,J}\in F$ and $r_{I,J}=0$ for {almost} all $I,J$.
\end{lemma}
    
    \begin{proof}
         Since $\{ t_i\mid i\in \I\}$ is a $p$-basis of $F$, {the group $ H_p^{m+1}(F)$ is additively generated \textcolor{blue}{by} the  elements  $ t^J r^p_{I,J} \frac{\dd t_I}{t_I}$,} where  $ I\in (T)_m$,  $ J\in T$ and $r_{I,J}\in F$. If $J=0$ or $l(J)\notin \tilde{I}$, then we are done. \\
         Let  $J\neq 0$ be such that   $ l(J)=i_l$ and $ J_{i_l}=p-r$ for some $r\in \{1,\ldots,p-1\}$. If $i_l\in \tilde{I}$, then we  have:
         \begin{align*}
             t^J r^p_{I,J} \frac{\dd t_I}{t_I}&= t^J r^p_{I,J} \frac{\dd t_{i_1}}{t_{i_1}}\wedge \ldots\wedge \frac{\dd t_{i_l}}{t_{i_l}}\wedge \ldots\wedge \frac{\dd t_{i_m}}{t_{i_m}}\\&= \frac{t^J  r^p_{I,J}}{r^p} \frac{\dd t_{i_1}}{t_{i_1}}\wedge \ldots\wedge \frac{\dd t_{i_l}^r}{t^r_{i_l}}\wedge \ldots\wedge \frac{\dd t_{i_m}}{t_{i_m}}\quad \text{ (since $ r^p=r$ and $ \frac{\dd t^r_{i_l}}{t^r_{i_l}}= r\frac{\dd t_{i_l}}{t_{i_l}}$)}\\&= t^J {a^p_{I,J}} \frac{\dd t_{i_1}}{t_{i_1}}\wedge\ldots\wedge \frac{\dd (t^r_{i_l}t^J)}{t^r_{i_l}t^J}\wedge\ldots\wedge \frac{\dd t_{i_m}}{t_{i_m}}\quad \text{(where $a_{I,J}= \frac{r_{I,J}}{p}$)}\\&=\sum_{k\in \tilde{J}\notin \tilde{I}} t^Ja^p_{I,J} \frac{\dd t_{i_1}}{t_{i_1}}\wedge \ldots\wedge \frac{\dd t_k^{J_k}}{t_k^{J_k}}\wedge \ldots\wedge \frac{\dd t_{i_m}}{t_{i_m}}\\&=\sum_{k\in \tilde{J}\notin \tilde{I}}t^J (J_k a_{I,J})^p \frac{\dd t_{i_1}}{t_{i_1}}\wedge \ldots\wedge \frac{\dd t_k}{t_k}\wedge \ldots\wedge \frac{\dd t_{i_m}}{t_{i_m}}\\&= \sum_{I'\in (T)_m\mid l(J)\notin \tilde{I'}}  t^J a^p_{I',J}\frac{\dd t_{I'}}{t_{I'}}
         \end{align*}  where   $\tilde{I'}=\{ i_1,\ldots,i_{l-1},i_{l+1},  \ldots,i_m\}\cup \{k\}$ for $k\in \tilde{J}\notin \tilde{I}$. Moreover,  $ a_{I',J}=  c_kJ_k a_{I,J}$, where  $c_k=1$ or $-1$. Therefore,  $a_{I',J}= c_{I',J}r_{I,J}$ for some $c_{I',J}\in \Z/ p\Z$. Thus the lemma follows.
          \end{proof}
          The following remark is   a direct consequence of  the previous lemma. 
     \begin{rem}\label{rim}
         Let   $\{ t_i\mid i\in \I\}$  be a $p$-basis of $F$ and  $\displaystyle{ A= \sum_{I\in (T)_m}\sum_{J\in T\mid l(J)\in \tilde{I}} t^J  r^p_{I,J}\frac{\dd t_I}{t_I}\in H_p^{m+1}(F)}$. Then,  repeating the proof  of Lemma \ref{lde}, $A$ can be decomposed as $$ \sum_{I'\in (T)_m}\sum_{J\in T\mid l(J)\notin \tilde{I'}} t^J  r^p_{I',J} \frac{\dd t_{I'}}{t_{I'}} $$  where $ r_{I',J}= c_{I',J} r_{I,J}$ for some $c_{I',J}\in \Z/p\Z$.
     \end{rem}   
     Now we  extend the result \cite[Lemma 1.4]{AJ} {to Kato-Milne cohomology in characteristic $p>0$.}

     \begin{lemma}\label{lnu}
         Let $\{ t_i\mid i\in \I\}$ be a $p$-basis of $F$. If $ \phi= \sum_{I\in (T)_m} a_I \frac{\dd t_I}{ t_I}=0\in  H_p^{m+1}(F)$,  then for all $I\in (T)_m$, we have $$  a_I\in \wp(F) + \sum_{J\in T\mid \tilde{J}\cap \tilde{I}\neq \emptyset} t^J F^p.$$
     \end{lemma}
     \begin{proof}
         Let $\I':=\{ i\in \I\mid  i\in \tilde{I}\text{ and } a_I\neq 0\}$ {which} is a finite subset of $\I$. Now we fix an element $K\in (T)_m$ such that $a_K\neq 0$. Let $\tilde{K}= \{K_1,\ldots,K_m\}$ be such that $K_1<\ldots<K_m$.  We define the field extension $L:= F( \sqrt[p]{t_i}\mid i\in \I'\setminus  \tilde{K})$.  Then we have  $$ \phi_L= a_K \frac{\dd t_K}{t_K}=0\in H_p^{m+1}(L)$$ Let $\I''$ be an ordered set whose elements are the elements of $\I$, with  the following  ordering  $$ i<  j, \quad \text{ for all } i\in \I\setminus (\I'\setminus \tilde{K}),\hspace{1mm} j\in \I'\setminus \tilde{K}$$ and  in other cases the  ordering in $\I''$ is same as $\I$. Hence, $\{ t_i'\mid i\in \I''\}$ is a $p$-basis of $L$, where $ t_i'=t_i$ for all $i\in \I\setminus ( \I'\setminus \tilde{K})$ and $ t_i'=  \sqrt[p]{t_i}$ for all $ i\in \I'\setminus \tilde{K}$. Since $a_K\frac{\dd t_K}{t_K}=a_K\frac{\dd t'_K}{t'_K}=0\in H_p^{m+1}(L)$, it follows from \cite[Proposition 3.3]{AAB} that  $$  a_K= \wp(u) + \sum_{i=1}^{m} \sum_{j=1}^{p-1} M_{i,j} (t'_{K_i})^j\in L$$ where $ u\in L$   and $M_{i,j}\in L^p( t'_r\mid r< K_i\in \I'') \subset F$. Since $ t'_{K_i}=t_{K_i}\in F$ for all $i$ and $L^p\subset F$, we have $ u\in F$.  Since $L^p=F^p( t_i\mid i\in \I'\setminus \tilde{K})$, we get $$  a_K\in \wp(F)+ \sum_{J\in T\mid \tilde{J}\cap \tilde{K}\neq \emptyset}t^J F^p$$ Thus the lemma is proved.\end{proof}

     The following lemma   extends  \cite[Proposition 1.5]{AJ} {to Kato-Milne cohomology in characteristic $p>0$}. 
\begin{lemma}\label{lu} {Let $v:F\xrightarrow{}\Z\cup\{\n\}$ be a valuation for which $F$ is complete, and let $\{ t_i\mid i\in \I\}$ be a $p$-basis.} Let $\overline{F}\subset F$ be the residue field and $\pi$ be a uniformizer of $F$.  Suppose $$ \phi= \sum_{I\in (T)_m} A_I \frac{\dd t_I}{ t_I}=0\in H_p^{m+1}(F)$$ where $  A_I = \sum\limits_{J\in T\atop l(J)\notin \tilde{I}}\sum\limits_{r\geq 1}   \frac{t^Ja_{r,I,J}^p}{\pi^{pr}} $ and $ a_{r,I,J}\in \overline{F}$. Moreover, $a_{r,I,J}=0$ for {almost} all $I$, $J$, $r$. Then $A_I=0$ for all $I$.  
\end{lemma}
\begin{proof}
     Let $ S:= \{ I\in (T)_m\mid  A_I\neq 0\}\subset T$ and suppose $S$ {is} nonempty. Let $M$ be the maximum element of $S$ with respect to the  ordering induced from $T$.  Now we can write $ A_M= \sum\limits_{J\in T\atop l(J)\notin  \tilde{M}} A_{M,J} $, where  \begin{equation}\label{5etr} A_{M,J}= \sum\limits_{r\geq 1} \frac{{t^Ja^p_{r,M,J}}}{\pi^{pr}},\quad \text{where } a_{r,M,J}\in \overline{F}.
             \end{equation}   Since $A_M\neq 0$, there exists an element $J\in {T}$ such that $l(J)\notin \tilde{M}$ and $A_{M,J}\neq 0$.   We define a new subset $\{t_i'\}_{i\in \I}$ of $F$ such that  $t_{l(J)}'= t^J$ and $t_i'=t_i$ for all $i\neq l(J)\in \I$.   Then, the set $\{ t'_i\mid i\in \I\}$ is a $p$-basis of $F$.   Since  $ l(J)\notin \tilde{M}$, we have $ \frac{ \dd {t_M'}}{t'_{M}}  =\frac{\dd {t_M}}{t_M} $. Let $K\in (T)_m$ be such that $l(J)\in \tilde{K}$. Then $\tilde{K}= \{l(J)\}\cup\tilde{Q}$, where $Q\in (T)_{m-1}$. Now in $H_p^{m+1}(F)$  we have\begin{align*}
     \frac{\dd t'_K}{t'_K}&= c_K \frac{ \dd t'_Q}{t'_Q} \wedge \frac{ \dd t'_{l(J)}}{t'_{l(J)}}\quad \text{(where $c_K\in \{-1,1\}$)}\\&=  c_K \frac { \dd t_Q} {t_Q}\wedge \frac {\dd t^J} {t^J}\\&=  \sum_{j\in \tilde{J}} c_KJ_j\frac{\dd  t_{Q}} {t_Q}\wedge \frac{\dd t_j}{t_j}\quad \text{(where $J_j\in \Z/{p\Z}$)} \\&=   J_{l(J)} \frac{\dd t_K}{t_K} + \sum_{j\in \tilde{J}\setminus (\tilde{Q}\cup \{l(J\})} a_{L_j}\frac{ \dd t_{L_j}} {t_{L_j}}\quad (\text{where } L_j\in (T)_m \text{ and }\tilde{L}_j=\tilde{Q}\cup \{j\},\hspace{1mm} a_{L_j}\in \Z/{p\Z})\\&= J_{l(J)}\frac{\dd t_K}{t_K}+ \sum_{j\in\tilde{J} \setminus \tilde{K}} a_{L_j}\frac{\dd t'_{L_j}}{t'_{L_j}} \quad \text{(since  $ l(J)\notin \tilde{L_j}$)}  \\&= J_{l(J)}\frac{\dd t_K}{t_K} + \sum_{j\in \tilde{J}\setminus \tilde{K}\atop L_j<K}  a_{L_j} \frac{\dd t'_{L_j}}{t'_{L_j}}. \end{align*}   The last equality holds because $\tilde{L_j}= \tilde{Q}\cup\{j\}$ and  $j<l(J)$ for all $j\in \tilde{J}\setminus \tilde{K}$. {Therefore, }\begin{align}\label{5ft} \frac {\dd t_K}{t_K} &= b_K \frac{\dd t'_K}{t'_K} +\sum_{j\in \tilde{J}\setminus \tilde{K}\atop L_j<K} b_{L_j}  \frac{\dd t'_{L_j}}{t'_{L_j}}, 
     \end{align} where  $ b_K,$ $ b_{L_j}\in \Z/{p\Z}$. Hence, we get in $H_p^{m+1}(F)$  \begin{align*}
         \phi&= \sum_{I\in (T)_m} A_I \frac{\dd t_I}{t_I}= \sum_{I<M\in (T)_m} A_I\frac{\dd t_I}{t_I}+ A_M\frac{\dd t_M}{t_M}= \sum_{I<M\in (T)_m} A'_I\frac{\dd t'_I}{t'_I}+  A_M\frac{\dd t'_M}{t'_M} \quad \text{(by (\ref{5ft}))}
     \end{align*} where $A'_I\in F$. By  Lemma \ref{lnu}, with respect to the $p$-basis of $\{t_i'\mid i\in \I\}$ of $F$, we have  \begin{equation*} A_M\in \wp(F)+ \sum_{I\in T\atop \tilde{I}\cap \tilde{M}\neq \emptyset } {t'}^IF^p.\end{equation*}  Now $ {t'}^I= t^J$, if $\tilde{I}=\{l(J)\}$ and $I_{l(J)}=1$, but then $\tilde{I}\cap \tilde{M}=\emptyset$. Therefore,  \begin{equation*}A_{M,J}\in \wp(F)- \sum\limits_{ K\neq J\in T\atop l(K)\notin \tilde{M}} A_{M,K}+ \sum\limits_{I\in T\mid  I\neq J\atop I\neq 0}  t^IF^p.\end{equation*} Since $-1= (-1)^p\in F^p$, we have  \begin{equation}\label{5sr}A_{M,J}\in \wp(F)+ \sum\limits_{I\in T\mid  I\neq J\atop I\neq 0}  t^IF^p.\end{equation} By the expression (\ref{5etr}) of  $ A_{M,J}$, we have  $v(A_{M,J})=s<0$. Now {$p\mid s$} if and only if $J_{\infty}=0$.  For any $K\in T$, let $K'$ denote the element of $T$ such that $K'_{i}=K_i$ for all $i\neq \infty\in \I$ and $K'_{\infty}=0$.  By (\ref{5etr}), we have \begin{equation}\label{5fr}
         \overline{{\pi^{-s}} A_{M,J}}\in  t^{J'} \overline{F}^p. 
     \end{equation}{ By (\ref{5sr}), we have    $ A_{M,J}= \wp(u)+ \sum\limits_{I\in T\mid I\neq J\atop I\neq 0} t^Ib^p_{I,J}$, where $u\in F$ and $b_{I,J}\in F$. Therefore,  \begin{equation*}
         \overline{{\pi^{-s}}  A_{M,J}} \in  \begin{cases}
             \overline{F}^p +& \sum\limits_{I\in T, I\neq J\atop I\neq 0, I_{\infty}=0} t^I \overline{F}^p\quad \text{when $p\mid s$,}\\ &\sum\limits_{I\in T\mid I'\neq J'\atop I_{\infty}=J_{\infty}\neq 0 }t^{I'} \overline{F}^p\quad\text{ when $p\nmid s$.}  
         \end{cases}
     \end{equation*} }Since $ \{t_i\mid i\neq \n\in {\I}\}$ is a $p$-basis of $\overline{F}$ and the above two expressions of $ \overline{p^{-s}  A_{M,J}}$ are different, $A_{M,J}=0$. Thus, $a_I=0$ for all $I$ which proves the lemma.
    
\end{proof}
  The following theorem  describes the group $W_1(H_p^{m+1}(F))$, where $F$ is a complete valued field of characteristic $p$.

  \begin{thm}\label{thd1}
      Let  $F$ be a field of characteristic $p>0$, and $v: F\xrightarrow{} \Z\cup{\{\infty\}}$ a valuation {for which} $F$ is complete. Let $\overline{F}$ and  $\pi$  be the residue field  and a uniformizer  of $F$, respectively.  If $\overline{F}\subset F$, then    any form $\phi\in H_p^{m+1}( F)$ can be expressed as  $$  \phi= \phi_1+ \psi+ \phi_2\wedge \frac{\dd \pi}{\pi}$$ where $\phi_1\in H_p^{m+1}(\overline{F})$, $\phi_2\in H_p^m(\overline{F})$ and $\psi= \sum\limits_{I\in (T)_m}\sum\limits_{J\in T\mid l(J)\notin \tilde{I} } t^J r^p_{I,J} \frac{\dd t_I}{ t_I}$, for $r_{I,J}\in \pi^{-1}\overline{F}[\pi^{-1}]$.  Moreover,  $\phi_1$, $\phi_2$ and $r_{I,J}$ are unique for each $I$, $J$. As instance there is a split exact sequence:$$  0\xrightarrow{} H_p^{m+1}(\overline{F})\xrightarrow{}  H_p^{m+1}(F)\xrightarrow{}  \R\oplus H_p^m(\overline{F})\wedge \frac{\dd \pi}{\pi}\xrightarrow{}0$$ 
  \end{thm}
\begin{proof}
    Similar to \cite[Theorem 2.5]{ex}.
\end{proof} 
  \begin{rem}\label{rex} Let  $v:F\xrightarrow{}\Z\cup\{\n\}$ be a valuation such that $F$ is complete and let $\overline{F}$ be the residue field of $F$. By Theorem \ref{thd1}, we have the following split exact sequence: $$  0\xrightarrow{} H_p^{m+1}(\overline{F})\xrightarrow{i_1} H_p^{m+1}(F)\xrightarrow{\partial_{v}} W_1(H_p^{m+1}(F))\xrightarrow{} 0$$  
      where $i_1$ is induced from a Teichm\"uller lifting and  $\partial_{v}$ is the  second residue map.  
  \end{rem}

  \section{The filtration of $ H_p^{m+1}( \F(x))$}\label{5sts}
  
For the rest of this {note},   let $\F$ be a field of characteristic $p$  with $\{ t_i\mid i\in \I\}$ as a $p$-basis, where $\I$ is an  ordered set, and {let  $F=\F(x)$ be the rational function field in one variable $x$ over $\F$}. Let $F_q$ and $\overline{F}_q$ be the completion and the  residue field of $F$, respectively, for the $q$-adic valuation, where $q=1/x$ or  $q\in \F[x]$ is a monic  irreducible polynomial. Now we make a filtration $(L_d)_{d\in \N}$ of  $H_p^{m+1}(\F(x))$ as was done {in \cite[Section 3]{ex}}   for $p=2$.   For any positive integer $d$,  we define the sets  $\F[x]_{<d}:= \{ f\in \F[x]\mid \deg f<d\}$  and $\F[x]_{\leq d}:= \{ f\in \F[x]\mid \deg f\leq d\}$.

  \begin{defn} (1) For $d\in \N_0$,  let $L_d$ be the subgroup of $H_p^{m+1}(\F(x))$ generated by the forms $$  \frac{h}{u^e} \frac{\dd f_1}{f_1}\wedge\ldots\wedge \frac{\dd f_m}{f_m}$$ where $ h\in \F[x]$, $e\geq 0$ and $u$, $f_i\in \F[x]_{\leq d}\setminus \{0\}$.\\

  \noindent{(2)} Let $L_0$ be the subgroup of  $ H_p^{m+1}(\F(x))$ generated by the forms $$  h \frac{\dd c_1}{c_1}\wedge \ldots\wedge \frac{\dd c_{m-1}}{c_{m-1}}\wedge \frac{\dd f_m}{f_m}$$ where $ c_i\in \F^*$ and $ f_m\in \F^*\cup\{x\}$. Moreover,   $h\in \F[x]$  (resp. $h\in x\F[x]$)  if $f_m\in F^*$ (resp. $ f_m=x$). 
  \end{defn}

  Clearly $(L_d)_{d\in \N}$ is an increasing  chain of subgroups of $H_p^{m+1}(\F(x))$ {satisfying} ${ H_p^{m+1}(\F(x))=\cup_{d\in \N} L_d}$.   {To prove  the subsequent results, we need the following relations:  \begin{align}\label{de}\frac{\dd(a+b)}{a+b}&= \frac{a}{a+b}\frac{\dd a}{a}+\frac{b}{a+b}\frac{\dd b}{b}\in \Omega_F^1 \\ \label{ef1}\frac{\dd a^r}{a^r}&= r\frac{\dd a}{a}\end{align} 
where $r\in \Z$ and $a$, $b\in F^*$ satisfying $a+b\neq 0$.}
  The following  lemma proves that    the second residue map $ \partial_q$ can be extended to $L_d/L_{d-1}$, for any monic  irreducible polynomial $q\in \F[x]$  of degree $d$.  

   \begin{lemma}\label{lvalu}
       Let $q\in \F[x]$ be a monic irreducible polynomial of degree $d$, and  $v_q$ the $q$-adic valuation  of $\F(x)$. Let $ s$, $r_1,\ldots,r_m\in \F(x)$. 
       \vspace{1mm}
       
       \noindent{(1)} If $v_q(s)\geq 0$ and $v_q(r_i)=0$ for all $i=1,\ldots,m$, then  $$  \partial_q\left( s\frac{\dd r_1}{r_1}\wedge \ldots\wedge \frac{\dd r_m}{r_m}\right)=0 $$

       \noindent{(2)}  For any $\phi\in L_{d-1}$, we have   $\partial_q(\phi)=0$. Therefore, $\partial_q$ extends to $L_d/ {L_{d-1}}$.
   \end{lemma}
   \begin{proof}
       It is similar to \cite[Lemma 3.4]{ex}.
   \end{proof}

   \begin{defn}\label{definse}
        For any  monic  irreducible polynomial $q\in \F[x]$, let $S_q$ be the subgroup of $H_p^{m+1}(\F(x))$ generated by the elements $$  \frac{h}{q^e} \frac{\dd c_1}{c_1}\wedge\ldots\wedge \frac{\dd c_m}{c_m},\quad \text{ where } h\in \F[x], \hspace{1mm} e\geq 0  \text{ and }c_i\in \F^*.$$
   \end{defn}
   The following lemma describes some generators of  the group $L_d$ for $d\geq1$.
   \begin{lemma}\label{lgener}(1) For any positive integer $d$, we have $L_d= \sum_{ \deg q\leq d}\left(S_q+S_q\wedge \frac{\dd x}{x}\right)$, where $S_q\wedge \frac{\dd x}{x}$ is the subgroup of $H_p^{m+1}(\F(x))$ generated by the elements: $$  \frac{h}{q^e} \frac{\dd c_1}{c_1}\wedge \ldots\wedge \frac{\dd c_{m-1}}{c_{m-1}}\wedge \frac{\dd x}{x}$$ where $ h\in \F[x]$, $c_i\in \F$ and $e\geq 0$.

       \noindent{(2)} $H_p^{m+1}(\F(x))= \sum_{q}\left( S_q+ S_q \wedge \frac{\dd x}{x}\right)$, where $q$ varies over  all monic   irreducible polynomials over $\F$.
       \vspace{1mm}

       \noindent{(3)}  For each $d\geq 1$, we have $ L_d= \sum_{ \deg q=d}\left( S_q+S_q\wedge \frac{\dd x}{x}\right)+L_{d-1}$.
   \end{lemma}

   \begin{proof}
       It is similar to \cite[Lemma 3.6]{ex}.
   \end{proof}

\textbf{Notation:} Let $q= \sum_{j=0}^{d} q_j x^{d-j}\in \F[x]$ be a monic irreducible polynomial of degree $d$.  Since $ \{ t_i\mid i\in \I\}$ is a $p$-basis of $\F$, we have $q_j= \sum_{I\in T}t^I q_{j,I}^p$ for some $q_{j,I}\in \F$, where $T= \oplus_{|\I|}\Z/p\Z$. \\ Let \label{pagf}$\I':= \{ i\in \I\mid  \exists I\in T \text{ such that } i\in \tilde{I}\text{ and } q_{j,I}\neq0\text{ for some }j\}$.    Let us define the subset $T_1:= \oplus_{|\I'|} \Z/p\Z$ of $T$. If $\I'=\emptyset$, then we take $T_1:=\{0\}\subset T$. If $\I'\neq \emptyset $,  then we take $\I':=\{j_1,\ldots,j_n\}\subset \I$ with $j_1<j_2<\ldots<j_n$.   Thus, for all $j\in \{0,\ldots,d\}$, we have: \begin{equation}\label{ee}q_j= \sum_{K\in T_1}t^K q_{j,K}^p \quad \text{ for some  $q_{j,K}\in \F$}\end{equation}  where $$t^K=\begin{cases}t_{j_1}^{K_1}\ldots t_{j_n}^{K_n}\quad \text{if $\I'\neq \emptyset$}\\ 1\hspace{18mm} \text{if $ \I'=\emptyset$.} \end{cases}$$ Now $T_1$ is {lexicographically ordered}  like $T$.  For all $K\in T_1$,  we define \begin{equation}\label{edq}   {Q_K = \sum_{j=0}^{d} (d-j)t^{K} q^p_{j,K} x^{d-j} \in\F[x]}.\end{equation} Let  \label{dm} $M$ be the matrix  over $\overline{F}_q$  whose $(I,J)$-th entry is {$\overline{Q}_{I-J}$} for all $I,J\in T_1$. We prove that $M$ is invertible whenever $q$ is separable.

\begin{lemma}\label{linv}
    Let $T$ be a commutative  group with $ p^n$ elements  for some $n\in \N$, and let $Q_K$ the elements of a field of characteristic $p$, for $K\in T$. Let $M$ be the matrix of order $p^n$ whose $(I,J)$-th element is $Q_{I-J}$ for all $I$, $J\in T$.  If $\sum_{K\in T} Q_K\neq 0$, then $M$ {is invertible.}  
\end{lemma}

\begin{proof} {The lemma was proved in \cite[Lemma 2.6]{AJ} for $p=2$}. Let us suppose $p> 2$. If $n=0$, then the lemma is trivial. Let $n\in \N_0$ and  $ \mathcal{T}:=\{ \sigma: T\xrightarrow{} T \text{ is a bijection}\}$. Then,  we have $$  \det M= \sum_{\sigma\in \mathcal{T}}{\rm sgn }(\sigma) \prod_{I\in T} Q_{I-\sigma(I)}.$$ Let $K$ be an arbitrary element of $T$. Let     $f_K: T\xrightarrow{} T$ be such that $ f_K( I)= I-K$ for all $I\in T$. Clearly $f_K\in \mathcal{T}$ for all $K\in T$. Moreover, \begin{align} \label{5b1}\prod_{I\in T} Q_{I-f_K(I)}=(Q_K)^{p^n}.\end{align}
    We define an  action  $\Phi: T\times  \mathcal{T}\xrightarrow{} \mathcal{T}$, such that $\Phi (J, \sigma)= \sigma^J$, where $ \sigma^J:T\xrightarrow{} T$ is a map defined by $$\sigma^J( I) = \sigma( I-J) +J,\quad \text{ for all $ I\in T$.}$$ Since $T$ is commutative, $\Phi$ is an action.  Now  for any $\sigma\in \mathcal{T}$ and  $ J\in T $, we have \begin{align}\label{5b2}
        \prod_{I\in T} Q_{I-\sigma^J(I)}=  \prod_{I\in T} Q_{I- \sigma(I-J)-J}= \prod_{I\in T} Q_{I-J- \sigma(I-J)}= \prod_{L\in T} Q_{L-\sigma(L)}.
    \end{align}  Let $\sigma$ be an  element of $\mathcal{T}$ such that  $\sigma^J=\sigma$ for all $J\in T$. Then  
    
    {\begin{eqnarray*}
        \sigma^J(J) =\sigma(J),\,\text{ for all } J\in T & \Longleftrightarrow & \sigma(0)+J = \sigma(J),\,\text{ for all } J\in T\\ &\Longleftrightarrow &  \sigma = f_{-\sigma(0)}\in \mathcal{T}.
    \end{eqnarray*}} Let $\sigma\in \mathcal{T}$ be such that $\sigma\neq f_K$ for each $K\in T$. Then, the orbit of $\sigma$ has $p^r$ elements for some $r\geq 1$. Moreover,  $\sigma^J= {f_J}^{-1}\circ \sigma\circ f_J$ and thus ${\rm sgn}(\sigma^J)= {\rm sgn}( \sigma)$, for all $J\in T$. Hence, \begin{align*}
         \sum_{\phi\in {\rm orbit}(\sigma)}{\rm sgn}(\phi) \prod_{I\in T} Q_{ I-\phi(I)}&=  p^r {\rm sgn}(\sigma) \prod_{I\in T} Q_{I-\sigma(I)}\quad (\text{by (\ref{5b2})})\\&=0. 
    \end{align*} Since   $ f_K\circ f_0\circ f_K= f_{2K}$ for all $K\in T$, we have  ${\rm sgn} (f_0)= {\rm sgn}(f_{2K})$ for all $K\in T$. Since   the cardinality of $T$ is {coprime} to $2$,  the following homomorphism is bijective: $$ \Psi: T\xrightarrow{}T,\quad \Psi(K)=2K, \quad \text{for all } K\in T.$$ Thus, $ {\rm sgn} (f_0)= {\rm sgn} (f_K)$ for all $K\in T$. Moreover, the orbit of $ f_K$ for the action $\Phi$ is $\{f_K\}$ for all $K\in T$. Therefore, we get \begin{align*}
         \det M&= {\rm sgn}(f_0) \sum_{K\in T}  (Q_K)^{p^n}\quad (\text{by (\ref{5b1})})\\&= {\rm sgn }(f_0) \left( \sum_{K\in T}Q_K\right)^{p^n}.
    \end{align*}Thus, $M$ is invertible if and only if $ \sum_{K\in T} Q_K\neq 0$.
\end{proof}

\begin{cor}\label{corsep}
    Let $q\neq x$ be a monic irreducible  polynomial of degree $d$ over $\F$.    If $q$ is separable, then  the matrix $  M$ as defined on page \pageref{dm} is invertible over $\overline{F}_q$.
\end{cor}

\begin{proof}
    Suppose $M$ is not invertible.  By Lemma \ref{linv},  we have   $ \sum_{K\in T_1} {\overline{Q}_K}=0\in \overline{F}_q$.  Since $\overline{x}$ is a unit in $\overline{F}_q$, we have \begin{align*} \sum_{K\in T_1}  \sum_{j=0}^{d-1} \overline{ (d-j) t^K q^p_{j,k} x^{d-j-1}}&=0\in \overline{F}_q\\\sum_{K\in T_1}  \sum_{j=0}^{d-1} { (d-j) t^K q^p_{j,k} x^{d-j-1}}&=0 \in \F[x]\quad \text{(since $\deg q=d$)}\end{align*}{Therefore,  for all $j\in \{0,\ldots,d-1\}$ such that  $p\nmid d-j$, we get $\sum_{K\in T_1} t^K q^p_{j,K}=0\in \F$. Hence, for all $j\in \{0,\ldots,d-1\}$ such that  $p\nmid d-j$, we have $q_{j,K}=0$ because $\{t_i\mid i\in \I\}$ is a $p$-basis of $\F$.} Now it follows from  (\ref{ee}) that  $q_j=0$ for all $j\in \{0,\ldots,d-1\}$ such that $p\nmid d-j$. It contradicts the separability of $q$. Thus, $M$ is invertible over $\overline{F}_q$.\end{proof}
    

\begin{lemma}\label{lepbas} Let $\F$ be a field with $\{t_i\mid  i\in \I\}$  as a $p$-basis, and  
       $ q$ a monic irreducible  inseparable polynomial over $\F$.  Then, there exists {an element} $l\in \I$ such that $ \{ \overline{t}_i\mid i\neq l\in \I\}\cup \{\overline{x}\}$ is a $p$-basis of $\overline{F}_q$ and $ \overline{\frac{\partial q}{\partial t_l}}\neq 0\in \overline{F}_q$. 
\end{lemma}
\begin{proof}  
{We put} $q= \sum_{j=0}^{d} q_jx^{d-j}\in \F[x]$, where $q_0=1$.  Let $\I'$ and $T_1$ be as defined in page \pageref{pagf}.  Since $q$ is inseparable, {it follows that $p$ divides $d$, and}  $q_j=0$ for all $j\in \{0,\ldots,d\}$ such that   $p\nmid j$.  Thus, $\I'\neq \emptyset$, otherwise  $q$ will be  reducible. By (\ref{ee}), for all $j\in \{0,\ldots,d\}$, we have  $$   q_{j}=\sum_{K\in T_1} t^K q^p_{j,K}\quad \text{for some $q_{j,K}\in \F$.}$$   Let $K$ be the smallest subfield of $\F$ containing $\{t_i\mid i\in \I'\}$ and $\{q_{j,K}\mid j\in \{0,\ldots, d\},\; K\in T_1\}$.  Then, the $p$-form\footnote[6]{We refer to \cite{Hoffmann1} for  the definition of a $p$-form.} $\otimes_{i\in \I'} \langle\langle \overline{t}_i\rangle \rangle$ is isotropic over $K(q):=K[x]/(q)$. By \cite[Corollary 4.10]{Hoffmann1},  there exists some $l\in \I'$ such that $ \overline{t}_l\in K(q)^p( \overline{t}_i\mid i\in\I'\setminus\{l\})$. Since $K(q)\subset \overline{F}_q$,  the set $S:=\{ \overline{t}_i\mid i\in \I\setminus\{l\}\}\cup \{\overline{x}\}$  generates $\overline{F}_q$ over $\overline{F}_q^p$. Since $l\in \I'$ and $q_j=0$ for all $j$ such that $p\nmid d-j$, we have $\overline{\frac{\partial q}{\partial t_l}}= \sum_{j=1}^{d}\overline{\frac{\partial q_j}{\partial t_l} x^{d-j}}\neq 0\in \overline{F}_q$. Now it suffices to prove that  $S$ is $p$-independent over $\overline{F}_q$.  To this end, we {distinguish between two cases where} $\I$ is finite or not.\\

\noindent{\textbf{Case 1.}}  Suppose $\I$ is finite.
    Now we write $q= f(x^{p^n})$, where $n\geq 1$ and  $f$ is a monic separable irreducible polynomial over $\F$.   Let $\alpha$ be a root of $q$ {in the algebraic closure of $\F$}, and $a= \alpha^{p^n}$. Then $\overline{F}_q\cong \F(\alpha)= K'(\alpha)$, where $K'=\F(a)$.  Since $K'/{\F}$ is separable, it follows from \cite[Lemma 2.1]{AAB}  that $\{ t_i\mid i\in \I\}$ is a $p$-basis of $K'$.  Then, $a= \sum_{I\in T} t^I a_I^p$ for unique $a_I\in K'$, where $T= \oplus_{|\I|}\Z/p\Z $. Let $B:=\{ i\in \I\mid \exists I\in T\text{ such that }i\in \tilde{I}\text{ and } a_I\neq 0\}$. Since $\sqrt[p]{a}\notin {K'}$, we have $B\neq \emptyset$. Let $k$ be an arbitrary element of $B$.  
     Then, the set   $\{ t_i\mid {i\in \I\setminus\{k\}}\} \cup \{ a\}$ is  $p$-independent over $K'$.  Otherwise,  {there exists $C$ a subset of $\I\setminus\{k\}$ such that} the $p$-form $\otimes_{i\in C} \langle\langle t_i\rangle\rangle\otimes \langle\langle a\rangle\rangle $ is isotropic over $K'$. If $C=\emptyset$, then the $p$-form $\langle\langle a\rangle\rangle$  is isotropic over $K'$. Then by  \cite[Corollary 4.10]{Hoffmann1}, we have $\sqrt[p]{a}\in K'$ which is a contradiction. Therefore, $C\neq \emptyset$. Then by Lemma \ref{lan},  $ a\in {K'}^p( t_i\mid {i\in \I\setminus\{k\}})$ which contradicts  $k\in B$.  Since $\I$ is a  finite set,  $\{ t_i\mid i\in \I\setminus\{k\}\} \cup \{a\}$ is a $p$-basis of $K'$ {because it has the same cardinality as the $p$-basis $\{ t_i\mid i\in \I\}$ of $K'$}. Therefore, $ \{ t_i\mid  {i\in \I\setminus\{k\}}\}\cup \{ \alpha\}$ is a $p$-basis of $K'(\alpha)=\F(\alpha)$. Since $\I$ is  finite and $S=\{ \overline{t}_i\mid {i\in \I\setminus\{l\}}\}\cup\{ \overline{x}\}$ is a $p$-generating set of $\overline{F}_q$, it follows that $S$ is a $p$-basis of $\overline{F}_q$. \\
     
\noindent{\textbf{Case 2.}} {Suppose that $\I$ is infinite}.   For any finite subset $\mathcal{J}$ of $\I$, we define $ T_{\mathcal{J}}:=\oplus_{|\mathcal{J}|} \Z/p\Z$. For each $I\in T_{\mathcal{J}}$, {let $t^I$ denote} the element $\prod_{i\in {\mathcal{J}}}t_i^{I_i} $, where $I_i$ are the components of $I$. {Suppose that $S$ is} $p$-dependent over $\overline{F}_q$. Then, there exists  a finite subset $M\subset \I\setminus \{l\}$ such that $\{\overline{t}_i\mid i\in M'\}$ is $p$-dependent over $\overline{F}_q$, where $M'=M\cup \{i_x\}$ and  $t_{i_x}=\overline{x}$. Then, there exist $ f_J\in \F[x]_{<d}$  not all zero polynomial such that \begin{equation}\label{eex} \sum_{J\in T_{M'}}\overline{f_J^pt^J}=0\in \overline{F}_q.\end{equation} Let $ G$ be a  subfield of $\F$ containing  $K$ and  $\{t_i\mid i\in M\}$ such that $G$ has a finite $p$-basis and $f_J\in G[x]$ for all $J\in T_{M'}$.  Since $\{t_i\mid i\in M\cup \I'\}$ is $p$-independent over $G$, it  extends to a  finite $p$-basis  (say) $ B$  of $G$.    Since $K(q)\subset G(q):= G[x]/(q)$, we have     $\overline{t}_l\in G(q)^p(\overline{t}_j\mid {j\in \I'\setminus\{l\}})$. Since $G$ has a finite $p$-basis, it follows from  {the first case} that $\{\overline{t}\mid {t\in B\setminus\{t_l\}}\}\cup \{\overline{x}\}$  is a  $p$-basis of $G(q)$. Therefore,  $\{ \overline{t}_i\mid i\in M'\}$ is $p$-independent over $G(q)$. Then, it follows from (\ref{eex}) that ${f}_J=0\in G[x]\subset \F[x]$ for all $J\in T_{M'}$ which  contradicts our assumption. Therefore,  $S:=\{\overline{t}_i\mid i\neq l\in \I\}\cup\{\overline{x}\}$ is a $p$-basis of $\overline{F}_q$. Thus, the lemma follows.
\end{proof}

\begin{defn}\label{dein}
    If $q$ is inseparable, then by  Lemma \ref{lepbas}  there exists some $l\in \I'\subset \I$ such that $ \{ \overline{t}_i\mid {i\in \I\setminus\{l\}} \}\cup \{\overline{x}\}$ is a $p$-basis of $\overline{F}_q$. {Without}  loss of generality, we take $ l=j_n$. Now we define a subgroup $\tilde{S}_q\subset  H_p^{m+1}(\F(x))$ generated by the forms $$ \displaystyle{ \frac{h}{q^s} \frac{\dd t_{i_1}}{t_{i_1}}\wedge \ldots\wedge \frac{\dd t_{i_m}}{t_{i_m}}}$$ where $ \{ t_{i_1}, \ldots, t_{i_m}\}\subset \{ t_i\mid i\in (\I\setminus \{j_n\})\cup \{i_x\}\}$ with   $i<i_x$ for all $i\in \I\setminus\{j_n\}$ and $t_{i_x}=x$. Moreover, $h\in \F[x]$ and $s\geq 0$.
\end{defn}
The following proposition  describes the relations among $ S_q\wedge \frac{\dd q}{q}$, $\tilde{S}_q\wedge \frac{\dd q}{q}$ and ${S}_q\wedge \frac{\dd x}{x}$.
\begin{prop}\label{prsei} Let $q$ be a  monic irreducible polynomial over $\F$. Then, we have:\\

\noindent{(1)}   $\displaystyle{S_q+ S_q\wedge \frac{\dd q}{q} \subset S_q+S_q\wedge \frac{\dd x}{x}}$.\\

\noindent{(2)}  If $q$ is separable, then $ \displaystyle{S_q+ S_q\wedge \frac{\dd x}{x}+ L_0\wedge \frac{\dd x}{x}= S_q+ S_q\wedge \frac{\dd q}{q}+ L_0\wedge \frac{\dd x}{x}}$.\\

\noindent{(3)} If $q$ is inseparable, then $ \displaystyle{S_q+ S_q\wedge \frac{\dd x}{x}+L_0+ L_0\wedge \frac{\dd x}{x}= \tilde{S}_q+\tilde{S}_q\wedge \frac{\dd q}{q}+L_0+L_0\wedge \frac{\dd x}{x}}$.
    \end{prop}

    \begin{proof}
        (1)    Let  $q= \sum_{i=0}^{d} q_i x^{d-i}\in \F[x]$, where $q_0=1$.     Now the additive generators of $ {S_q\wedge \frac{\dd q}{q}}$ are  $ \frac{h}{q^s} \frac{\dd c_1}{c_1}\wedge \ldots\wedge \frac{\dd c_{m-1}}{c_{m-1}}\wedge \frac{\dd q}{q}$, where $ h\in \F[x]$, $s\geq 0$ and $c_i\in \F$. By  relation (\ref{de}),  we have: \begin{align*}
            \frac{h}{q^s}\frac{\dd c_1}{c_1}\wedge \ldots\wedge \frac{\dd c_{m-1}}{c_{m-1}}\wedge \frac{\dd q}{q}&= \sum_{j=0}^{d}\frac{ q_jhx^{d-j}}{q^{s+1}}\frac{\dd c_1}{c_1}\wedge\ldots\wedge \frac{\dd c_{m-1}}{c_{m-1}}\wedge \frac{\dd (q_jx^{d-j})}{q_jx^{d-j}}\in H_p^{m+1}(\F(x))\end{align*}
           Using relation (\ref{ef1}), we get: 
            \begin{align*} \frac{h}{q^s}\frac{\dd c_1}{c_1}\wedge \ldots\wedge \frac{\dd c_{m-1}}{c_{m-1}}\wedge \frac{\dd q}{q}
            &\equiv   \sum_{j=0}^{d}\frac{(d-j) q_jhx^{d-j} }{q^{s+1}} \frac{\dd c_1}{c_1}\wedge \ldots\wedge \frac{\dd c_{m-1}}{c_{m-1}}\wedge \frac{\dd (q_j x)}{q_j x}\pmod{S_q}\end{align*}
        { Since $ q_j=\sum_{K\in T_1} t^K q^p_{j,K} $, we have:}\begin{align*} \frac{h}{q^s}\frac{\dd c_1}{c_1}\wedge \ldots\wedge \frac{\dd c_{m-1}}{c_{m-1}}\wedge \frac{\dd q}{q}&\equiv  \sum_{K\in T_1} \sum_{j=0}^{d}\frac{ (d-j) t^Kq^p_{j,K}h x^{d-j}}{ q^{s+1}} \frac{\dd c_1}{c_1}\wedge \ldots\wedge \frac{\dd c_{m-1}}{c_{m-1}}\wedge \frac{\dd (t^K x)}{t^Kx}\pmod{S_q}\end{align*}Using the expression of $Q_K$  of (\ref{edq}), we get: \begin{align}\label{5e1}
            \frac{h}{q^s}\frac{\dd c_1}{c_1}\wedge \ldots\wedge \frac{\dd c_{m-1}}{c_{m-1}}\wedge \frac{\dd q}{q}&\equiv \sum_{K\in T_1} \frac{Q_Kh}{q^{s+1}} \frac{\dd c_1}{c_1}\wedge \ldots\wedge \frac{\dd c_{m-1}}{c_{m-1}}\wedge \frac{\dd (t^K x)}{t^Kx}\pmod{S_q} 
        \end{align} Since $t^K\in \F$,   it follows from (\ref{5e1}) that  $\displaystyle{ S_q\wedge \frac{\dd q}{q}\subset S_q+ S_q\wedge \frac{\dd x}{x}}$. Thus statement (1) is proved.\\

        \noindent{(2)}  Similar to \cite[Proposition 3.9, (2)]{ex}. We  have to use    relation (\ref{5e1}).
        
\noindent{(3)} Let $q= \sum_{i=0}^{d}q_i x^{d-i}\in \F[x]$ be  monic irreducible  inseparable of degree $d$. Then,  $p\mid d$, $q_0=1$ and $q_i=0$ for all $i\in \{1,\ldots,d\}$ such that {$p\nmid d-i$}.   Now we write  $\displaystyle{q_i=\sum_{s=0}^{p-1} q_{i,s} t_{j_n}^s } $ for some $q_{i,s}\in \F^p( t_k\mid {k\in \I'\setminus\{j_n\}})$ for all $i\in \{0,\ldots,d\}$. Therefore$$ q_{i,s}= \sum_{K\in T_1,\, K_{j_n}=0} q_{i,s,K}^p t^K,$$ where $q_{i,s,K}\in \F$ for all $i$ and $s$. Now for each $K\in T_1$ with $K_{j_n}=0$, we define \begin{align}\label{sd2} {\mathcal{Q}_K= \sum_{i=0}^{d}\sum_{s=1}^{p-1} {  s t_{j_n}^{s-1} q^p_{i,s,K}t^K x^{d-i}}\in \F[x]}\end{align} Taking sum on $K\in T_1$ with $K_{j_n}=0$, we get \begin{align*}
     { \sum_{K\in T_1,\,K_{j_n}=0}\overline{\mathcal{Q}}_K}&=\sum_{K\in T_1,\, K_{j_n}=0}\sum_{i=0}^{d}\sum_{s=1}^{p-1} \overline{  s t_{j_n}^{s-1} q^p_{i,s,K} t^Kx^{d-i}}= \sum_{i=0}^{d} \overline{\frac{\partial q_i}{\partial t_{j_n}} x^{d-i}}= \overline{ \frac{\partial q}{\partial t_{j_n}}}\neq 0\in \overline{F}_q.
 \end{align*} Let $M$ be the matrix  of order $p^{n-1}$ with $(I,J)$-th element as {$\overline{\mathcal{Q}}_{I-J}\overline{t}_{j_n}$} for all $I,J\in T_1$ with $I_{j_n}=J_{j_n}=0$. Then, by Lemma \ref{linv}, $M$ is invertible over $\overline{F}_q$.  
 
  Now we have to repeat the proof of \cite[Proposition 3.9, (3)]{ex}. At first we do some calculations.
  
 \noindent{\textbf{Calculations:}} Let $\frac{\dd t_I}{t_I}$ denote the element $ {\frac{\dd t_{i_1}}{t_{i_1}}\wedge \ldots\wedge \frac{\dd t_{i_{m-1}}}{t_{i_{m-1}}}}$, where $\tilde{I}=\{i_1,\ldots,i_{m-1}\}\subset \I\setminus \{j_n\}$ with $i_1<\ldots<i_{m-1}$.  Let    $e\geq 0$ and  $h_L\in \F[x]$  for   $L\in T_1$ with $L_{j_n}=0$.  By (\ref{de}), we have:\begin{align*}
     \frac{h_L}{q^e}\frac{\dd t_I}{t_I}\wedge \frac{\dd (t^Lq)}{t^Lq}& = \sum_{i=0}^{d}\frac{h_Lq_ix^{d-i}}{q^{e+1}} \frac{\dd t_I}{t_I}\wedge \frac{\dd (t^Lq_i)}{t^Lq_i}\quad (\text{since $ q_i=0$ for $p\nmid d-i$})\\&= \sum_{i=0}^{d} \sum_{s=0}^{p-1} \frac{ h_L q_{i,s} t^s_{j_n} x^{d-i}}{q^{e+1}}\frac{\dd t_I}{t_I}\wedge \frac{\dd (t^L q_{i,s} t_{j_n}^s)}{t^L q_{i,s} t^s_{j_n}} \hspace{2mm} \text{(using $q_i=\sum_{s=0}^{p-1}q_{i,s}t^s_{j_n}$)}\\&\equiv \sum_{i=0}^{d} \sum_{s=1}^{p-1} \frac{ h_L q_{i,s} t^s_{j_n} x^{d-i}}{q^{e+1}}\frac{\dd t_I}{t_I}\wedge \frac{\dd (t^L q_{i,s} t_{j_n}^s)}{t^L q_{i,s} t^s_{j_n}}\pmod{\tilde{S}_q}\\ &\equiv \sum_{i=0}^{d} \sum_{s=1}^{p-1} \frac{s h_L q_{i,s} t_{j_n}^s x^{d-i}}{q^{e+1}}\frac{\dd t_I}{t_I}\wedge \frac{\dd (t^L q_{i,s} t_{j_n})}{t^L q_{i,s}t_{j_n}}\pmod{\tilde{S}_q}\quad \text{(using (\ref{ef1}))}\intertext{ Since $ q_{i,s}=\sum\limits_{K\in T_1,\, K_{j_n}=0} q_{i,s,K}^p t^K $, we have: }  \frac{h_L}{q^e}\frac{\dd t_I}{t_I}\wedge \frac{\dd (t^Lq)}{t^Lq}& \equiv  \sum_{i=0}^{d} \sum_{s=1}^{p-1}\sum_{K\in T_1,\, K_{j_n}=0} \frac{s h_L q_{i,s,K}^p t^K t^s_{j_n} x^{d-i}}{q^{e+1}}\frac{\dd t_I}{t_I}\wedge \frac{\dd (t^{K+L} t_{j_n})}{t^{K+L}  t_{j_n}}\pmod{\tilde{S}_q}\\&= \sum_{K\in T_1,\, K_{j_n}=0} \frac{ {\mathcal{Q}_K h_Lt_{j_n}}}{q^{e+1}}\frac{\dd t_I}{t_I}\wedge \frac{\dd (t^{K+L} t_{j_n})}{t^{K+L} t_{j_n}}\quad \text{(by (\ref{sd2}))}
 \end{align*} Now repeating the proof of \cite[Proposition 3.9, (3)]{ex}, we  can prove statement (3).
 \end{proof}

 \section{Milnor's exact sequence}\label{5s4}

 If $q=1/x$ or $q$ is a monic irreducible separable polynomial over $\F$, then it follows from  \cite[Lemma 2.1]{AAB} that the set   $\{ \overline{t}_i\mid i\in \I\}$ is a $p$-basis of $\overline{F}_q$, and  $ \{ t_i\mid i\in \I\} \cup \{q\}$ is a $p$-basis of $F_q$.  If $q$ is inseparable, {then with the same notations as in Definition  \ref{dein}, the set}  $ \{ \overline{t}_i\mid i\neq j_n\in \I\}\cup\{\overline{x}\}$ is a $p$-basis of $\overline{F}_q$. Consequently, $ \{ t_i\mid i\neq j_n\in \I\}\cup\{ x, q\}$ is a $p$-basis of $F_q$ for $q$ inseparable. Now we fix some notations related to these $p$-bases of $\F$, $\overline{F}_q$ and $F_q$.

\begin{defn}\label{dpb}(1)  Let  $T:= \oplus_{|\I|}\Z/p\Z $ be such that $\{ t^I\mid I\in T\}$ is  a $\F^p$-basis of $\F$, which comes from the fixed $p$-basis $ \{ t_i\mid i\in \I\}$ of $\F$.\\

{\noindent{(2)} Let $ \overline{T}_q:= \oplus_{|\I|}\Z/p\Z$ be such that  $ \{ t^I\mid I\in \overline{T}_q\}$ is a $\overline{F}_q^p$-basis of $\overline{F}_q$, which comes from the $p$-basis $\{ \overline{t}_i\mid i\in \I\}$ (\resp $ \{ \overline{t}_i\mid i\neq j_n\in \I\}\cup\{\overline{x}\}$) if $q$ is separable or $q=1/x$(\resp $q$ is inseparable).}\\

{\noindent{(3)} Let $T_q:= \oplus_{|\I|+1}\Z/p\Z$ be such that  $ \{ t^I\mid I\in T_q\}$  is a $F_q^p$-base of $F_q$, which comes from the $p$-basis $\{ t_i\mid i\in \I\} \cup \{q\}$ (\resp $ \{ t_i\mid i\neq j_n\in \I\}\cup\{ x, q\}$) if $q$ is separable or $q=1/x$ (\resp $q$ is inseparable).}
\end{defn}
 Now we describe the group  $ W_1(H_p^{m+1}(F_q))$ for $q\in \F[x]$ monic irreducible or $q=1/x$.

 \begin{thm}\label{thd2}
    (1) Let $q$ be a monic irreducible polynomial of degree $d$ over $\F$. Then, each element $\phi\in W_1( H_p^{m+1}(F_q))$ can be decomposed as $ \phi=\psi+ \phi'\wedge \frac{\dd q}{q} $, where $ \phi'\in H_p^{m}(\F(q))$ and $$  \psi= \sum_{I\in (T_q)_m}\sum_{J\in T_q,\, l(J)\notin \tilde{I}} t^J u^p_{I,J}\frac{\dd t_I}{t_I}$$ with $ u_{I,J}= \sum_{l\geq 1} \frac{u_{l,I,J}}{q^l}$, $u_{l,I,J}\in\F[x]_{\leq d-1} $ and $\F(q):= \F[x]/(q)$. Moreover, $u_{l,I,J}$ and $\phi'$ are unique for each $\phi$. Here $(T_q)_m=\{ I\in T_q\mid I_i\in \{0,1\},\hspace{1mm}  |\tilde{I}|=m\}$ for $m\in \N_0$.\\
    
    \noindent{(2)} If $q=1/x$, then each element $\phi\in H_p^{m+1}(F_{1/x})$ can be decomposed as $\phi= \psi- \phi'\wedge \dd x/x$, where $\phi'\in H_p^m(\F)$ and $$ \psi=\sum_{I\in (T_{\frac{1}{x}})_m}\sum_{J\in T_{\frac{1}{x}},\, l(J)\notin \tilde{I}}t^J u^p_{I,J}\frac{\dd t_I}{t_I}, \quad { u_{I,J}\in x\F[x]}$$ where $u_{I,J}$ and $\phi'$ are unique for each $\phi$.
 \end{thm}

 \begin{proof}(1)  It is similar to \cite[Theorem 4.2]{ex}. We have to use Theorem \ref{thd1}  instead of \cite[Theorem 2.5]{ex}.  To prove uniqueness, we use the residue maps $\zeta$ and $\chi_q$ of Definition \ref{defre}.\\

 \noindent{(2)}  Since $ \overline{F}_{1/x}\cong \F\subset F_{{1/x}}$ and $\dd (1/x)/(1/x)= -\dd x/x$, statement (2) is a direct consequence of Theorem \ref{thd1}.
      \end{proof}

      Using the above decomposition, we introduce the Milnor's splitting $\tau_q: W_1(H_p^{m+1}(F_q))\xrightarrow{} L_d/ L_{d-1}$. For any monic irreducible polynomial $q$ of degree $d$ and $g\in \F[x]$, we use $\overline{g} $ to denote the unique polynomial of degree $\leq d-1$ such that $ \overline{g}\equiv g\pmod{q}$.

      \begin{thm}\label{thmil}
          Let $q\in \F[x]$ be a monic and irreducible polynomial of degree $d$. By Theorem \ref{thd2}, each element $\phi\in W_1(H_p^{m+1}(F_q))$ can be decomposed uniquely as $$ \phi=\psi+ \phi'\wedge\frac{\dd q}{q},$$ where  $$  \psi= \sum_{I\in (T_q)_m}\sum_{J\in T_q,\, l(J)\notin \tilde{I}} t^Jr^p_{I,J} \frac{\dd t_I}{t_I}$$ such that $ r_{I,J}\in q^{-1} \F[x]_{< d}[q^{-1}]$ and  $\phi'\in H_p^m(\F(q))$. Moreover, $\phi'$ is unique over $ \F(q)$.  Let $$  \phi'= \sum_{i}\overline{r}_i\frac{\dd \overline{s}_{i_1}}{\overline{s}_{i_1}}\wedge\ldots\wedge\frac{\dd \overline{s}_{i_{m-1}}}{\overline{s}_{i_{m-1}}} $$ such that $ r_i$,  $s_{i_j}\in \F[x]_{\leq d-1}$. Then we have a homomorphism $\tau_q: W_1( H_p^{m+1}(F_q))\xrightarrow{}L_d/L_{d-1}$ given by $$ \tau_q ( \phi)= \psi+ \sum_{i}  r_i \frac{\dd s_{i_1}}{s_{i_1}}\wedge \ldots\wedge \frac{\dd s_{i_{m-1}}}{s_{i_{m-1}}}\wedge \frac{\dd q}{q}+L_{d-1} $$
      \end{thm}

      \begin{proof}
          It suffices to prove that the map $ \eta_q: H_p^m(\F(q))\xrightarrow{} L_d/ L_{d-1}$,  given  by $$  \eta_q\left( \sum_{i}\overline{r}_i\frac{\dd \overline{s}_{i_1}}{\overline{s}_{i_1}}\wedge\ldots\wedge\frac{\dd \overline{s}_{i_{m-1}}}{\overline{s}_{i_{m-1}}}\right)=\sum_{i}  r_i \frac{\dd s_{i_1}}{s_{i_1}}\wedge \ldots\wedge \frac{\dd s_{i_{m-1}}}{s_{i_{m-1}}}\wedge \frac{\dd q}{q}+L_{d-1} $$ is well defined, where $s_i$, $r_{i_j}\in \F[x]_{<d}$. By \cite[Lemma 4.4]{ex},    we have $$H_p^m(\F(q)) \cong \F(q)\otimes \underbrace{\F(q)^{*}\otimes\ldots\otimes \F(q)^*}_{m-1 \text{ times}}/J$$ where $J$ is the subgroup generated by the following elements:    \\
          
          \noindent{(1)} $ ( \overline{r}^p- \overline{r}) \otimes \overline{s}_1\otimes\ldots\otimes \overline{s}_{m-1}$, where $r$, $s_i\neq 0\in \F[x]_{< d}$,\\

          \noindent{(2)} $\overline{r}\otimes \overline{r}\otimes \overline{s}_2\otimes \ldots\otimes \overline{s}_{m-1}$,\\

          \noindent{(3)} $ \overline{r} \otimes \overline{s}_1\otimes \ldots\otimes \overline{s}_{m-1}$, where $\overline{s}_i= \overline{s}_j$ for some $i\neq j$.\\

          Let  $\zeta_q: \F(q)\times \underbrace{\F(q)^*\times \ldots\times \F(q)^*}_{m-1}\xrightarrow{} L_d/ L_{d-1}$ be defined by  $$  \zeta_q( \overline{r}\times \overline{s}_1\times \ldots\times \overline{s}_{m-1})= r \frac{\dd s_1}{s_1}\wedge \ldots\wedge\frac{\dd s_{m-1}}{s_{m-1}}\wedge \frac{\dd q}{q}+L_{d-1}. $$ Now   using the same {arguments as in} \cite[Lemma 4.5]{ex},  we prove that $\zeta_q$ is linear with respect to each quadrant. Hence, we induce the homomorphism $$ \zeta_q':  \F(q)\otimes \underbrace{\F(q)^*\otimes\ldots\otimes \F(q)^*}_{m-1}\xrightarrow{} L_d/L_{d-1} $$ defined on the generators $$ \zeta_q'( \overline{r}\otimes \overline{s}_1\otimes \ldots\otimes \overline{s}_{m-1})=r \frac{\dd s_1}{s_1}\wedge \ldots\wedge\frac{\dd s_{m-1}}{s_{m-1}}\wedge \frac{\dd q}{q}+L_{d-1}. $$  Clearly $\zeta_q'( x)=0$, where $x$ is a generator of $J$ of type $(2)$ or $(3)$.  For the  generators of type $(1)$, let $ \overline{r}^p=\overline{f}$, where $r$, $f\in \F[x]_{< d}$. Hence, $ r^p= f+qh$ for some $h\in \F[x]$. Therefore, we have:\begin{align*}
              \zeta_q'( (\overline{r}^p-\overline{r})\otimes \overline{s}_1\otimes \ldots\otimes \overline{s}_{m-1})&=  (f-r) \frac{\dd s_1}{s_1}\wedge \ldots\wedge \frac{\dd s_{m-1}}{s_{m-1}}\wedge \frac{\dd q}{q}+L_{d-1}\\&= ( r^p-r-qh) \frac{\dd s_1}{s_1}\wedge \ldots\wedge \frac{\dd s_{m-1}}{s_{m-1}}\wedge \frac{\dd q}{q}+L_{d-1}\\&= (-qh) \frac{\dd s_1}{s_1}\wedge \ldots\wedge \frac{\dd s_{m-1}}{s_{m-1}}\wedge \frac{\dd q}{q}+L_{d-1}\intertext{ Taking  $ A= -\frac{\dd s_1}{s_1}\wedge \ldots\wedge \frac{\dd s_{m-1}}{s_{m-1}}$   and  $q= \sum_{i=0}^{d} q_i x^{d-i}$, we get} \zeta_q'( (\overline{r}^p-\overline{r})\otimes \overline{s}_1\otimes \ldots\otimes \overline{s}_{m-1})&= \sum_{i=0}^{d} hq_ix^{d-i}A\wedge \frac{\dd (q_i x^{d-i})}{q_i x^{d-i}}+L_{d-1}\quad \text{(by (\ref{de}))}\\&= \sum_{i=0}^{d}\left( hq_ix^{d-i} A\wedge \frac{\dd q_i}{q_i}\right) +\sum_{i=0}^{d-1}\left( h(d-i)q_i x^{d-i}A\wedge \frac{\dd x}{x}\right)+L_{d-1}\\&=0
          \end{align*}  Thus $\zeta_q'(J)=0$. Therefore, $\zeta_q'$ can be extended on $ \F(q)\otimes \underbrace{\F(q)^{*}\otimes\ldots\otimes \F(q)^*}_{m-1 \text{ times}}/J$. Since $ H_p^{m}(\F(q))\cong \F(q)\otimes \underbrace{\F(q)^{*}\otimes\ldots\otimes \F(q)^*}_{m-1 \text{ times}}/J$,   the map $\eta_q$ is well defined. Hence, $\tau_q$ is well defined and it is a group homomorphism.  Thus,  the theorem is proved.
      \end{proof}

      \begin{lemma}\label{lim}Let $q= \sum_{i=0}^{d} q_i x^{d-i}$ be  a monic irreducible  polynomial of degree $d$ over $\F$.\\
      
           \noindent{(1)}  If $q$ is separable, then $ \displaystyle{S_q+S_q\wedge \frac{\dd q}{ q}}+L_{d-1}\subset {\rm Im} ( \tau_q) $. \\

           \noindent{(2)} If $q$ is inseparable, then $  \displaystyle{\tilde{S}_q+ \tilde{S}_q\wedge \frac{\dd q}{q}+ L_{d-1}\subset {\rm Im}( \tau_q)}$.
      \end{lemma}

      \begin{proof}
           It is similar to \cite[Lemma 4.7]{ex}. We   have to use  Remark \ref{rim}  instead of \cite[Remark 2.4]{ex}.
      \end{proof}

 Now we prove  that the Milnor's splitting  $\bigoplus\limits_{\deg q =d}\tau_q$ is an isomorphism with the  second residue map $\bigoplus\limits_{\deg q=d}\partial_q$ as  inverse.
      \begin{thm}\label{5ismt}For any positive integer $d$,
    the {Milnor's splitting} $\bigoplus\limits_{\deg q=d} \tau_q: \bigoplus\limits_{\deg q=d} W_1 (H_p^{m+1}(F_q))\xrightarrow{} L_d/{L_{d-1}}$ is    an  isomorphism.
\end{thm}
\begin{proof}
      Repeating the proof  of \cite[Lemma 4.8]{ex}, we  can prove that if $q$, $q'$ are two distinct monic irreducible polynomials over $\F$ of degree $d$, then  $\partial_{q'}\circ \tau_{q}(\phi)=0$ for all $\phi\in  W_1(H_p^{m+1}(F_{q}))$.  Clearly $\partial_q\circ \tau_q(\phi)=\phi$ for all $\phi\in W_1(H_p^{m+1}(F_q))$.  Therefore, we have $$  \bigoplus_{\deg q=d}\partial_q\circ \bigoplus_{\deg q=d} \tau_q={\rm Id}$$ proving    the  injectivity of  $\bigoplus\limits_{\deg q=d} \tau_q$. Then repeating the proof of \cite[Theorem 4.9]{ex} and using {Lemmas} \ref{lgener}, \ref{lim} and Proposition \ref{prsei}, we  prove the surjectivity of this map. 
\end{proof}

{Now we get the following short exact sequence.}

\begin{thm}\label{fsec} The  short sequence 
      $0\xrightarrow{} L_0 \xrightarrow{\textcolor{black}{i_1}} H_p^{m+1}(\mathcal{F}(x))\xrightarrow{\bigoplus\partial_q }\bigoplus\limits_{q}W_1(H_p^{m+1}({F}_q)) \xrightarrow{}0$ is  exact, where   $\textcolor{black}{i_1}$ is the inclusion map and $q$ varies over all monic irreducible polynomials over $\F$.
\end{thm}

\begin{proof}
     It is similar to \cite[Theorem 4.15]{ex}.
\end{proof}
By the decomposition of the group of $W_1(H_p^{m+1}(F_{1/x}))$ of Theorem \ref{thd2}, we have the following result. 
\begin{thm}
    \label{ssec}The following  sequence  
      \[ 0\xrightarrow{} H_p^{m+1}(\mathcal{F}) \xrightarrow{\textcolor{black}{i_2}}L_0\xrightarrow{\textcolor{black}{\partial_{1/x}}} W_1( H_p^{m+1}(F_{\frac{1}{x}}))\xrightarrow{\eta} H_p^{m}(\mathcal{F})\xrightarrow{} 0 \] is  exact, where $\textcolor{black}{i_2}$  is induced by the  inclusion $\F\xrightarrow{} \F(x)$,   $ \eta(\phi)= \phi'$ and $\phi'$  comes from the decomposition of $\phi$, stated in   Theorem  $\ref{thd2}\hspace{1mm} (2)$.  
\end{thm}

\begin{proof}
      Recall that $F=\F(x)$.  Let  $v_x:F\xrightarrow{}\Z\cup\{\n\}$ be the $x$-adic valuation and  $ H_p^{m+1}(F)'$  be the subgroup of $ H_p^{m+1}(F)$ generated by the elements $ a \frac{\dd b_1}{ b_1}\wedge \ldots\wedge \frac{\dd b_m}{b_m}$, where $ b_i\in F^*$ and $a\in F$ such that $v_x(a)\geq 0$.     Let $\beta_x: H_p^{m+1}( F)'\xrightarrow{}   H_p^{ m+1}({\F})$ be the  residue map as defined in Definition \ref{defre}.  For any $a\in \F$ and $b_i\in \F^*$, we have \begin{align*}
           \beta_x\circ i_2\left( a\frac{\dd b_1}{b_1}\wedge \ldots\wedge \frac{\dd b_m}{b_m}\right)&=  \beta_x\left(  a\frac{\dd b_1}{b_1}\wedge \ldots\wedge \frac{\dd b_m}{b_m}\right)\\& =   a\frac{\dd b_1}{b_1}\wedge \ldots\wedge \frac{\dd b_m}{b_m}\in H_p^{m+1}(\F).\\ \beta_x\circ i_2&= \rm{Id}.
      \end{align*} Therefore,  the  homomorphism $i_2$ is injective.  By Theorem \ref{thd2} $(2)$, the map $\eta$ is onto.  \\
       Repeating  the proof of  \cite[Theorem 4.16]{ex}, we  can prove that any element $\phi\in L_0$ can be decomposed as $\phi= \phi_1+ \psi$, where $\phi_1\in H_p^{m+1}(\F)$ and $$\psi = \sum_{ I\in ( T_{\frac{1}{x}})_m}\sum_{J\in T_{\frac{1}{x}},\, l(J)\notin \tilde{I}} t^J r^p_{I,J} \frac{\dd t_I}{t_I},\quad r_{I,J}\in x\F[x].$$   By Theorem \ref{thd2} (2), the form $\psi$ is unique. Therefore, $ {\rm Im}(i_2)= \ker(\partial_{1/x}) $ and $ {\rm Im}(\partial_{1/x})= \ker(\eta)$.  Thus the theorem is proved.
\end{proof}
Now we prove the Milnor's exact sequence in the {setting}  of Kato-Milne cohomology in characteristic $p$.
\begin{thm}\label{rst}Let $\mathcal{F}$ be a field of characteristic $p$ and $F=\mathcal{F}(x)$ be the rational function field {in one variable $x$ over $\F$}. Then, the following \textcolor{black}{short sequence is   exact}:\[ 0\xrightarrow{} H_p^{m+1}(\mathcal{F})\xrightarrow{\textcolor{black}{i}}  H_p^{m+1}(F)\xrightarrow{\displaystyle{\oplus_{q}\partial_q\oplus\partial'_{{1}/{x}}}}\displaystyle{ \oplus_{q}W_1(H_p^{m+1}(F_q))\oplus} \left({ W_1(H_p^{m+1}( F_{\frac{1}{x}}))}/{ H_p^m(\mathcal{F})\wedge { \frac{\mathrm{d}x}{x}}}\right)\xrightarrow{}0,\] where $\textcolor{black}{i}$ is induced by 
 the inclusion $ \mathcal{F}\xrightarrow{}F$ and $q$ varies  over  all monic irreducible polynomials  of $ \mathcal{F}[x]$. Moreover, $\partial'_{{1}/{x}}$ is the composition of $\partial_{{1}/{x}}$ and the projection map. 
  \end{thm}

  \begin{proof}
       Since $i= i_1\circ i_2$, where $ i_1$, $i_2$ come from Theorems \ref{fsec} and \ref{ssec}, respectively, and both $i_1$, $i_2$ are injective, it follows that $i$ is  injective.   Moreover, we have \begin{align*}
            \ker ( \oplus_{q} \partial_q\oplus \partial_{1/x}')= \ker( \oplus_q\partial_q)\cap \ker\partial_{1/x}'&= L_0\cap \ker\partial_{1/x}'\quad \text{(by Theorem \ref{fsec})}\\&=  {\rm Im}(i_2)\quad \text{(by Theorem \ref{ssec})}\\&= {\rm Im}(i).
       \end{align*} Now we prove that $  \oplus_{q} \partial_q\oplus\partial_{1/x}'$ is surjective. Let $ (   u, \overline{v})$ be an  arbitrary element of $ \oplus_{q}W_1(H_p^{m+1}(F_q))\oplus \left({ W_1(H_p^{m+1}( F_{\frac{1}{x}}))}/{ H_p^m(\mathcal{F})\wedge { \frac{\mathrm{d}x}{x}}}\right)$, where $ u\in \oplus_{q}W_1(H_p^{m+1}(F_q))$ and $v\in W_1(H_p^{m+1}( F_{\frac{1}{x}})) $. By Theorem \ref{fsec}, there exists an element $\phi\in  H_p^{m+1}( F)$ such that $$  \bigoplus_{q\in A} \partial_q( \phi)=u$$ where   $A$ is the set of  monic irreducible polynomials of $\F[x]$.  By Theorem \ref{thd2} $(2)$, $v$ can be written as  $v= v'-v''$  such that $ v''\in H_p^{m}(\F)\wedge \dd x/x $ and $$  v'= \sum_{I\in (T_{\frac{1}{x}})_m}\sum_{J\in T_{\frac{1}{x}},\, l(J)\notin \tilde{I}}  t^J r^p_{I,J} \frac{\dd t_I}{t_I}, \quad r_{I,J}\in x\F[x].$$ Let $\partial_{1/x}(\phi)= v_1-v_2$, where $ v_2\in H_p^{m}(\F)\wedge \frac{\dd x}{x}$ and $$  v_1= \sum_{I\in (T_{\frac{1}{x}})_m}\sum_{J\in T_{\frac{1}{x}},\, l(J)\notin \tilde{I}}  t^J {u}^p_{I,J} \frac{\dd t_I}{t_I}, \quad u_{I,J}\in x\F[x].$$ By Theorem \ref{ssec}, there exists an element $ \phi'\in L_0$ such that $\partial_{1/x}(\phi')= v'-v_1$. Now we have \begin{align*}
           \partial_{1/x}'( \phi'+\phi)= \overline{ \partial_{1/x}(\phi'+\phi)}= \overline{v'-v_2}=\overline{v'}=\overline{v}.
       \end{align*} Since $  \phi'\in L_0$,  it follows from Theorem \ref{fsec} that   $\bigoplus\limits_{q\in A} \partial_q(\phi')=0$.\\ Therefore, we have  $ \bigoplus\limits_{q\in A}\partial_q\oplus\partial'_{1/x}(\phi+ \phi')=(u,\overline{v})$, as desired. 
  \end{proof}

\section{A decomposition of  $W_1(H_p^{m+1}(F_q))$ }\label{5sad}

Let $T$, $\overline{T}_q$ and $T_q$ be  as defined in Definition \ref{dpb}, where $q=1/x$ or $q\in \F[x]$ is a monic irreducible polynomial.  For any element $I\in T_q$, the notation $I_{\infty}$ denotes $q$-th component of $I$. 
Now $\overline{T}_q\subset T_q$ and  for $m\geq 1$, we have $$ (\overline{T}_q)_m:=\{I\in (T_q)_m\mid I_{\n}=0\}\quad (T_q)_m\setminus (\overline{T}_q)_m:= \{ I\in (T_q)_m\mid I_{\n}=1\}.$$
{Take into account the definition of $(\overline{T}_q)_m$ introduced {on page} \pageref{p1} without invoking the uniformizer.}

{For any monic irreducible polynomial $q$ over $\F$, let $L_q/F_q$ be a finite separable unramified extension.} By \cite[Chapter 5]{R}, the $q$-adic valuation can be extended to $L_q$ uniquely, and  $L_q$ will be a complete field with respect to this valuation. Let $\overline{L}_q$ be the residue field of $L_q$. Then, $\overline{L}_q/ \overline{F}_q$ is a separable extension. By \cite[Lemma 3.1]{AJ},  $\{ t^I\mid I\in T_q\}$ is {an} $L_q^p$-basis of $L_q$ and $ \{ t^I\mid I\in \overline{T}_q\}$ is {an}  $\overline{L}_q^p$-basis of $\overline{L}_q$. From now on $\overline{g}$ denotes the unique polynomial of degree $<d$ such that $g\equiv \overline{g}\pmod{q}$, where $g\in \F[x]$ and $\deg q=d$. Using these notations, we prove the following proposition:

\begin{prop}\label{pru} Let  $L_q$ and  $\overline{L}_q$ be as stated above for a monic irreducible polynomial $q\in \F[x]$ of degree $d$.   Let $\phi\in H_p^{m+1}(L_q)$  be such that $$ \phi= \sum_{K\in (T_q)_m\setminus (\overline{T}_q)_m} a_K \frac{\dd t_K}{t_K}+ \sum_{L\in (\overline{T}_q)_m} a_L\frac{\dd t_L}{t_L} ,$$ where \begin{itemize}
\item For all $K\in ({T}_q)_m\setminus (\overline{T}_q)_m$, we have $a_{{K}}= \sum\limits_{J\in \overline{T}_q,\,l(J)\notin  \tilde{K}}\sum\limits_{r\geq 1}\frac{\overline{t^Jv_{r,K,J}^p}}{q^{pr}}\in F_q$ for $ v_{r,K,J}\in \F[x]_{< d}$,
\item For $L\in (\overline{T}_q)_m$, we have$$a_{{L}}=\sum\limits_{J\in \overline{T}_q} \sum\limits_{r\geq 1 \atop p\nmid r}\frac{\overline{ t^{J} u_{r,L,J}^p}} {q^{r}}+\sum\limits_{J\in \overline{T}_q,\,l(J)\notin \tilde{ L}}\sum\limits_{r\geq 1}\frac{\overline{t^Jv_{r,L,J}^p}}{q^{pr}}\in F_q,$$for $u_{r,L,J}$, $v_{r,L,J}\in \F[x]_{< d}$,  
\end{itemize}
such that $u_{r,I,J}=v_{r,I,J}=0$ for {almost} all $r$, $I$, $J$. If $\phi=0\in H_p^{m+1}(L_q)$, then $a_K=a_L=0\in F_q$ for all $K$, $L\in (T_q)_m$.
    \end{prop}

    \begin{proof} {The same proof as {of} Lemma} \ref{lu}.
    \end{proof}

Now we prove the required  decomposition of $W_1(H_p^{m+1}(F_q))$    to define the  the transfer map $s_q^*$ in the next section.
\begin{thm}\label{thd3}
     Let $q\in\F[x]$ be a monic irreducible polynomial of degree $d\geq 1$.   Then, any element  $\phi\in W_1(H_p^{m+1}(F_q))$ can be  decomposed as:
     \begin{align*}
         \phi= \sum_{K\in (T_q)_m\setminus (\overline{T}_q)_m} a_K \frac{\dd t_K}{t_K}+ \sum_{L\in (\overline{T}_q)_m} a_L \frac{\dd t_L}{t_L}+ \phi'\wedge \frac{\dd q}{q},
     \end{align*} where 
     \begin{enumerate}
     \item $\phi'\in H_p^m(\F(q))$ is unique.
         \item   For all $K\in (T_q)_m\setminus (\overline{T}_q)_m$, we have $  a_K= \displaystyle{{\sum\limits_{J\in \overline{T}_q, \,l(J)\notin \tilde{K}}} \sum\limits_{r\geq 1 } \frac{ \overline{t^J v^p_{r,K,J}}}{ q^{pr}}}$ with $v_{r,K,J}\in \F[x]_{<d}$.

         \item  For all $L\in (\overline{T}_q)_m$, we have  $$ a_L= \sum\limits_{J\in  \overline{T}_q } \sum\limits_{r\geq 1\atop p\nmid r} \frac{ \overline{t^{J} u^p_{r,L,J}}}{ q^r}+ {\sum_{ J\in \overline{T}_q,\, l(J)\notin \tilde{L}}}\sum_{r\geq 1}\frac{\overline{t^J v^p_{r,L,J}}}{ q^{pr}} $$ with  $u_{r,L,J}, v_{r,L,J}\in \F[x]_{<d}$. Moreover,  $a_K$ and $a_L$ are unique.  Consequently  $u_{r,I,J}$, $v_{r,I,J}$ are unique for all $r$, $I$ and $J$.
     \end{enumerate}
\end{thm}
\begin{proof}
    \textbf{Existence of the decomposition.} Let $\alpha: \overline{F}_q\xrightarrow{} F_q$   be a Teichm\"uller lifting map such that \begin{enumerate}
        \item  for $q$ separable: $ \alpha(\overline{t}_i)=t_i$ for all $i\in \I$,
        \item for $q$ inseparable: $ \alpha(\overline{ t}_i)= t_i$ for all $ i\in \I\setminus \{j_n\}$ and $\alpha(\overline{x})=x$. 
    \end{enumerate}
     Since $ \alpha(\overline{F}_q)\subset F_q$ is the residue field of $F_q$, it follows from Theorem \ref{thd1} that   any  element $ \phi\in W_1(H_p^{m+1}(F_q))$ can be    uniquely  expressed as $ \phi= \psi+ \phi_2\wedge \frac{\dd q}{q}$, where $\phi_2\in H_p^m(\alpha(\overline{F}_q))$ and  $$\psi= \sum\limits_{I\in (T_q)_m}\sum\limits_{J\in T_q, \,l(J)\notin \tilde{I}}   t^J r^p_{I,J}  \frac{\dd t_I}{t_I}\quad \text{ with $ r_{I,J}\in q^{-1}\alpha( \overline{F}_q)[q^{-1}]$}. $$  Since $  t^J r^p_{I,J}= \sum\limits_{i=1}^{l} \frac{\alpha( \overline{g}_i)}{q^i}$, where   $g_i\in \F[x]$ and $ \alpha(\overline{g}_i)= \sum_{j\geq 0}g_{i,j}q^j$  for some $g_{i,j}\in \F[x]_{<d}$, we have $$  \psi= \sum_{I\in (T_q)_m}\sum_{l\leq i} \frac{g_{l,I}}{q^l}\frac{\dd t_I}{t_I}, \quad \text{for some }  g_{l,I}\in \F[x]_{<d} \hspace{2mm} \text{ and } i\in \Z .$$
We prove the decomposition in two steps.\\

 \noindent{\textbf{Step 1:}} Now we prove that  for any $ f\in \F[x]$, $I\in (T_q)_m$ and $l\in \Z$, we have in $  H_p^{m+1}(F)$:

 \begin{equation}\label{eq1}
     \frac{f}{q^l}\frac{\dd t_I}{t_I}= \sum_{K\in (T_q)_m\setminus (\overline{T}_q)_m}a_K \frac{\dd t_K}{t_K}+\sum_{L\in (\overline{T}_q)_m}a_L \frac{\dd t_L}{t_L}+ \sum_{M\in (T_q)_m} f_M \frac{\dd t_M}{t_M}\end{equation}where $f_M\in \F[x]$ and $a_K$, $a_L$  are as defined in Theorem \ref{thd3}.   \\

 If $l\leq 0$, then it is obvious.  Suppose {that} (\ref{eq1}) holds for all  $i<l$ for some  integer  $l\geq 1$. Now we prove (\ref{eq1}) for $l$.  For any $f\in \F[x]$,  there exist $f_J\in \F[x]_{<d}$ and $ k\in \F[x]$ such that:

  \begin{equation}\label{eq2}
       f=\sum_{J\in \overline{T}_q} t^J f^p_J+qk.\end{equation}  We have  two cases depending on {$p\mid l$} or $p\nmid l $.\\

 \noindent{(i)} Let $ l= pr-j$ for some $j\in \{1,\ldots,p-1\}$ and $r\geq 1$.  Now using (\ref{eq2}) we get \begin{align*}
     \frac{f}{q^l}\frac{\dd t_I}{t_I} = \sum_{J\in \overline{T}_q} \frac{t^J f^p_J}{q^l}\frac{\dd t_I}{t_I}+ \frac{k}{q^{l-1}}\frac{\dd t_I}{t_I}= \sum_{J'\in {T}_q\atop t^{J'}=t^Jq^j} \frac{t^{J'} f^p_J}{q^{pr}}\frac{\dd t_I}{t_I}+ \frac{k}{q^{l-1}}\frac{\dd t_I}{t_I}.
 \end{align*}  By Remark \ref{rim},  we  replace $I$ by $I'$ such that $l(J')\notin \tilde{I'}$. Since $J'_{\infty}=j\neq 0$, the element  $I'\in (\overline{T}_q)_m$.  Thus we have:  \begin{align*}
      \frac{f}{q^l}\frac{\dd t_I}{t_I}&= \sum_{I'\in (\overline{T}_q)_m}\sum_{J'\in T_q\atop t^{J'}=t^Jq^j} \frac{ t^{J'} f^p_{I',J}}{q^{pr}}\frac{\dd t_{I'}}{t_{I'}}+ \frac{k}{q^{l-1}}\frac{\dd t_I}{t_I}\quad \text{(for some $ f_{I',J}\in \F[x]_{<d}$)}\\&=  \sum_{I'\in (\overline{T}_q)_m}\sum_{J\in \overline{T}_q} \frac{ t^{J} f^p_{I',J}}{q^l}\frac{\dd t_{I'}}{t_{I'}}+ \frac{k}{q^{l-1}}\frac{\dd t_I}{t_I}\\&= \sum_{I'\in (\overline{T}_q)_m}\sum_{J\in \overline{T}_q} \frac{ \overline{t^{J} f^p_{I',J}}}{q^{l}}\frac{\dd t_{I'}}{t_{I'}}+   \sum_{I'\in (\overline{T}_q)_m}\sum_{J\in \overline{T}_q} \frac{  h_{I',J}}{q^{l-1}}\frac{\dd t_{I'}}{t_{I'}}+ \frac{k}{q^{l-1}}\frac{\dd t_I}{t_I}
 \end{align*}where    $ t^J f^p_{I',J}=\overline{t^Jf^p_{I',J}}+qh_{I',J}$ and $ h_{I',J}\in \F[x]$.  The sum in the  right side is in our required form by induction hypothesis.\\

 \noindent{(ii)}  Let $ l=pr$ for some positive integer $r$. Again by (\ref{eq2}) we get \begin{align*}
     \frac{f}{q^l}\frac{\dd t_I}{t_I}&=  \sum_{J\neq 0\in \overline{T}_q} \frac{t^J f^p_J}{q^{pr}}\frac{\dd t_I}{t_I}+  \frac{f_0}{q^{r}}\frac{\dd t_I}{t_I}+ \frac{k}{q^{l-1}}\frac{\dd t_I}{t_I}\\&= \sum_{I'\in (T_q)_m}\sum_{J\in \overline{T}_q,\, l(J)\notin \tilde{I}} \frac{t^J f^p_{I',J}}{q^{pr}}\frac{\dd t_{I'}}{t_{I'}}+  \left(\frac{f_0}{q^r}+ \frac{k}{q^{l-1}}\right)\frac{\dd t_I}{t_I}\quad ( \text{by Remark } ( \ref{rim}))\\&= \sum_{I'\in (T_q)_m}\sum_{J\in \overline{T}_q,\, l(J)\notin \tilde{I}} \frac{\overline{t^Jf^p_{I',J}}}{q^{pr}}\frac{\dd t_{I'}}{t_{I'}}+ \sum_{I'\in (T_q)_m}\sum_{J\in \overline{T}_q,\, l(J)\notin \tilde{I}} \frac{h_{I',J}}{q^{l-1}}\frac{\dd t_{I'}}{t_{I'}}+ B
 \end{align*} where $f_{I',J}\in \F[x]_{<d}$, $h_{I',J}\in \F[x]$ and $\displaystyle{B= \frac{f_0}{q^r}+\frac{k}{q^{l-1}}}$. Moreover,  $t^J f^p_{I',J}=\overline{t^Jf^p_{I',J}}+qh_{I',J} $. The first sum of right side is in our needed form and other two sums are also in required form by induction hypothesis.  Thus, any polynomial $f$, $I\in (T_q)_m$ and $l\in \Z$ satisfy (\ref{eq1}).  By Remark \ref{rv},  we have  $$  \phi= \sum_{K\in (T_q)_m\setminus (\overline{T}_q)_m} a_K \frac{\dd t_K}{t_K}+ \sum_{L\in (\overline{T}_q)_m} a_L \frac{\dd t_L}{t_L}+ \sum_{M\in (T_q)_m} f_M\frac{\dd t_M}{t_M}+ \phi_2\wedge \frac{\dd q}{q}\in W_1(H_p^{m+1}(F_q)),$$ where  $f_M\in \F[x]_{<d}$ and  $\phi_2\in H_p^m(\alpha(\overline{F}_q))$. 
 Then, repeating the argument of \cite[ Theorem 2.3, Step 2]{Laghribi-Maiti},  the element  $\phi\in W_1(H_p^{m+1}(F_q))$ can be expressed as $$ \phi= \sum_{K\in (T_q)_m\setminus (\overline{T}_q)_m} a_K \frac{\dd t_K}{t_K}+ \sum_{L\in (\overline{T}_q)_m} a_L \frac{\dd t_L}{t_L}+  \phi'\wedge \frac{\dd q}{q}$$ where $\phi'\in H_p^m(\F(q))$.\\

 \noindent{\textbf{Uniqueness of the decomposition:}}
 
 It  is similar to  the uniqueness part of \cite[Theorem 2.3]{Laghribi-Maiti}. We  have to use Proposition \ref{pru} instead of \cite[Proposition 2.2]{Laghribi-Maiti}.
   \end{proof}

   \section{Introducing the transfer map $s_q^*$ }\label{5sit}
   
   {If $q\in \F[x]$ is  a monic irreducible   inseparable  polynomial,   and $I\in T_q$, then} the notation $I_x$ denotes    $x$-th component   of $I$.  Now we define   some subgroups of  $H_p^{m+1}(\F(x))$ depending on  the separability of a monic irreducible polynomial $q\in \F[x]$.

   \begin{defn}[\textbf{The separable case}]\label{d1f} Let $q$ be {monic irreducible separable polynomial} and $r$ be a positive integer.\\
   \noindent{(i)} For $ p\nmid r$, {let $S'_{q,r}$ be the subgroup of}  $S_q+ S_q\wedge \frac{\dd q}{q}$  generated by the forms  $\displaystyle{ \frac{\overline{t^J s^p_{r,I,J}}}{q^r}\frac{\dd t_I}{t_I}}$, where $ I\in (\overline{T}_q)_m$, $J\in \overline{T}_q$ and $s_{r,I,J}\in \F[x]_{<d}$.\\

   \noindent{(ii)} For  {$p\mid r$,  let $ S'_{q,r}$ be the subgroup of} $S_q+ S_q\wedge \frac{\dd q}{q}$ generated by the forms    $ \displaystyle{\frac{\overline{t^J s^p_{r,I,J}}}{q^r}\frac{\dd t_I}{t_I}}$, $\displaystyle{\frac{ \overline{t^{J'} s^p_{r,I',J'}}}{q^r} \frac{\dd t_{I'}}{t_{I'}}\wedge \frac{\dd q}{q}}$,  where $ I\in (\overline{T}_q)_m$, $I'\in (\overline{T}_q)_{m-1}$ and $J$, $J'\in \overline{T}_q$ such that $ l(J)\notin \tilde{I}$, $l(J')\notin \tilde{I'}$ and  $s_{r,I,J}$,  $s_{r,I',J'}\in \F[x]_{<d}$. 
        \end{defn}
\medskip

 \begin{defn}[\textbf{The inseparable case}]\label{d2f} Let $q$ be {monic irreducible inseparable polynomial} and $r$ be a positive integer.\\
 \noindent{(i)} If $p\nmid r$, then we define the subgroups:
 \begin{itemize}
 \item  $  S'_{q,r}\subset \tilde{S}_q+ \tilde{S}_q\wedge \frac{\dd q}{q}$ generated by  the forms $ \displaystyle{\frac{ \overline{t^J s^p_{r,I,J}}}{ q^r}\frac{\dd t_I}{t_I}}$,  where  $I\in (\overline{T}_q)_m$, $J\in \overline{T}_q$ and $s_{r,I,J}\in \F[x]_{<d}$.
 \item  $S^0_{q,r}\subset S'_{q,r}$ generated by the forms $\displaystyle{\frac{\overline{t^Js^p_{r,I,J}}}{ q^r}\frac{\dd t_I}{t_I}}$, where $ J\in \overline{T}_q$, $ s_{r,I,J}\in \F[x]_{<d}$ and $I\in (\overline{T}_q)_m$ such that $I_x=0$.
\end{itemize}
 \noindent{(ii)} If {$p\mid r$}, then we define the subgroups:
 \begin{itemize}
 \item $ S'_{q,r}\subset \tilde{S}_q+ \tilde{S}_q\wedge \frac{\dd q}{q}$ generated by the forms    $ \displaystyle{\frac{\overline{t^J s^p_{r,I,J}}}{q^r}\frac{\dd t_I}{t_I}}$, $\displaystyle{\frac{ \overline{t^{J'} s^p_{r,I',J'}}}{q^r} \frac{\dd t_{I'}}{t_{I'}}\wedge \frac{\dd q}{q}}$,  where $ I\in (\overline{T}_q)_m$, $I'\in (\overline{T}_q)_{m-1}$ and $J$, $J'\in \overline{T}_q$ such that $ l(J)\notin \tilde{I}$, $l(J')\notin \tilde{I'}$ and  $s_{r,I,J}$,  $s_{r,I',J'}\in \F[x]_{<d}$. 
  \item  $ S^0_{q,r}\subset S'_{q,r}$  generated by the forms    $ \displaystyle{\frac{\overline{t^J s^p_{r,I,J}}}{q^r}\frac{\dd t_I}{t_I}}$ and $\displaystyle{\frac{ \overline{t^{J'} s^p_{r,I',J'}}}{q^r} \frac{\dd t_{I'}}{t_{I'}}\wedge \frac{\dd q}{q}}$,  where $J$, $J'\in \overline{T}_q$ and  $ I\in (\overline{T}_q)_m$, $I'\in (\overline{T}_q)_{m-1}$  with $I_x=I'_x=0$  such that $ l(J)\notin \tilde{I}$, $l(J')\notin \tilde{I'}$. Moreover,  $s_{r,I,J},$  $s_{r,I',J'}\in \F[x]_{<d}$. 
 \end{itemize}
\end{defn}
\medskip

 \begin{defn}\label{du}
 For any monic irreducible polynomial  $q$, we define  the subgroup $U_q:= \sum_{r\geq 1} S'_{q,r}\subset H_p^{m+1}(F)$. If $q$ is inseparable, then we define $ U^0_q:= \sum_{r\geq 1} S^0_{q,r}$.\\
 
  Let $ \overline{U}_q$ and $\overline{U^0_q}$ be the images of $ U_q$ {and $U^0_q$}  in $W_1(H_p^{m+1}(F_q))$, respectively.
\end{defn}
Now we estsblish some relations  between $ U_q$ and $ S_q\wedge \frac{\dd x}{x}$.
 \begin{lemma}\label{lp} Let $q$ be a monic irreducible polynomial over $\F$. Then, we have the following:\\
             
 \noindent{(i)} If $q$ is separable, then $ S_q+ S_q\wedge \frac{\dd q}{q}\subset U_q+L_0+ L_0\wedge \frac{\dd q}{q}$.\\

\noindent{(ii)} If $q$ is separable, then  $ S_q+ S_q\wedge \frac{\dd x}{x}\subset U_q+L_0+ L_0\wedge \frac{\dd x}{x}+L_0\wedge \frac{\dd q}{q}$.\\

 \noindent{(iii)} If $q$ is inseparable, then $ \tilde{S}_q+ \tilde{S}_q\wedge \frac{\dd q}{q}\subset U_q+L_0+ L_0\wedge \frac{\dd x}{x}+ L_0\wedge \frac{\dd q}{q}+ L_0\wedge \frac{\dd x}{x}\wedge \frac{\dd q}{q}$.\\

\noindent{(iv)} If $q$ is inseparable, then  $  S_q+ S_q\wedge \frac{\dd x}{x}\subset  U_q+ L_0+L_0\wedge \frac{\dd x}{x}+L_0\wedge \frac{\dd q}{q}+ L_0\wedge \frac{\dd x}{x}\wedge \frac{\dd q}{q}$.\\

\noindent{(v)} The groups  $U_q= \oplus_{r\geq 1} S'_{q,r}$ and $ U_q^0= \oplus_{r\geq 1} S^0_{q,r}$ are direct sums.\\

\noindent{(vi)}  Considering $H_p^m(\F(q))\wedge \frac{\dd q}{q}$ as a subgroup of $H_p^{m+1}(\F(x))$, we get $$  W_1(H_p^{m+1}(F_q))/ \overline{U}_q\cong H_p^m(\F(q))\wedge \frac{\dd q}{q}.$$ 
\end{lemma}

 \begin{proof}
It is similar to \cite[Lemma 3.4]{Laghribi-Maiti}. We have to use Proposition \ref{prsei},  equation (\ref{eq1}) and Theorem \ref{thd3}. \end{proof}

 Now we define some maps which  are needed to define the transfer maps $s_q^*$.\\

\noindent{(1)}  For any monic irreducible polynomial  $q$, let  $ \theta_q: W_1(H_p^{m+1}(F_q))\xrightarrow{}  H_p^m(\F(q))$  be the map defined as:    for any $\phi\in W_1(H_p^{m+1}(F_q))$, $\theta_q(\phi)$ is the unique element of  $H_p^m(\F(q))$ such that $ \phi\equiv \theta_q(\phi)\wedge \frac{\dd q}{q}\pmod{\overline{U}_q}$.\\
               
  Moreover, we take  $\theta_{1/x}:= \eta$ as defined in Theorem \ref{ssec}.\\

\noindent{(2)}  Now we define  the map $t_q^*: H_p^m(\F(q))\xrightarrow{} H_p^m(\F)$ depending on   $q$ separable or not.\begin{itemize}
\item If $q$ is separable, then we take {$t_q^*:= \rm{Tr}_*$}  defined in \cite[Corollary 2.9]{AAB} {and induced by the usual trace map}.
\item Let $q$ be inseparable. Then $q= f(  x^{p^r})$, where $r\geq 1$ and $f(y)\in \F[y]$ is monic irreducible and separable. Let  $K:=\F(f)$ and $g=x^{p^r}-\overline{y}\in K[x]$. Now   the map $ \zeta:\F(q)\xrightarrow{}  K(g)$ defined by $\zeta(c)=c$ for all $c\in \F$ and $\zeta(\overline{x})=\overline{x}$ is a field isomorphism. Then, we take $$ t_q^*:= t_f^*\circ s_*\circ \zeta^*$$ where $s_*: H_p^m(K(g))\xrightarrow{} H_p^m(K)$ is {the transfer map defined in \cite[Proposition 3.4]{ALO}, and }$\zeta^*: H_p^{m}(\F(q))\xrightarrow{} H_p^m(K(g))$ is induced by $\zeta$.           
\item If $q=1/x$, then we take $ t_{q}^*: H_p^m(\F)\xrightarrow{} H_p^m(\F)$ as the identity map.\end{itemize}
Finally we are able to define the transfer map $s_q^*$.

\begin{defn}[{\bf Transfer maps}]\label{dfntr} 
 \noindent{(i)} If $q$ is separable or $1/x$, then {the transfer map $s_q^*: W_1(H_p^{m+1}(F_q))\xrightarrow{} H_p^m(\F)$ is defined as follows:}  $ s_q^*(\phi)= t_q^*( \theta_q(\phi))$ for all $\phi\in W_1(H_p^{m+1}(F_q))$.\\

 \noindent{(ii)} If $q$ is inseparable, then  by Theorem \ref{thd3} we decompose   $\phi- \theta_q(\phi)\wedge \frac{\dd q}{q}$ modulo $\overline{U^0_q}$  as:   \begin{align*} {\sum\limits_{I\in (\overline{T}_q)_m\atop  I_x=1}\sum\limits_{J\in \overline{T}_q}\sum\limits_{i\geq 1\atop p\nmid i}{ \frac{\overline{ t^J u_{{i},I,J}^p}}{q^{i}}\frac{\dd t_I}{t_I}}\;+\sum\limits_{I\in (\overline{T}_q)_m\atop  I_x=1}\sum_{J\in \overline{T}_q\atop l(J)\notin \tilde{ I}}\sum\limits_{i\geq 1}{ \frac{\overline{ t^J v_{i,I,J}^p}}{q^{pi}}\frac{\dd t_I}{t_I}}}\\+{\sum\limits_{I\in (\overline{T}_q)_{m-1}\atop I_x=1}\sum_{J\in \overline{T}_q\atop l(J)\notin\tilde{I}} \sum\limits_{i\geq 1} { \frac{\overline{t^Js^p_{{i},I,J}}}{q^{pi}}\frac{\dd t_I}{t_I}\wedge \frac{\dd q}{q}},}  \end{align*}
 where $u_{i,I,J}$, $v_{i,I,J}$, $s_{i,I,J}\in \F[x]_{<d}$ are unique. For each   $I\in (\overline{T}_q)_m$ with $I_x=1$, we define $ I'\in (\overline{T}_q)_{m-1}$ such that  $t^I= t^{I'}x$.  Now we define the transfer map $ s_q^*:W_1(H_p^{m+1}(F_q))\xrightarrow{} H_p^m(\F) $ {as follows:}   {\begin{align*}  s_q^*(\phi)={{{\textbf{-}}}\sum\limits_{I\in (\overline{T}_q)_m\atop   I_x=1}\sum\limits_{J\in \overline{T}_q}\sum\limits_{i\geq 1\atop p\nmid i}{ \frac{\left(\overline{ t^J u_{{i},I,J}^p}\right)_c}{q_c^{i}}\frac{\dd t_{I'}}{t_{I'}}}\; {{\textbf{-}}}\sum\limits_{I\in (\overline{T}_q)_m\atop  I_x=1}\sum_{J\in \overline{T}_q\atop l(J)\notin \tilde{ I}}\sum\limits_{i\geq 1}{ \frac{\left(\overline{ t^J v_{i,I,J}^p}\right)_c}{q^{pi}_c}\frac{\dd t_{I'}}{t_{I'}}}}\\+{\sum\limits_{I\in (\overline{T}_q)_{m-1}\atop I_x=1}\sum_{J\in \overline{T}_q\atop l(J)\notin\tilde{I}} \sum\limits_{i\geq 1} { \frac{\left(\overline{t^Js^p_{{i},I,J}}\right)_c}{q^{pi}_c}\frac{\dd t_{I'}}{t_{I'}}\wedge \frac{\dd q_c}{q_c}}+ t_q^*(\theta_q(\phi))}  \end{align*}} where $h_c$ denotes the constant part of any polynomial $h\in \F[x]$.
               \end{defn}

  \section{The reciprocity law for $L_0+ L_0\wedge \frac{\dd q}{q}$ }\label{5srl}

  Let $q=\sum_{i=0}^{d}q_ix^{d-i}$ be a monic irreducible polynomial of degree $d$ over $\F$.    We use $x$ to denote the elements of both fields $\F(q)$ and $\F(x)$. Now we define  $\gamma_i\in \F$ such that  for all $i\geq 0$, $$ x^{d+i-1}= \gamma_ix^{d-1}+g_i\in \F(q)$$ where $  g_i\in \F[x]_{\leq d-2}$ are unique and  $\gamma_i=0$ for all $i<0$.  Clearly $ \gamma_0=1$ and $ \gamma_1= -q_1$.   Repeating the arguments of \cite[Section 4]{Laghribi-Maiti}, we get  \begin{equation}\label{eq3} \gamma_i=-\sum_{k=1}^{d} q_k \gamma_{i-k}, \quad \text{for all } i\geq 1.
  \end{equation}   By (\ref{eq3}), we get  the relation \begin{equation}\label{eq4}(1+ q_1x^{-1}+\ldots+ q_dx^{-d}) ( 1+ \gamma_1x^{-1}+ \gamma_2 x^{-2}+\ldots)=1\in \F(x)_{1/x}.   \end{equation} Since $q= x^d( 1+ q_1x^{-1}+\ldots+ q_d x^{-d})$, it follows from (\ref{eq4})    \begin{equation}\label{eq5}q^{-1}= x^{-d}( 1+ \gamma_1x^{-1}+\gamma_2x^{-2}+\ldots)\in \F(x)_{1/x}.
                           \end{equation} 

  \begin{lemma}\label{l81}
  Let $ q= \sum_{i=0}^{d} q_i x^{d-i}$ be a monic irreducible  separable polynomial of degree $d$.  Then $${\rm Tr}_{\F(q)/\F}( x^j)= \sum_{i=0}^{d-1} (d-i)q_i \gamma_{j-i}\;\;\text{for all}\; j\geq 0.$$
  \end{lemma}

 \begin{proof}
Let $ {\rm Tr}_{\F(q)/{\F}}(x^j)=S_j$ for all $j\geq 0$. We prove this lemma by induction on $j$. For $j=0$, we have  $$S_0=d= \sum_{i=0}^{d-1} (d-i)q_i \gamma_{0-i} $$and for $j=1$, it is easy to see $S_1= -q_1= \sum_{i=0}^{d-1} (d-i)q_i \gamma_{1-i} $. \\

Let $k\geq 2$ be an integer and {suppose} $S_j= \sum_{i=0}^{d-1}(d-i) q_i \gamma_{j-i}$ for all $0\leq j<k$.  We have two cases $ k\leq d$ or not.\\

 \noindent{\textbf{Case 1}} {Suppose $k\leq d$.} {By Newton's identities \cite[Proposition 12.6]{M}}, we have \begin{align*}
 S_k&= -\sum_{j=1}^{k-1} q_jS_{k-j}- kq_k\\&=-\sum_{j=1}^{k-1} \sum_{i=0}^{d-1} (d-i)q_j q_i \gamma_{k-j-i}-kq_k\quad ( \text{by induction hypothesis})\\&= -\sum_{j=1}^{k-1} dq_jq_0\gamma_{k-j}- \sum_{j=1}^{k-1}\sum_{i=1}^{k-1}(d-i) q_jq_i\gamma_{k-j-i}-kq_k\quad (\text{since for  $i>k-1$,  $\gamma_{k-j-i}=0$})\\&= d(\gamma_k+ q_k)+\sum_{i=1}^{k-1}(d-i)q_i\gamma_{k-i}-kq_k  \quad (\text{using  (\ref{eq3}) and  $ k-1\geq k-i$ for  $1\leq i\leq k-1$})\\&= \sum_{i=0}^{k}(d-i) q_i \gamma_{k-i}= \sum_{i=0}^{d-1} (d-i) q_i\gamma_{k-i}.
                                \end{align*}

 \noindent{\textbf{Case 2}} {Suppose $k>d$. Again by Newton's identities,} we have\begin{align*}
S_k&= -\sum_{j=1}^{d} q_j S_{k-j}\\&= -\sum_{j=1}^{d} \sum_{i=0}^{d-1}(d-i)q_jq_i\gamma_{k-j-i}\quad \text{(by induction hypothesis)}\\&= \sum_{i=0}^{d-1}(d-i) q_i\gamma_{k-i}\quad \text{(by (\ref{eq3})).}  
 \end{align*}    Thus, by induction hypothesis, we have $ {\rm Tr}_{\F(q)/{\F}}(x^j)= \sum_{i=0}^{d-1} (d-i) q_i \gamma_{j-i}$ for all $j\geq 0$.                             
                           \end{proof}
The following lemma  is crucial    to prove the reciprocity law of Lemma \ref{l85}.
 \begin{lemma}\label{l82}
     Let $q=\sum_{i=0}^{d}q_ix^{d-i}$ be a monic irreducible polynomial of degree $d$.  Then, we have {$$t_q^*\left( ax^j \frac{\dd t_I}{t_I}\wedge \frac{\dd x}{x}\right)= -a\sum_{i=1}^{d} \gamma_{j-i} q_i \frac{\dd t_I}{t_I}\wedge \frac{\dd q_i}{q_i},$$}where $a\in \F$,  $j\geq 1$ and $I\in (T)_{m-2}$.
 \end{lemma}     

 \begin{proof}
      We prove this lemma in two cases depending on $q$ is separable or not. \\

      \noindent{(1)} {Suppose that $q$ is separable, and let }$K:=\F(q)$.  Let  ${\rm Tr}: \Omega_K^m\xrightarrow{} \Omega_{\F}^m$ be  the map induced from the trace map ${\rm Tr}_{K/\F}$, defined in \cite[Section 2]{AAB}.     First of all we compute $ {\rm Tr}( ax^{j-1}\dd x)$.  Since $ q=\sum_{i=0}^{d}q_ix^{d-i}$, we have $\sum_{i=0}^{d} q_i x^{d-i}=0\in K$. Hence, we get 
\begin{align*}
\sum_{i=1}^{d}x^{d-i} \dd q_i&+ \sum_{k=0}^{d-1} (d-k)q_kx^{d-k-1}\dd x=0\in \Omega_K^1,
\end{align*}and thus$$\dd x= -\left( \sum_{k=0}^{d-1} (d-k) q_kx^{d-k-1}\right)^{-1} \left(\sum_{i=1}^{d} x^{d-i} \dd q_i\right).$$
Since $q$ is separable, there exist 
    $l_r\in \F$ for $r\in \{0,\ldots,d-1\}$ such that \begin{equation}\label{eq6}
     \left(\sum_{k=0}^{d-1}(d-k) q_k x^{d-1-k}\right)\left(\sum_{r=0}^{d-1} l_rx^r\right)=1\in K.
 \end{equation}  Now we have \begin{align*}
      {\rm Tr} ( a x^{j-1} \dd x)&= {\rm Tr} \left( -a\sum_{i=1}^{d} \sum_{r=0}^{d-1} l_r x^{d+r+j-i-1} \dd q_i\right)\\&= -a \sum_{i=1}^{d}{\rm Tr}_{K/\F}\left( \sum_{r=0}^{d-1} l_r x^{d+r+j-i-1}\right)\dd q_i.
 \end{align*}{We claim} that  $ {\rm Tr}_{K/\F}\left( \sum_{r=0}^{d-1} l_r x^{ d+r+j-i-1}\right)=\gamma_{j-i}$ for  $j\geq 1$ and $1\leq i\leq d$. {In fact, on the one hand we have} \begin{align*}
     {\rm Tr}_{K/\F}\left( \sum_{r=0}^{d-1} l_r x^{d+r+j-i-1}\right)&= \sum_{r=0}^{d-1} \sum_{k=0}^{d-1}(d-k)l_rq_k\gamma_{d+r+j-i-1-k}\quad \text{(by Lemma \ref{l81})}\\&= \sum_{k=0}^{d-1}\sum_{r=\mu_k}^{d-1}(d-k) l_r q_k \gamma_{d+r+j-i-1-k}\quad \text{(since $\gamma_i=0$ for $i<0$)}\end{align*} where $\mu_k={\rm Max}\{ 0, k+1+i-d-j\}$. 
     
 {on the other hand, multiplying} both sides of (\ref{eq6}) by $ x^{d+j-i-1} $ for $j\geq 1$ and $1\leq i\leq d$, we get$$\sum_{k=0}^{d-1}\sum_{r=0}^{d-1} (d-k)l_rq_kx^{d+(d+r+j-i-1-k)-1}= x^{d+j-i-1}\in K.$$It follows that$$\sum_{k=0}^{d-1}\sum_{r=\mu_k}^{d-1}(d-k) l_r q_k x^{d+(d+r+j-i-1-k)-1}+G_{i,j} = x^{d+j-i-1}\in K,$${where} $G_{i,j}\in \F[x]_{\leq d-2}$. {Since, $\gamma_{l}=0$ for $l<0$,} we get$$\sum_{k=0}^{d-1}\sum_{r=\mu_k}^{d-1} (d-k)l_r q_k \gamma_{d+r+j-i-1-k} x^{d-1}+G'_{i,j}= 
             \gamma_{j-i}x^{d-1}+H_{i,j}\in K,$$where $G'_{i,j}$, $H_{i,j}\in \F[x]_{\leq d-2}$. Since $\{1, x,\ldots,x^{d-1}\}$ is an  $\F$-basis of $K$,   for $j\geq 1$ and $1\leq i\leq d$, we have$$\sum_{k=0}^{d-1}\sum_{r=\mu_k}^{d-1}(d-k) l_r q_k \gamma_{d+r+j-i-1-k} = \gamma_{j-i}\in \F.$$ Therefore, ${\rm Tr}_{K/\F} \left( \sum_{r=0}^{d-1} l_r x^{d+r+j-i-1}\right) =\gamma_{j-i}$. Consequently, ${\rm Tr}(ax^{j-1} \dd x)= -a\sum_{i=1}^{d} q_i\gamma_{j-i}\frac{\dd q_i}{q_i}$. 
             
Finally we get \begin{align*}
             t_q^*\left( ax^j\frac{\dd t_I}{t_I}\wedge \frac{\dd x}{x}\right)& =\frac{\dd t_I}{t_I}\wedge  {\rm Tr}\left(  ax^{j-1}\dd x\right)\quad \text{(by \cite[Proposition 2.4]{AAB})}\\&= -a\sum_{i=1}^{d} q_i\gamma_{j-i}\frac{\dd t_I}{t_I}\wedge \frac{\dd q_i}{q_i}.
         \end{align*}
\noindent{(2)} {Suppose that $q$ is inseparable.}  Then  $q= f(x^{p^r})$, where $r\geq 1$ and $f(y)\in \F[y]$ is a separable polynomial of degree $n\geq 1$. Hence, $\deg q=d= np^r$. We write $f(y) =\sum_{i=0}^{n} f_i y^{n-i}$, where $f_i\in \F$. Then, $ f_i= q_{p^ri}$ for all $i\in \{0,1,\ldots,n\}$ and $q_l=0$ for all $l$ such that $p^r\nmid l$. We use $y$ to denote the elements of $\F(f):= \F[y]/(f)$. For any non-negative integer $i$, we define $\gamma_i^f\in \F$ as $$  {y}^{n+i-1}=\gamma_i^f{y}^{n-1}+G_i\in \F(f)$$  where $G_i\in F[y]_{\leq n-2}$ are unique. Moreover, we take $\gamma_i^f=0$ for all $i<0$.
\medskip

\noindent{{{\bf Claim.} We have that $\gamma_i^f= \gamma_{ip^r}$ for all $i\in \Z$, and $\gamma_i=0$ for all $i$ such that $p^r\nmid i$.}}  

The first part of the claim is obviously true for $i\leq0$. Let it be true for  all $i\leq s-1$ for some $s\geq 1$.  By (\ref{eq3}), we have \begin{align*}
     \gamma_s^f &= -\sum_{i=1}^{n} f_i \gamma^f_{s-i}\\&= -\sum_{i=1}^{n} q_{ip^r}\gamma_{(s-i)p^r} \quad \text{(by induction hypothesis)}\\&= -\sum_{k=1}^{np^r} q_k \gamma_{sp^r-k}\quad \text{(since $q_l=0$ for $p^r\nmid l$ and  $k=ip^r$ )}\\&= \gamma_{sp^r}\quad \text{(by (\ref{eq3}))}.
 \end{align*} 
{Moreover, the second claim is true for all $i\leq 1$. Suppose that} $\gamma_j=0$ for all $j<i$  such that $p^r\nmid j$, where  $i$ is a positive integer such that  $p^r\nmid i$. Now we have \begin{align*}
     \gamma_i&= -\sum_{l=1}^{d=np^r} q_l \gamma_{i-l}=0\quad \text{(by induction hypothesis and $q_l=0$ for $p^r\nmid l$).}
 \end{align*}  
 
 Let $g= x^{p^r}-{y}\in \F(f)[x]$. Let $s_*: H_p^m(K(g))\xrightarrow{} H_p^m(K)$ be the  map as defined in \cite[Proposition 3.4]{ALO} where $K= \F(f)$.  Now $\F(q)\cong K(g)$. For $a\in \F$ and $I\in (T)_{m-2}$,  we first compute the following:\begin{align*}
     s_*\left( ax^j\frac{\dd t_I}{t_I}\wedge \frac{\dd x}{x}\right)= s_*\left( a^{p^r}x^{jp^r} \frac{\dd t_I}{t_I}\wedge \frac{\dd x}{x}\right)&= s_*\left( a^{p^r} {y}^j\frac{\dd t_I}{t_I}\wedge \frac{\dd x}{x}\right)\\&= a^{p^r} {y}^j\frac{\dd t_I}{t_I}\wedge \frac{\dd {y}}{{y}}\quad \text{(by \cite[Proposition 3.4]{ALO})}.
 \end{align*} Thus, we have \begin{align*}
     s_*\left( ax^j \frac{\dd t_I}{t_I}\wedge \frac{\dd x}{x}\right)=\begin{cases}
         a {y}^{(j/p^r)} \frac{\dd t_I}{t_I}\wedge \frac{\dd {y}}{{y}}&\text{if $p^r\mid j$,}\\0&\text{if $p^r\nmid j$.} \end{cases}
 \end{align*}   If {$p^r\mid j$}, then by the definition of $t_q^*$ we have   \begin{align*}
     t_q^*\left( ax^j \frac{\dd t_I}{t_I}\wedge \frac{\dd x}{x}\right)&= t_f^*\left( a{y}^{(j/p^r)} \frac{\dd t_I}{t_I}\wedge \frac{\dd {y}}{{y}}\right)\\&= -a\sum_{i=1}^{n}f_i \gamma^f_{(j/p^r)-i} \frac{\dd t_I}{t_I}\wedge \frac{\dd f_i}{f_i}\quad \text{(by (i))}\\&= -a\sum_{i=1}^{n} q_{ip^r} \gamma_{j-ip^r}\frac{\dd t_I}{t_I}\wedge \frac{\dd q_{ip^r}}{q_{ip^r}}\quad \text{(since $f_i= q_{ip^r}$ and $\gamma^f_i= \gamma_{ip^r}$)}\\&= -a\sum_{k=1}^{np^r} q_k \gamma_{j-k}\frac{\dd t_I}{t_I}\wedge \frac{\dd q_k}{q_k}\quad \text{(since $q_l=0$ for $p^r\nmid l$ and $k=ip^r$)}.
 \end{align*}  Thus, for all $j\geq 1 $, $a\in \F$ and $I\in (T)_{m-2}$,  we have \begin{align*}
      t_q^*\left( ax^j \frac{\dd t_I}{t_I}\wedge \frac{\dd x}{x}\right)&= \begin{cases} -a \sum\limits_{i=1}^{d} q_i\gamma_{j-i} \frac{\dd t_I}{t_I}\wedge \frac{\dd q_i}{q_i}\quad \text{ if $p^r/j$}\\ 0\hspace{37mm} \text{if $p^r\nmid j$}\end{cases}\\&=-a\sum_{i=1}^{d=np^r} q_i \gamma_{j-i}\frac{\dd t_I}{t_I}\wedge \frac{\dd q_i}{q_i} \quad \text{(since $q_l=\gamma_l=0$ for $p^r\nmid l$)}.
 \end{align*}\end{proof}

 \begin{lemma}\label{l83} Let $q=\sum_{i=0}^{d}q_ix^{d-i}$ be a monic irreducible and inseparable polynomial  of degree $d$.   Then ${t_q^*\left( a\frac{\dd t_I}{t_I}\wedge \frac{\dd x}{x}\right)= a\frac{\dd t_I}{t_I}\wedge \frac{\dd q_d}{q_d}}$ for all $a\in \F$ and $I\in (T)_{m-2}$.\end{lemma}
\begin{proof}
     Since $q$ is inseparable, we write $q= f(x^{p^r})$ for some separable polynomial $f(y)\in \F[y]$. Let $K= \F(f)$ and $g= x^{p^r}-{y}\in K[x]$.    By definition, \begin{align*}
         t_q^*\left( a\frac{\dd t_I}{t_I}\wedge \frac{\dd x}{x}\right) &=  t_f^*\left( s_*\left(\zeta^*\left(  a\frac{\dd t_I}{t_I}\wedge \frac{\dd x}{x}\right)\right)\right)\\&= t_f^*\left(s_*\left(  a\frac{\dd t_I}{t_I}\wedge \frac{\dd x}{x}\right)\right)\\&=  t_f^*\left( a\frac{\dd t_I}{t_I}\wedge \frac{\dd {y}}{{y}}\right)\quad \text{(by \cite[Proposition 3.4, (ii)]{ALO})} \\&= a\frac{\dd t_I}{t_I}\wedge \frac{\dd q_d}{q_d}\quad \text{(by \cite[Lemma 2.5, Proposition 2.4]{AAB})}.  \end{align*}\end{proof}
    Now we  compute $s^*_{1/x}(\partial_{1/x}(\phi))$ for all $\phi\in L_0\wedge \frac{\dd q}{q}$  and for any monic irreducible polynomial $q$.
\begin{lemma}\label{l84}
    Let $q=\sum_{i=0}^{d} q_ix^{d-i}$ be a monic   irreducible polynomial of degree $d$ over $\F$. \\
 \noindent{(i)} If $\phi= ax^j\frac{\dd t_I}{t_I}\wedge \frac{\dd q}{q}\in H_p^{m+1}(\F(x))$, where $j\geq 0$, $I\in (T)_{m-1}$, then $$s^*_{1/x}\circ \partial_{1/x}(\phi)= -\sum_{i=0}^{d-1}(d-i) aq_i\gamma_{j-i}\frac{\dd t_I}{t_I} .$$ In particular,   if $q$ is inseparable,  then $s^*_{1/x}\circ \partial_{1/x}(\phi)=0$.\\
 
 \noindent{(ii)} If $\phi= ax^j\frac{\dd t_I}{t_I}\wedge \frac{\dd x}{x}\wedge \frac{\dd q}{q}\in H_p^{m+1}(\F(x))$, where $j\geq 1$ and $I\in (T)_{m-2}$, then  $$  s^*_{1/x}\circ \partial_{1/x}(\phi)= \sum_{i=1}^{d} aq_i \gamma_{j-i} \frac{\dd t_I}{t_I}\wedge \frac{\dd q_i}{q_i}.$$
\end{lemma}
 \begin{proof}
     (i)  By (\ref{de}), we have in $H_p^{m+1}(\F(x))$: \begin{align*}
         \phi= ax^j \frac{\dd t_I}{t_I}\wedge \frac{\dd q}{q}&= \sum_{i=0}^{d}  \frac{aq_i x^{j+d-i}}{q}\frac{\dd t_I}{t_I}\wedge \frac{\dd (q_ix^{d-i})}{q_ix^{d-i}}\\&= \sum_{i=0}^{d}  \frac{aq_i x^{j+d-i}}{q}\frac{\dd t_I}{t_I}\wedge \frac{\dd q_i}{q_i}+ \sum_{i=0}^{d-1}  \frac{(d-i)aq_i x^{j+d-i}}{q}\frac{\dd t_I}{t_I}\wedge \frac{\dd x}{x}
     \end{align*} Then by repeating the proof of \cite[Lemma 4.3, (i)]{Laghribi-Maiti} and using (\ref{eq5}) we can prove the lemma.\\   

     If  $q$ is inseparable, then $d=np^r$ for some positive integers $n$, $r$.   We have already proved in Lemma \ref{l82} that   $ q_i=\gamma_i=0$  for all $i\in \{0,\ldots,d\}$ such that  $p^r\nmid i$. Thus $ s^*_{1/x}\circ \partial_{1/x}(\phi)=0$.\\
     
     \noindent{(ii)}  It is similar to \cite[Lemma 4.3, (ii)]{Laghribi-Maiti}.

 \end{proof}

 \begin{lemma}\label{l85}
     Let $q=\sum_{i=0}^{d} q_ix^{d-i}$ be a monic  irreducible polynomial of degree $d$. Then we have:\\
     \noindent{(i)} $\sum_{q'} s_{q'}^*( \partial_{q'}(\phi))=0$ for all $\phi\in L_0+ L_0\wedge \frac{\dd q}{q}$, where $q'$ varies over $1/x$ and monic irreducible polynomials over $\F$.\\
     \noindent{(ii)} If $q$ is inseparable, then  $\sum_{q'}s_{q'}^*( \partial_{q'}(\phi))=0$ for all $\phi\in L_0\wedge \frac{\dd x}{x}\wedge \frac{\dd q}{q}$.
 \end{lemma}
 \begin{proof}
     (i)    Recall that $L_0$ is generated by   two following forms   $ ax^i \frac{\dd t_I}{t_I}$ and $bx^j \frac{\dd t_J}{t_J}\wedge \frac{\dd x}{x}$, where $i\geq 0$, $j\geq 1$, $I\in (T)_m$ and $J\in (T)_{m-1}$.  Let $\phi= ax^i \frac{\dd t_I}{t_I}\in L_0$  be a first type generator.  By  Lemma \ref{lvalu},   we have $ \partial_{q'}(\phi)=0$, for any monic irreducible polynomial $q'$.  By Theorem \ref{ssec},   $$  s_{1/x}^*(\partial_{1/x}(\phi))=0\quad (\text{since $s^*_{1/x}=\eta$} )$$Let $\phi= bx^j \frac{\dd t_J}{t_J}\wedge \frac{\dd x}{x}$ be a second type generator of $L_0$. By Lemma \ref{lvalu}, $\partial_{q'}(\phi)=0$ for all monic irreducible polynomial $q'\neq x$. Since $j\geq 1$, $\partial_x(\phi)=0$. By Theorem \ref{ssec}, we have: $$ s_{1/x}^*(\partial_{1/x}(\phi))=0.$$  Thus we get $$ \sum_{q'} s^*_{q'}(\partial_{q'}(\phi))=0, \quad \text{for all } \phi\in L_0.$$
     Moreover, $L_0\wedge \frac{\dd q}{q}$ is generated by    the forms $  ax^i \frac{\dd t_I}{t_I}\wedge \frac{\dd q}{q}$ and $bx^j \frac{\dd t_J}{t_J}\wedge \frac{\dd x}{x}\wedge \frac{\dd q}{q}$, where $i\geq 0$, $j\geq 1$, $I\in (T)_{m-1}$ and $J\in (T)_{m-2}$.  Let $\phi= ax^i \frac{\dd t_I}{t_I}\wedge \frac{\dd q}{q}$ be a first type generator of $L_0\wedge \frac{\dd q}{q}$.      By Lemma \ref{lvalu}, $ \partial_{q'}(\phi)=0$  for any monic irreducible polynomial $q'\neq q$.   Now we have \begin{align*}
         s_q^*(\partial_q(\phi))&= t_q^*\left(  ax^i \frac{\dd t_I}{t_I}\right)\\&= \begin{cases}
            a \sum_{k=0}^{d-1}(d-k)q_k\gamma_{i-k}\frac{\dd t_I}{t_I},\hspace{1mm}\text{if $q$ is separable}\hspace{2mm}\text{(by Lemma \ref{l81})}\\0,\hspace{20mm}\text{if $q$ is inseparable}\hspace{2mm} \text{(by \cite[Proposition 3.4]{ALO})}
         \end{cases}\\&= -s_{1/x}^*(\partial_{1/x}(\phi)),\quad \text{(by Lemma \ref{l84})}
     \end{align*} Let $\phi= bx^j \frac{\dd t_J}{t_J}\wedge \frac{\dd x}{x}\wedge \frac{\dd q}{q}$ be the second type generator of $L_0\wedge \frac{\dd q}{q}$.   By Lemma \ref{lvalu}, $\partial_{q'}(\phi)=0$, for any monic irreducible polynomial $q'\neq q$, $x$. Since $j\geq 1$, $\partial_{x}(\phi)=0$.  Moreover, 
     \begin{align*}
         s_q^*(\partial_q(\phi))&= t_q^*\left(  bx^j \frac{\dd t_J}{t_J}\wedge \frac{\dd x}{x}\right)\\&= -b\sum_{i=1}^{d} \gamma_{j-i} q_i \frac{\dd t_J}{t_J}\wedge \frac{\dd q_i}{q_i}\quad \text{(by Lemma \ref{l82})}\\ &= -s_{1/x}^*(\partial_{1/x}(\phi))\quad \text{(by Lemma \ref{l84})}
     \end{align*} Thus we have  $\sum_{q'} s^*_{q'}(\partial_{q'}(\phi))=0$ for all $\phi\in L_0+ L_0\wedge \frac{\dd q}{q}$.\\

     \noindent{(ii)}It is  similar to \cite[Lemma 4.4, (ii)]{Laghribi-Maiti}. We have to use Lemma \ref{l83} instead of \cite[Lemma 4.2 (iii)]{Laghribi-Maiti}.
 \end{proof}

 \section{Milnor-Scharlau's exact sequence}\label{5sms}
 
 Here we prove   the Milnor-Scharlau exact sequence in the notion of Kato-Milne cohomology group in characteristic $p$. To this end we  extend the reciprocity law for $H_p^{m+1}(\F(x))$.
 \begin{thm}\label{thms}     For any differential form $\phi\in H_p^{m+1}(\F(x))$, we have $ \sum_{q}s^*_{q}(\partial_q(\phi))=0$, where $q$ varies over  $1/x$ and all monic irreducible polynomials over $\F$.
     \end{thm}

     \begin{proof}
           By Lemmas \ref{lgener}, \ref{lp}, \ref{l85} and Definition \ref{du}, it is sufficient to prove that  \[ \sum_{q} s_q^*(\partial_q(\phi))=0\] where $\phi$ is a generator   of  $ S'_{f,r}$ as stated in Definitions \ref{d1f}, \ref{d2f}, where $r\in \N_0$ and $f$ is a monic irreducible polynomial over $\F$.   We prove this in $4 $ cases.\\

           \noindent{(i)} Let ${\phi={\frac{\overline{t^J s^p_{r,I,J}}}{f^r}\frac{\dd t_I}{t_I}}}\textcolor{black}{\in S'_{f,r}}$, where  $ J\in \overline{T}_f,r\geq 1 $,  $s_{r,I,J}\in \F[x]_{< \deg f}$ and  $I\in ({T}_f)_{m}$ such that $t^I\in \F$.    Moreover,  $l(J)\notin \tilde{I}$ when $p\mid r$.     Repeating the proof of \cite[Theorem 5.1,  (i)]{Laghribi-Maiti}, we  prove that  \[ \sum_{q} s_q^*(\partial_q(\phi))=0.\]

           \noindent{(ii)} Let  ${\phi={\frac{\overline{t^J s^p_{r,I,J}}}{f^r}\frac{\dd t_I}{t_I}\wedge \frac{\dd f}{f}}}\textcolor{black}{\in S_{f,r}'} , $ where  $p\mid r$,  $ s_{r,I,J}\in \F[x]_{< \deg f}$ and $I\in ({T}_f)_{m-1}$ such that $t^I\in \F$. Moreover, $J\in \overline{T}_f\subset T_f$ such that $l(J)\notin \tilde{I}$.  Repeating the proof of  \cite[Theorem 5.1, (ii)]{Laghribi-Maiti}, we prove that  \[ \sum_{q} s_q^*(\partial_q(\phi))=0.\]  \\

The elements  of (i) and (ii)  generate $U_f$ for  $f$ separable  and $  U_f^0$ for $f$ inseparable.  Now we check for the remaining generators of $U_f\setminus U^0_f$ for $f$ inseparable.\\

            \noindent{(iii)} Let $ {\phi={\frac{\overline{t^J s^p_{r,I,J}}}{f^r}\frac{\dd t_I}{t_I}\wedge \frac{\dd x}{x}}}\textcolor{black}{\in S'_{f,r}},$ where $J\in \overline{T}_f, r\geq 1$, $s_{r,I,J}\in \F[x]_{< \deg f}$ and  $I\in ({T}_f)_{m-1}$ such that $t^I\in \F$.
      If  $p/r$, then  $l(J)\notin \tilde{I'}$, where $t^{I'}=t^I x$.  By Lemma \ref{lvalu},  $\partial_q(\phi)=0$ for all  $q$ except $f$, $x$ and $1/x$.  Since  $\deg (\overline{t^J s^p_{r,I,J}})< \deg (f)$, it follows from Remark \ref{rv} that $\partial_{1/x}(\phi)=0$. Now we have  $ \displaystyle{\frac{ \overline{t^Js^p_{r,I,J}}}{f^r}= \frac{ \overline{(t^Js^p_{r,I,J})}_c}{f^r_c}+u}\in \F((x))$, where $   v_x(u)\geq 1$ and $h_c$ denotes the constant term of a polynomial  $h\in \F[x]$.  Therefore, we have \begin{align*}
         s_x^*(\partial_x(\phi))= s_x^*\left( { \frac{\left(\overline{t^J s^p_{r,I,J}}\right)_c}{f_c^r}\frac{\dd t_I}{t_I}\wedge \frac{\dd x}{x}}\right)&= { \frac{\left(\overline{t^J s^p_{r,I,J}}\right)_c}{f_c^r}\frac{\dd t_I}{t_I}}\\&= - s_f^*(\partial_f(\phi))
\quad \text{(by Definition \ref{dfntr} (ii))}     \end{align*}   
            Thus we get $\sum_qs_q^*(\partial_q(\phi))=0$.\\

            \noindent{(iv)} Let  $ {\phi={\frac{\overline{t^Js^p_{r,I,J}}}{f^r}\frac{\dd t_I}{t_I}\wedge \frac{\dd x}{x}\wedge \frac{\dd f}{f}}}\textcolor{black}{\in S'_{f,r}},$ where  $p\mid r$, $ s_{r,I,J}\in \F[x]_{< \deg f}$ and $I\in ({T}_f)_{m-2}$ such that $t^I\in \F$. Moreover, $J\in \overline{T}_f$ such that $l(J)\notin \tilde{I'}$, where $t^{I'}=t^Ix$.   As in part (iii),  we have  $\partial_q(\phi)=0$ for all $q$ except $f$ and $x$.  Let $f=\sum_{i=0}^{d} f_i x^{d-i}$, where $f_0=1$. Let $h$ denote the polynomial $ \overline{t^Js^p_{r,I,J}}$ of degree $\leq d-1$. Then by (\ref{de}), we get \begin{align*}
                \phi&= \sum_{i=0}^{d} { \frac{hf_ix^{d-i}}{f^{r+1}}\frac{\dd t_I}{t_I}\wedge \frac{\dd x}{x}\wedge \frac{\dd f_i}{f_i}}\\&= -\sum_{i=1}^{d} { \frac{hf_ix^{d-i}}{f^{r+1}}\frac{\dd t_I}{t_I}\wedge \frac{\dd f_i}{f_i}\wedge \frac{\dd x}{x}} \quad\textcolor{black}{(\text{as }  f_0=1)}.
            \end{align*}  Therefore, we have:  \begin{align*}
                 s_x^*(\partial_x(\phi))&=   s_x^*\left(-\frac{h_cf_d}{f_d^{r+1}}\frac{\dd t_I}{t_I}\wedge \frac{\dd f_d}{f_d}\wedge \frac{\dd x}{x}\right)\quad (\text{since  } v_x\left( \frac{hf_ix^{d-i}}{\textcolor{black}{f^{r+1}}}\right)\geq 1, \text{ for } \hspace{1mm} i<d ).\\&=  - \frac{h_c}{f_d^r}\frac{\dd t_I}{t_I}\wedge \frac{\dd f_d}{f_d}= -s_f^*(\partial_f(\phi)).
            \end{align*}  Finally we get   $\sum_qs_q^*(\partial_q(\phi))=0$ which proves the theorem. 
     \end{proof}
Finally we are able to  prove  the main result of this note.
     \begin{thm}\label{thml} Let $\F$ be a field of characteristic $p$   and $F=\F(x)$ be the rational function field in one variable $x$ over $\F$. Then we have the following exact sequence:
     $$  0\xrightarrow{} H_p^{m+1}(\F)\xrightarrow{i} H_p^{m+1}(F)\xrightarrow{\oplus \partial_q} \bigoplus_{q}W_1(H_p^{m+1}(F_q))\xrightarrow{\oplus s_q^*} H_p^m(\F)\xrightarrow{} 0 $$  where $ i$ is induced by the inclusion $\F\xrightarrow{} F$. Moreover, $q$  varies over  $1/x$ and monic irreducible polynomials over $\F$.
          \end{thm}

          \begin{proof}
               The map $i=i_1\circ i_2$, where $i_1$, $i_2$ come from Theorem \ref{fsec}, \ref{ssec} respectively. Since  $i_1$, $i_2$ are injective, $i$ is injective. Moreover, it follows from  Theorems \ref{fsec} and \ref{ssec} that $$\ker(\oplus \partial_q)= {\rm Im}(i_2)={\rm Im}(i).$$   By Theorem \ref{ssec}, the map $s^*_{1/x}$ is onto and thus $\oplus s_q^*$ is onto.    By the  previous theorem,  ${\rm Im}(\oplus \partial_q)\subset {\ker}(\oplus s_q^*)$. Now repeating the proof of Theorem \cite[Theorem 5.2]{Laghribi-Maiti},  we can prove  $ \ker (\oplus s_q^*)\subset {\rm Im} (\oplus \partial_q)$.   Thus the theorem is proved.
          \end{proof}

           \section{Norm Theorem}\label{5sn}
           
           Here we aim to give an application of Theorem \ref{thml} regarding the norm theorem in the notion of Kato-Milne cohomology in characteristic $p>0$. Now we define hyperbolicity  and  norm of a differential form,  which were  introduced by Mukhija \cite[Definition A.1, A.2]{Appendix}.  
  \begin{defn}
 \label{dfh} Let $F$ be a field of characteristic $p$ and $w\in \Omega_F^m$ be a differential form. Then $w$     is said to be hyperbolic if $w=0\in H_p^{m+1}(F)$, i.e., $w\in \wp(\Omega_F^m)+\dd \Omega_F^{m-1}$.
  \end{defn}
  This definition relates the Witt group of nonsingular quadratic forms  $W_q(F)$ with  the Kato-Milne cohomology $H_p^{m+1}(F)$.  Now we extend the  notion of norm to differential forms.   Let $c\in F^*$ and $\phi$ be a nonsingular quadratic form over $F$. Then $c$ is said to be a norm of $\phi$ if $\phi\cong c\phi$, i.e., the quadratic form $ \langle 1,c \rangle_b\otimes \phi$ is hyperbolic, where $\langle 1,c \rangle_b$ is the  $F$-bilinear form $ X^2+cY^2$. For $p=2$,  it follows from  Kato's isomorphism \cite[Equation (17)]{Appendix}  that $ \overline{\langle 1,c \rangle_b}\in I(F)/I^2(F)$ corresponds to $ \frac{\dd c}{c}\in \nu_F(1) $. This motivates the following definition.
\begin{defn}\label{dn}
     Let $c\in F^*$ and $w\in \Omega_F^m$. The element  $c$ is said to be a norm of $w$ if $w\wedge \frac{\dd c}{c}\in \Omega_F^{m+1}$ is hyperbolic. 
\end{defn}
Let $F$ be a field of characteristic $p$ and $K/F$ be a finite extension. Then the notation $H_p^{m+1}(K/F)$ denotes the group $\ker(i:H_p^{m+1}(F)\xrightarrow{} H_p^{m+1}(K))$, where $i$ is induced from the inclusion $F\xrightarrow{}K$. The following theorem extends \cite[Theorem A.16]{Appendix} to characteristic $p$.
     \begin{thm}\label{thns} Let  $K/F$ be a field extension such that $p\nmid [K:f]$. Then the group   $H_p^{m+1}(K/F)=0$.
          \end{thm}
          \begin{proof}
           Let $w\in \Omega_F^m$  be  a differential form such that $ w\in H_p^{m+1}(K/F)$.  Since $p\nmid [K:F]$, the extension $K/F$ is separable.   By applying the transfer map ${\rm Tr}_*: H_p^{m+1}(K)\xrightarrow{} H_p^{m+1}(F)$, defined in \cite[Section 2]{AAB}, we get $ {\rm Tr}_*(w_K)=0\in H_p^{m+1}(F)$.  Since $ {\rm Tr}_*(w_K)= [K:F] w$ and $p\nmid [K:F]$, we have $w=0\in H_p^{m+1}(F)$.
          \end{proof} The norm theorem  was first introduced by Knebusch for bilinear forms, which states that: \textit{If $\F$ is a field of arbitrary characteristic and $q\in \F[x_1,\ldots,x_n]$ is a normed  irreducible polynomial in $n$ variables, then an anisotropic bilinear form $B$ over $\F$ becomes  metabolic over $\F(q)$ if and only if $ B\cong qB$ over $\F(x_1,\ldots,x_n)$, where $\F(q)$ is the quotient field of $\F[x_1,\ldots,x_n]/(q)$}, \cite[Theorem 4.2]{Knebusch}. Thus we get the norm theorem for quadratic forms in characteristic different from $2$. Later Baeza    proved  the norm theorem for nonsingular quadratic forms in characteristic $2$, \cite{Baeza}.   Then first author \cite{Laghribi1} and Hoffman \cite{Hoffmann1}  proved the  norm theorem  independently for singular quadratic forms in characteristic $2$.  Finally,    first author and Mukhija proved  the   norm theorem for  semisingular quadratic forms in characteristic $2$ in \cite{LM}.  Using the transfer defined by Aravire, Laghribi and O'Ryan for purely inseparable extensions and Definitions \ref{dfh}, \ref{dn}, Mukhija \cite[Theorem A.3]{Appendix} proved the norm theorem in the setting of Kato-Milne cohomology for  inseparable polynomials. Recently by proving the Milnor-Scharlau exact sequence,  we prove  the norm theorem for Kato-Milne cohomology in characteristic $2$ for any normed irreducible polynomial, \cite[Theorem 5.8]{ex}. Now we aim to prove this theorem in any positive characteristic $p$.      At first in the following theorem we deal with the case of  monic irreducible polynomial in one variable  and then the general case in Theorem \ref{n2th}. 
            \begin{thm}\label{thno} Let $\F$ be a field of characteristic $p$, $F=\F(x)$ and $ w\in \Omega_{\F}^m$. Let $q\in \F[x]$ be a monic irreducible polynomial. Then the following statements are equivalent:\\
          \noindent{(1)}  $w$ is hyperbolic over $\F(q)$.\\
          \noindent{(2)} $q$ is a norm of $w_F$. 
           \end{thm}

           \begin{proof}
                (2)$\Rightarrow{}$(1): Let $H_p^{m+1}(F)'$ be the subgroup of $H_p^{m+1}(F)$ generated by the elements $ a\frac{\dd b_1}{b_1}\wedge\ldots\wedge \frac{\dd b_m}{b_m}$, where $b_i\in F^*$ and $a\in F$ such that $v_q(a)\geq 0$.   Let $\beta:H_p^{m+2}(F)'\xrightarrow{} H_p^{m+1}(\F(q))$ be the map  as defined in Definition \ref{defre} with respect to  $q$-adic valuation.  If $q$ is a norm of $w_F$, then we have:\begin{align*}
                    \beta\left(w_F\wedge \frac{\dd q}{q}\right)= w_{\F(q)}=0.
                \end{align*} Therefore, $w$ is hyperbolic over $\F(q)$.\\
                
                \noindent{(1)$\Rightarrow{}(2)$:} Let $w$ be hyperbolic over $\F(q)$.  If $p\nmid \deg q$, then it follows from  Theorem \ref{thns} that $w=0\in H_p^{m+1}(\F)$. Therefore, $w_F\wedge \frac{\dd q}{q}=0\in H_p^{m+2}(F)$. \\ Let $\deg q= pr$, where $r\geq 1$.  Since $w=0\in H_p^{m+1}(\F(q))$, it follows from Teichm\"uller lifting that  $w=0\in H_p^{m+1}(F_q)$. Therefore $$\partial_q\left( w_F\wedge \frac{\dd q}{q}\right)=0.$$ Moreover, for any monic irreducible polynomial $q'\neq q$, we have $$ \partial_{q'}\left(  w_F\wedge \frac{\dd q}{q}\right)=0\quad \text{(by Lemma \ref{lvalu})}$$Now we write  $q= \sum\limits_{i=0}^{pr} q_i x^{pr-i} \in \F[x]$ where $q_0=1$. Thus $q=  x^{pr} y$, where $ y= \sum\limits_{i=0}^{pr} q_i x^{-i}$. Therefore, in $H_p^{m+2}(F_{1/x})$ we have: \begin{align*}
                       w\wedge \frac{\dd q}{q}= w\wedge \frac{\dd y}{y}&= \sum_{i=0}^{pr} \frac{ q_ix^{-i}}{y}w\wedge \frac{\dd (q_i x^{-i})}{q_ix^{-i}} \\&= \frac{q_0}{y} w\wedge \frac{\dd q_0}{q_0}\quad \text{(since $v_{1/x}( {q_ix^{-i}}/{y})\geq 1$ for all $i\geq 1$)}\\&= 0\quad \text{(since $q_0=1$)}
                \end{align*} Hence,  $\partial_{1/x}\left( w_{F}\wedge \frac{\dd q}{q}\right)=0$. Therefore, $w_{F}\wedge \frac{\dd q}{q}\in \ker(\oplus_{q'}\partial_{q'})$, where $q'$ varies over $1/x$ and monic irreducible polynomials over $\F$. By Theorem \ref{thml},  $ w_F\wedge \frac{\dd q}{q}=i(\alpha)$ for some $\alpha\in H_p^{m+2}(\F)$.   Applying $\beta_{1/x}$ of Definition \ref{defre} we get: \begin{align*}
                     \beta_{1/x}( i(\alpha))&= \beta_{1/x}\left(w_F\wedge \frac{\dd q}{q}\right)\\ \alpha &=\beta_{1/x}\left( w_F\wedge \frac{\dd y}{y}\right)\\&= w\wedge \frac{\dd \overline{y}}{\overline{y}}\\&=0\quad \text{(since $\overline{y}=1\in \overline{F}_{1/x}$)}
                \end{align*} Hence, $ w_F\wedge 
                \frac{\dd q}{q}=i(\alpha)=0\in H_p^{m+2}(F)$, as required.
           \end{proof}
           The following result is crucial to prove the norm theorem for arbitrary number of  variables.
 \begin{prop}
     \label{prn}
     Let $w\in \Omega_{\F}^m$ be a differential form and  $f\in \F[x_1,\ldots,x_n]$  be a polynomial.   Let $F=\F(x_1,\ldots,x_n)$ and  $q$ be a normed irreducible  polynomial such that  $ f =q^r g$, where $ p\nmid r\in \N_0$ and  $q\nmid g\in \F[x_1,\ldots,x_n]$. If $ f$ is a norm of $w_F$, then $q$ will  be  a norm of $w_F$. 
 \end{prop}
\begin{proof}We prove the proposition by induction on $n$. Let $n=1$. 
    Since $f$ is a norm of $w_F$, we have:\begin{align*}
        w_F\wedge \frac{\dd f}{f}=r w_F\wedge \frac{\dd q}{q}+w_F\wedge \frac{\dd g}{g}=0
    \in H_p^{m+2}(F).\end{align*} Let $\beta: H_p^{m+2}(F)'\xrightarrow{} H_p^{m+1}(\F(q))$ be the residue map of  Definition \ref{defre} for $q$-adic valuation. By applying the residue map $\beta$, we have \begin{align*}
        \beta\left( w_F\wedge \frac{\dd f}{f}\right)= r w_{\F(q)}&= 0\in H_p^{m+1}(\F(q))\\ w_{\F(q)}&=0\quad \text{(as $p\nmid r$)}.
    \end{align*}  Since $w$ is hyperbolic over $\F(q)$, it follows from  Theorem \ref{thno} that  $ q$ is a norm of $w_F$.\\
 Let $n\geq 2$. This case will be similar to  \cite[Proposition A.15, Step 2]{Appendix}.

\end{proof}
 Now  we prove the norm theorem   in arbitrary number of variables.
\begin{thm}\label{n2th} 
Let $\mathcal{F}$ be a field of characteristic $p$, $F=\F(x_1,\ldots,x_n)$ and $w\in \Omega_{\mathcal{F}}^m$. Let $q\in \mathcal{F}[x_1,\ldots,x_n]$ be a normed irreducible polynomial. Then, the following statements are equivalent:\\(1) $w$ is hyperbolic over $\mathcal{F}(q)$.\\(2) $q$ is a norm of $w_F$.
\end{thm}
\begin{proof}
The case $n=1$ is given by Theorem \ref{thno}.  Now we take  $n>1$. 

(2) $\Longrightarrow$ (1): This direction  will be similar to Theorem \ref{thno}.

(1) $\Longrightarrow$ (2): Renaming the variables if required, we assume that $ q\in K[x_1]$ is a polynomial of degree  $\geq 1$, where $K:=\F(x_2,\ldots,x_n)$. Let $g\in \F[x_2,\ldots,x_n]$ be the leading coefficient of $q\in K[x_1]$.  Since $g$ is a unit in $K$, it follows that $\F(q)\cong K(q)\cong K(g^{-1}q)$. Since $g^{-1}q\in K[x_1]$ is a monic irreducible polynomial, it follows from Theorem \ref{thno} that $ g^{-1}q$ is a norm of $w_{K(x_1)}$.  Therefore, $ g^{p-1} q$ is a norm of $w_{K(x_1)}$. By Proposition \ref{prn}, $q$ is a norm of   $w_{K(x_1)}=w_F$, as desired.
  \end{proof}

\section{Appendix}\label{S6}

{Our aim in this appendix is to give} an  explicit description of  the Teichm\"uller lifting map {that played an important role in {the previous sections}.  The results we  present here are
  from  \cite[chapter 7]{Sc}.   But  for the sake of completeness, we are providing the proofs. 
Throughout  this {section},  $F$ is a complete  field of {characteristic $p>0$} with respect to a discrete valuation $v: F\xrightarrow{} \mathbb{Z}\cup\{\n\}$, {except stated.} Let    $\overline{F}$ be   the residue field  of $ {F}$ and $A_F$ the valuation ring of $F$. Moreover, {let $ A_F^{\times}$ denote}  the group of units of $A_F$. To define the Teichm\"uller lifting, we need the following results.  

\begin{thm}\label{p1t}
     If   the residue field  $\overline{F}$ of the complete field $F$ is perfect, then there exists a unique field
     monomorphism  
     $ \displaystyle{\beta: \overline{F}\xrightarrow{}} F$ such that  ${\rm Im}\displaystyle{(\beta)\subseteq A_F}$, and $\displaystyle{  \overline{\beta({a})}={a}}$ for all $\displaystyle{a\in \overline{F}}$.  
\end{thm}
 
\begin{proof}
      \textbf{Existence:} Let  $ {a}$ be an arbitrary element of $ \overline{F}.$ Since  $\overline{F}$ is  perfect, $ {a}^{p^{-n}}\in \overline{F}$ for all $n\in \mathbb{N}_0.$ Suppose $  a_n\in A_F$ such 
      that \begin{equation}\label{peq1} \overline{a}_n= a^{p^{-n}}\hspace{2mm}\forall\, n\in \mathbb{N}_0.\end{equation} Consequently, we have the following relations: 
      \begin{align*}
    {\overline{a}^p_{n+1}} &=\overline{a}_n \qquad \forall\, n\in \mathbb{N}_0.\\  a^p_{n+1}&\equiv a_n\pmod{P} \qquad (P \hspace{2mm} \text{is   the  
 prime ideal  of }\hspace{1mm} A_F).\end{align*} \textcolor{black}{Moreover, \cite[Chapter 7, Lemma 33]{Sc} implies that:} \begin{equation}\label{peq2}
  a^{p^{n+1}}_{n+1}\equiv a^{p^n}_n\pmod{P^{n+1}} . \end{equation} 
      Let  $k$, $m$, $n\in \mathbb{N}_0$ be such that  $n\geq m\geq k$. Then, we have: \[ a_n^{p^n}- a_m^{p^m}= a_n^{p^n}-a^{p^{n-1}}_{n-1}+ a_{n-1}^{p^{n-1}}- \ldots + a_{m+1}^{p^{m+1}}-a_m^{p^m} \in P^k \quad\text{(by (\ref{peq2}))}.  \] Therefore, $$ \displaystyle{v\left( a_n^{p^n}-a_m^{p^m}\right)\geq k\quad \forall\, n\geq m\geq k}.$$ 
      
Thus, $(a_n^{p^n})_n$ is a {Cauchy} sequence. Since  $F$ is complete, there exists an element  $x_a\in A_F$ such that $ \underset{n\to \infty}{\lim}a_n^{p^n}=x_a$. Moreover, (\ref{peq1}) implies that    $\overline{x}_a=a\in \overline{F}$.  Now we prove that  this $x_a$ does not depend on the choice of  the elements  $a_n$. If possible, let  $a'_n\in A_F$ be such that $ \overline{a'}_n= a^{p^{-n}}$ \ for all $n\in \mathbb{N}_0.$ Then \begin{align*}
          \overline{a'}_n&= \overline{a}_n, \qquad \forall n\in \mathbb{N}_0,\quad \text{i.e.,} \quad a'_n\equiv a_n\pmod{P}\hspace{2mm} \forall n\in \N_0.\end{align*}
Therefore, ${a'}_n^{p^n}\equiv  a_n^{p^n} \pmod{P^{n+1}}$ by \cite[Lemma \hspace{1mm} 33]{Sc}. Hence,  $\underset{n\to \infty}{\lim} {a'}_n^{p^n}= \underset{n\to \infty}{\lim}a_n^{p^n}=x_a.$   Now for each $a\in \overline{F}$, we can define such $x_a\in A_F$. Therefore,   we have the map $ \beta:\overline{F}\xrightarrow{} F$ defined by: \[\beta( a)=x_a\hspace{2mm} \forall\, a\in \overline{F}.\] In particular, if {$ a_n\in A_F$} such that $(\overline{a}_n)^{p^n}= a $ for all $n\in \N_0$, then we have  \begin{align}\label{peq3} \beta(a)=\lim_{n\to \infty} a_n^{p^n} \hspace{2mm} \forall\, a\in \overline{F}.\end{align}  Now we  prove that  $\beta$  is a field homomorphism. Let $ a,b\in \overline{F}$ and  $ a_n,b_n\in A_F$ such that
       $ (\overline{a}_n)^{p^n}=a$ and $(\overline{b}_n)^{p^n}=b$  for all $n\in \N_0$. Then, $({\overline{ a_nb_n}})^{p^n}= ab$ and $ (\overline{a_n+b_n})^{p^n}= a+b$, for all $n\in \N_0$. Therefore, we have  \begin{align*}\beta( ab)&= \lim_{n\to \infty}(a_nb_n)^{p^n}\\ &= \lim_{n\to \infty}a_n^{p^n}\lim_{n\to \infty} b_n^{p^n}
      \qquad \text{(the multiplication  is  a continuous function in $F$) }\\ & = \beta( a)\beta(b). \end{align*}
Again  we get \begin{align*}
          \beta( a+b)= \lim_{n\to \infty} (a_n+b_n)^{p^n}&= \lim_{n\to \infty} (a_n^{p^n}+b_n^{p^n})  \qquad \text{(as char$(F)=p$)}\\&=\lim_{n\to \infty} a_n^{p^n}+ \lim_{n\to \infty} b_n^{p^n}\qquad \text{(the addition  is a continuous function in $F$)}\\ &= \beta(a)+ \beta(b).
      \end{align*}
 Moreover, $ \beta(1)=1.$
      Thus, $\beta$ is  a field homomorphism such that $\overline{\beta(a)}=a$ for all $a\in \overline{F}.$ Hence, $\beta$ is injective.
\vskip1.5mm

\textbf{Uniqueness:} 
 Let $\beta': \overline{F}\xrightarrow{} F$ be    another field monomorphism such that   Im$(\beta')\subseteq A_F $ and $\overline{\beta'(a)} = a$ for all $a\in \overline{F}.$  We will prove that $ \beta=\beta'$. In fact, let $a\in \overline{F}.$  Let $ \beta'( a)= y$ and $ \beta'( a^{p^{-n}})= y_n$ for all $n\in \mathbb{N}_0$. 
Hence, \begin{align}\label{ue1} \beta'( a)&= ({y}_n)^{p^n}\quad \forall n\in \N_0\\\label{ue2}&=y.\end{align} Now by the definition of $\beta',$  we have $ \overline{y}=a$. Using (\ref{ue1}) and (\ref{ue2}) we get \begin{align*}
    \overline{y}_n= \overline{ \beta'(a)}^{p^{-n}}= {\overline{y}}^{p^{-n}}= a^{p^{-n}}  \quad \forall n\in \N_0.
\end{align*} Now by the definition of $ \beta$ from (\ref{peq3}), we have:\begin{align*} \beta(a)&= \lim_{n\to \infty} ({y_n})^{p^n}\quad \text{(since ${a}=(\overline{y_n})^{p^n} \forall\, n\in \mathbb{N}_0$)} \\ &= y\quad \text{(using (\ref{ue1}) and (\ref{ue2}))} \\&= \beta'(a) \quad \forall\, a\in \overline{F}.\end{align*}Thus, the Theorem is proved.
\end{proof}
 The next theorem gives a description of  $F$ using the map $\beta $ of Theorem \ref{p1t}.
\begin{thm}\label{p2t}
    Let $F$, $\overline{F}$ and $\beta$ be as defined in Theorem \ref{p1t}. Then, $ F=\beta(\overline{F})(( \pi)),$ where $\pi$ is a uniformizer.
\end{thm}
\begin{proof}
    Clearly, $\beta(\overline{F})((\pi)) $ is a complete subfield  of $F.$  \textcolor{black}{Moreover, the value groups and residue fields of
 $ F$ and $\beta(\overline{F})((\pi))$ are  same}. Let $ a\in F^*$.  Then, there exists an element $l\in \beta(\overline{F})((\pi))$  such that $ v(a)=v(l).$ Therefore, $v(a/l)=0.$ Then, there exists an element $ c_0\in \beta( \overline{F})$ such that : $$ a/l\equiv c_0\pmod{P_F^{n_1}}, \hspace{1mm}\text{where}\hspace{1mm} P_F=\langle \pi\rangle\hspace{1mm}\text{ is the prime ideal of} \hspace{1mm}
A_F \hspace{1mm}\text{and}\hspace{1mm}  v( a/l-c_0)=n_1\geq 1.$$  Then, $ v\left(\frac{a/l-c_0}{\pi^{n_1}}\right)=0.$  Again there exists an element  $ c_1\in \beta(\overline{F})$  such that$$ \frac{ a/l-c_0}{\pi^{n_1}}\equiv c_1\pmod{P_F^{n_2}}, \hspace{2mm} \text{where}\hspace{1mm} v\left( \frac{a/l-c_0}{\pi^{n_1}}-c_1\right)= n_2\geq 1. $$  Then, $  v\displaystyle{\left( \frac{a/l-c_0-c_1\pi^{n_1}}{\pi^{n_1+n_2}}\right)}=0$.  Similarly, there exists an element $c_2\in \beta(\overline{F})$ such that:$$  \frac{a/l-c_0-c_1\pi^{n_1}}{\pi^{n_1+n_2}}\equiv c_2\pmod{P_F^{n_3}},\hspace{1mm} \text{where}\hspace{1mm} v\left( \frac{a/l-c_0-c_1\pi^{n_1}}{\pi^{n_1+n_2}}-c_2\right)=n_3\geq 1.$$ 
Thus, by induction  we have a  Cauchy sequence $ \left(\sum\limits_{k=0}^{r} c_k\pi^{n_0+n_1+\ldots+n_k}\right)_{r\in \N}$, where $n_0=0$ and $c_i\in \beta(\overline{F})$.
 Since $ \beta(\overline{F})((\pi))   $ is a complete field,  there exists an 
element $b\in \beta(\overline{F})((\pi))$ such that $$b= \sum_{k=0}^{\infty} c_k\pi^{n_0+n_1+n_2+\ldots+n_k}.$$

Since $a/l  \equiv c_0+ c_1\pi^{n_1}+\ldots + c_k\pi^{n_1+n_2+\ldots +n_k}\pmod{ {P_F^{n_1+\ldots+n_{k+1}}}}$ for all $k\in \N$, it follows that  $a/l=b$. Hence, $a\in \beta(\overline{F})((\pi))$ since $l$, $b\in \beta(\overline{F})((\pi))$. Therefore, $F=\beta(\overline{F})((\pi))$.
\end{proof}
\vspace{2mm}

In the sequel, we  \textcolor{black}{generalize} this notion for a complete valued  field $F$ whose residue field $\overline{F}$ is  imperfect.
 If $\overline{F}$ is  imperfect, then it has a nontrivial $p$-basis. Let $ S$ be a  $p$-independent subset of the residue   field $\overline{F}$.  We consider a subset $S'\subseteq A_F$ which is made of exactly one representative for each element in $S.$  Then, $S'$ will be $p$-independent in $F$.
  
\begin{lemma}\label{np1l}
Suppose $K= F(b_{1,1}, b_{1,2}, \ldots, b_{1,r})$, where $ b^p_{1,i}=b_i$ for $b_i\in S'$. So $K$ is a finite extension of $F.$ Hence, we can extend the valuation $v$  uniquely on $K.$ If $\overline{K}$ is the residue field of $K$, then we have:\\(1) $[K:F]= p^r$.\\(2) $[\overline{K}:\overline{F}]=p^r$ and  $\mathcal{T}(K)=\mathcal{T}(F)$, where $ \mathcal{T}(K)$ and $\mathcal{T}(F)$ are the value groups of $K$ and $F$, respectively. 
\end{lemma}
    \begin{proof}
          Since $ K=F( b_{1,1},b_{1,2},\ldots , b_{1,r})$ is generated by $r$ elements over $F$ and each $b_{1,i}$ satisfies a $p$-degree  equation over $F$, it \textcolor{black}{follows} that  $ [K:F]\leq p^r.$  \textcolor{black}{Moreover,}  $\overline{K}= \overline{F}\left(\overline{b}_{1,1}, \overline{b}_{1,2},\ldots, \overline{b}_{1,r}\right) $ and $ \overline{b}_{1,i}^p= \overline{b}_i\in S.$ Since 
$\{\overline{b}_1, \overline{b}_2, \ldots , \overline{b}_r\}$ is a $p$-independent set in $\overline{F}$, we get 
$[\overline{K}:\overline{F}]=p^r.$   \textcolor{black}{Now} by \cite[Chapter 2, Theorem 11]{Sc}, $ [\mathcal{T}(K):\mathcal{T}(F)][\overline{K}:\overline{F}]=[K:F]$. Therefore, $ [K:F]=p^r$ and $[\mathcal{T}(K):\mathcal{T}(F)]=1.$ So the lemma is proved.
   \end{proof}

\begin{defn}\label{mpd}Let $F$ be a field of characteristic $p$, and $F^{p^{-1}}$ the extension of $F$ which is obtained by adjoining the $p$-th roots of all elements of $F$. Then, $\cup_{n\in \mathbb{N}_0} F^{p^{-n}}$ is the minimal perfect extension of $F$. We denote it by $F^{p^{-\n}}$. 
\end{defn}

   \begin{lemma}\label{np2l}
       Let $F$, $\overline{F}$, $S$ and $S'$ be as defined above. \textcolor{black}{Suppose  \textcolor{black}{that} $\overline{F}$}  is  imperfect. Then, there exists a complete extension   $L/F$  such that the residue field of $L$ is $\overline{L}= \overline{F}\left( S^{{p}^{-\infty}}\right)$.  Moreover, if $S$ is a $p$-basis of $\overline{F},$ then $ \overline{L}=\overline{F}^{p^{-\n}}$, the minimal perfect extension of $\overline{F}$. \end{lemma}

\begin{proof}
Let $T$ be a finite subset of $S'$. Then, we can form the field $F\left(T^{\frac{1}{p}}\right)$ as in  Lemma \ref{np1l}.  Consider $K_{1,1}= {\bigcup}_T F\left(T^{\frac{1}{p}}\right)= F\left({S'}^{\frac{1}{p}}\right),$  where the union is over all finite subsets  $ T\subseteq S'.$ Since $K_{1,1}$ is an algebraic extension of the complete field $F$,   we  extend the valuation $v$ uniquely  on $K_{1,1}$, say $v_{K_{1,1}}$.   Now by  Lemma \ref{np1l},    $ \mathcal{T}\left(F\left(T^{\frac{1}{p}}\right)\right)= \mathcal{T}(F), $ for all  finite subsets  $T\subseteq S'$. So $  \mathcal{T}(K_{1,1})=\mathcal{T}(F)$.  Let $ K_1$ denote  the completion  of $ K_{1,1}$ with respect to the valuation $v_{K_{1,1}}$. Then, the valuation $v_{K_{1,1}}$ can be extended on $K_{1}$ uniquely say, $v_{K_1}$  and $ \mathcal{T}(K_1)= \mathcal{T}(K_{1,1})=\mathcal{T}( F)$. Moreover, by Lemma \ref{np1l}, the residue field of $ K_1$ is given by $\overline{K}_1=\overline{F}\left( S^{\frac{1}{p}}\right)$.
   
   Note that   $ S^{\frac{1}{p}}$ is a $p$-independent subset of $  \overline{K}_1= \overline{F}\left(S^{\frac{1}{p}}\right).$     Take  $M_1= {S'}^{\frac{1}{p}}$ and we repeat the same process. Let $K_{2,1}= K_1\left( M_1^{\frac{1}{p}}\right)$ and the completion of $K_{2,1}$ is $K_2.$  Again the residue field of $K_2$ is $\overline{K}_2= \overline{K}_1\left( S^{\frac{1}{p^2}}\right)= \overline{F}\left(S^{\frac{1}{p^2}}\right)$. By repeating this process,  we  get a  chain of fields $F=K_0\subseteq K_1\subseteq K_2\subseteq\ldots$ satisfying the two conditions:
\begin{enumerate}
\item[(1)] The value groups  $\mathcal{T}(K_n)=\mathcal{T}(F)$ for all $n\in \N$.
\item[(2)]  The residue field of $K_n$ is $\overline{K}_n= \overline{F}\left( S^{{\frac{1}{p^n}}}\right)$ for all $n\in \N_0$.
\end{enumerate}
Let us consider the field $L_0= {\bigcup}_{n\in \N} K_n$. Then, its residue field  is $ \overline{L}_0= {\bigcup}_{n\in \N_0} \overline{K}_n= \overline{F}\left( S^{p^{-\infty}}\right) $ and $ \mathcal{T}(L_0)=\mathcal{T}(F).$ Finally, let $L$ be the
 completion of $L_0$, with the residue field $\overline{L}$. Then $\overline{L}= \overline{F}\left( S^{p^{-\infty}}\right)$ and  $\mathcal{T}(L)=\mathcal{T}(F)$. Moreover, if $S$ is a $p$-basis of $ \overline{F},$ then $ \overline{K}_1= \overline{F}\left( S^{\frac{1}{p}}\right)= \left(\overline{F}\right)^{\frac{1}{p}} $ \textcolor{black}{and}     $ S^{\frac{1}{p}}$ is a $p$-basis of $ \overline{K}_1$. Therefore, $ \overline{K}_2= \overline{F}\left( S^{\frac{1}{p^2}}\right)  = \overline{K}_1\left( S^{\frac{1}{p^2}}\right)= \left(\overline{K}_1\right)^{\frac{1}{p}}  =\left( \overline{F}\right)^{\frac{1}{p^2}}.$ Similarly $ \overline{K}_n= \left(\overline{F}\right)^{\frac{1}{p^n}}$ for all $n \in \mathbb{N}$. So $ \overline{L}={\bigcup}_n\overline{K}_n= \left(\overline{F}\right)^{p^{-\infty}},$ which is the minimal perfect extension of $\overline{F}.$ Thus, the lemma is proved.\end{proof}
 \vspace{2mm}
 
Using Lemmas \ref{np1l} and \ref{np2l}, we prove the next theorem which extends Theorem \ref{p1t} {to the case of imperfect residue field}.

\begin{thm}\label{imp}
\textcolor{black}{Suppose \textcolor{black}{that} the} residue field $\overline{F}$ of $F$ is  imperfect. Let $S$ be a $p$-basis of $\overline{F}$ and $S'\subseteq A_F$,  made of  exactly one representative for each element of $S.$ Then, there exists a unique  field monomorphism  $\alpha: \overline{F}\xrightarrow{}F$ such that ${\rm Im}(\alpha)\subseteq A_F$, $\alpha(S)=S'$ and $ \overline{\alpha(a)}=a$ for all $a\in \overline{F}.$    
 \end{thm}

 \begin{proof} \textbf{Existence:} By Lemma \ref{np2l}, there exists a complete field extension $L$ of $F$ whose residue field is $\overline{L}=\left(\overline{F}\right)^{2^{-\infty}},$ the minimal perfect extension of $\overline{F}$, and $L$ contains $F\left(\left(S'\right)^{p^{-\infty}}\right).$ Since $\overline{L}$ is perfect, it follows from Theorem \ref{p1t} that there exists a unique field monomorphsim $\beta: \overline{L}\xrightarrow{} L$ such that $ \overline{ \beta( x)}= x$ for all $x\in \overline{L}.$
\vskip1.5mm
     
     \textbf{Step 1:}  Let $a\in S $ and $x\in S'$ be the representative of $a.$ Now  for each   $n\in \N$, $a^{p^{-n}}\in \overline{L}$ and $ x^{p^{-n}}\in L $  with $ \overline{x^{p^{-n}}}= a^{p^{-n}}$. Therefore$$a= \overline{x^{p^{-n}}}^{p^n}=  \overline{a_n}^{p^n},\;\;\text{where }\;\; a_n= x^{p^{-n}} \; \forall\, n\in \N.$$Then, we have$$\beta(a) = \lim_{n\to \infty} a_n^{p^n}=x,\;\; \text{(by definition of}\;\; \beta\;\; \text{from }\;(\ref{peq3})).$$ \textcolor{black}{Consequently,}  $\beta(S)=S'$.

     \textbf{Step 2:}
     Now we  have to prove that $ \beta(\overline{F})\subseteq F.$   Let $a\in \overline{F}^*.$
     Then, $ a^{p^{-n}}\in \overline{F}^{p^{-n}}=\overline{K}_n$ for all $n\in \mathbb{N}$, where  $\overline{K}_n$ was defined in  Lemma \ref{np2l}. Therefore,   there exist $a_n\in K_n$ such that$$ v_{K_n}(a_n)= 0\hspace{2mm} \text{and}\hspace{2mm}\overline{a}_n= a^{p^{-n}},\hspace{2mm}\forall n\in \N.$$  \textcolor{black}{Recall} \textcolor{black}{that} $K_n$ was constructed by adjoining $p^n$-th roots of $S'$ with $F$, and  then by completion. So for each $n\in \N$, there exists  a sequence $\{a_{k,n}\mid\hspace{1mm}k\in \N\}\subset F$ such that\begin{align*}  a_n&=\lim_{ k\to \infty}a_{k,n}^{p^{-n}}.\end{align*}Since $F$ is complete, we have for all $n\in \N$\begin{align*}(a_n)^{p^n}&= \lim_{k\to \infty} a_{k,n}\in F.\end{align*}    By the definition of $\beta$ from (\ref{peq3}),\hspace{2mm} $\beta(a)=\underset{n\to \infty}{\lim} (a_n)^{p^n}\in F. $ Thus $\beta(\overline{F})\subseteq F\cap A_L= A_F$, where $A_L,A_F$ are the valuation rings of $L$ and $F$ respectively. Therefore, we define the map $\alpha: \overline{F}\xrightarrow{} F$, such that  $\alpha= \beta_{ \overline{F}}.$ Now  $ \alpha$ is  a field monomorphism with   $ \alpha(S)=S'$ and $\overline{ \alpha(a)}=a, $ for all $a\in \overline{F} $. \\

\textbf{ Uniqueness:} Let   $\alpha_1: \overline{F}\xrightarrow{} F$ be a field monomorphism, such that Im$(\alpha_1)\subseteq A_F$, $\alpha_1(S)= S'$ and $\overline{\alpha_1(a)}= a$ for all $a\in \overline{F}.$ \textcolor{black}{Recall}  $\overline{L}= \overline{F}\left(S^{p^{-\infty}}\right)$ and $F\left((S')^{p^{-\infty}}\right)\subseteq L.$  Therefore, we  define a homomorphism $\beta': \overline{L} \xrightarrow{} L$,  \textcolor{black}{given by}  :$$ \beta'_{| \overline{F}}= \alpha_1\hspace{1mm} \text{and} \hspace{1mm} \beta'\left(s^{p^{-n}}\right)= (\alpha_1(s))^{p^{-n}},\hspace{1mm} \text{for all}\hspace{1mm} n\in \N \hspace{1mm}\text{ and for all}\hspace{1mm} s\in S.$$  Since $S$ is $p$-independent in $\overline{F}$, 
 $\beta'$ is a well-defined field homomorphism.  Clearly  Im$(\beta') \subseteq A_L$ and $ \overline{\beta'(x)}=x$ for all $x\in \overline{L}.$ So $\beta'$ satisfies  all the  conditions of $\beta$. By Theorem \ref{p1t}, $\beta$ is unique and it implies that  $\beta=\beta'$. Therefore,   $\alpha=\alpha_1.$ Thus the theorem is proved.    
 \end{proof}
 
The following theorem is a generalization of Theorem \ref{p2t} to the case of imperfect residue field.
 
\begin{thm}\label{im2t}
    Let $F$, $\overline{F}$, $S$,  $S'$ and $\alpha$  be as defined in  Theorem \ref{imp}. Then, 
$F=\alpha(\overline{F})((\pi)),$ where $\pi$ is a uniformizer. 
\end{thm}
\begin{proof}
{The same proof as that of Theorem  \ref{p2t}.} 
\end{proof}

\begin{defn}\label{Tfd}
The monomorphism $\alpha$ of Theorem \ref{imp} is called \textcolor{black}{a Teichm\"uller lifting, which depends on $S$ and $S'$}.
\end{defn}

\begin{cor}\label{T1c}
{The monomorphism $\alpha$ of Theorem \ref{imp} induces a natural homomorphism $\alpha': H_p^{m+1}(\overline{F})\xrightarrow{}  H_p^{m+1}(F)$. Moreover, $\alpha'$ is a monomorphism and for any $x\in A_F$ and $x_i\in A_F\setminus\pi A_F$,} we have: $$\displaystyle{\alpha'\left(  {\overline{x}\frac{\mathrm{d} \overline{x}_{1}}{\overline{x}_{1}}\wedge \ldots \wedge \frac{\mathrm{d}\overline{x}_{m}}{\overline{x}_{m}}}\right) =  { x\frac{\mathrm{d}x_{1}}{x_{1}}\wedge \ldots \wedge \frac{\mathrm{d}x_{m}}{x_{m}}}}.$$
 \end{cor}
 \begin{proof} Let $\alpha': H_p^{m+1}(\overline{F})\xrightarrow{}  H_p^{m+1}(F)$ be the homomorphism induced by the monomorphism $\alpha:\overline{F}\longrightarrow F$. By   Theorem \ref{imp}, for each $x\in A_F$ we have $b_x\in A_F$ such that $v(b_x)>0$ and$$ \alpha( \overline{x})= x+ b_x.$$Moreover, for any $x\in A_F\setminus\pi A_F$,  we have  $$(\alpha(\overline{x}))^{-1}= (x+ b_x)^{-1}= (y_x+c_x),$$  where  $\overline{x} \overline{y}_x=1\in \overline{F} $ and $v(c_x)>0.$ Therefore,  for each $x\in A_F\setminus  \pi A_F$ we get \begin{align}\label{t1em}
          \frac{ x}{\alpha(\overline{x})}= x( x+b_x)^{-1}= x(y_x+c_x)= 1+ d_x,
     \end{align} where $v(d_x)>0$.  Now for $ x\in A_F$ and $x_i\in A_F\setminus\pi A_F$,  we have:   \begin{center}\begin{align*}
         \alpha'\left( {\overline{x} \frac{\mathrm{d}\overline{x}_1}{\overline{x}_1}\wedge \ldots \wedge \frac{\mathrm{d}\overline{x}_m}{\overline{x}_m}}\right)&= { \alpha(\overline{x})\frac{\mathrm{d}(\alpha(\overline{x}_1))}{\alpha(\overline{x}_1)}\wedge \ldots \wedge \frac{\mathrm{d}(\alpha(\overline{x}_m))}{\alpha(\overline{x}_m)}}\\&= { ( x+ b_x)\frac{\mathrm{d}(\alpha(\overline{x}_1))}{\alpha(\overline{x}_1)}\wedge \ldots \wedge \frac{\mathrm{d}(\alpha(\overline{x}_m))}{\alpha(\overline{x}_m)}}\\&= { x\frac{\mathrm{d}(x_1+b_{x_1})}{x_1+ b_{x_1}}\wedge \ldots \wedge \frac{\mathrm{d}(\alpha(\overline{x}_m))}{\alpha(\overline{x}_m)}} \qquad \text{(as $ b_x\in \wp(F)$)}\\ &= {\frac{x x_1}{\alpha(\overline{x}_1)}\frac{\mathrm{d}x_1}{x_1}\wedge \ldots \wedge \frac{\mathrm{d}(\alpha(\overline{x}_m))}{\alpha(\overline{x}_m)} }+ { \frac{xb_{x_1}}{\alpha(\overline{x}_1)}\frac{\mathrm{d}(b_{x_1})}{b_{x_1}}\wedge \ldots \wedge \frac{\mathrm{d}(\alpha(\overline{x}_m))}{\alpha(\overline{x}_m)}}\\&= {\frac{xx_{1}}{\alpha(\overline{x}_1)}\frac{\mathrm{d}x_{1}}{x_1}\wedge \ldots \wedge \frac{\mathrm{d}(\alpha(\overline{x}_m))}{\alpha(\overline{x}_m)}} \qquad \left(\text{as}\hspace{1mm} \frac{xb_{x_1}}{\alpha(\overline{x}_1)}\in \wp( F)\right)\\& = { x(1+d_{x_1})\frac{\mathrm{d}x_1}{x_1}\wedge \ldots \wedge \frac{\mathrm{d}(\alpha(\overline{x}_m))}{\alpha(\overline{x}_m)}} \qquad\text{(by  (\ref{t1em}))} \\&=  { x\frac{ \mathrm{d}x_1}{x_1}\wedge \ldots \wedge \frac{\mathrm{d}(\alpha(\overline{x}_m))}{(\alpha(\overline{x}_m))}}\qquad \text{(since  $xd_{x_1}\in \wp(F)$)}.\end{align*}\end{center}{Repeating the  process for $x_2,\ldots, x_m$,   we  get: } $\displaystyle{\alpha'\left( { \overline{x}\frac{\mathrm{d}\overline{x}_1}{\overline{x}_1}\wedge \ldots \wedge \frac{\mathrm{d}\overline{x}_m}{\overline{x}_m}}\right) = { x\frac{\mathrm{d}x_1}{x_1}\wedge \ldots \wedge \frac{\mathrm{d}x_m}{x_m}}}.$

\noindent{Now} we \textcolor{black}{prove} {that} $\alpha'$ is injective. Let  $x_i\in A_F$ and $x_{i_j} \in A_F\setminus\pi A_F$ be such that  $\displaystyle{\sum_{i} { \overline{x}_i\frac{\mathrm{d}\overline{x}_{i_1}}{\overline{x}_{i_1}}\wedge \ldots \wedge \frac{\mathrm{d}\overline{x}_{i_m}}{\overline{x}_{i_m}}}\in \ker(\alpha')}. $ Then, we have: 
     \begin{align*}
         \alpha'\left( \sum_{i} { \overline{x}_i\frac{\mathrm{d}\overline{x}_{i_1}}{\overline{x}_{i_1}}\wedge \ldots \wedge \frac{\mathrm{d}\overline{x}_{i_m}}{\overline{x}_{i_m}}}\right)= \sum_{i} { x_i\frac{\mathrm{d}x_{i_1}}{x_{i_1}}\wedge \ldots \wedge \frac{\mathrm{d}x_{i_m}}{x_{i_m}}}&=0.  \end{align*}\text{Applying the residue map $\chi_{\pi}$ from Definition \ref{defre} (i), we get}  $\displaystyle{
             \chi_{\pi}\left( \sum_{i} { x_i\frac{\mathrm{d}x_{i_1}}{x_{i_1}}\wedge \ldots \wedge \frac{\mathrm{d}x_{i_m}}{x_{i_m}}}\right)=0}. $
        Using the definition of $\chi_{\pi}$, one gets that $ \displaystyle{\chi_{\pi}\left( \sum_{i} { x_i\frac{\mathrm{d}x_{i_1}}{x_{i_1}}\wedge \ldots \wedge \frac{\mathrm{d}x_{i_m}}{x_{i_m}}}\right)= \sum_{i} { \overline{x}_i\frac{\mathrm{d}\overline{x}_{i_1}}{\overline{x}_{i_1}}\wedge \ldots \wedge \frac{\mathrm{d}\overline{x}_{i_m}}{\overline{x}_{i_m}}}}=0\in H_p^{m+1}(\overline{F}) $.
      Thus $ \alpha'$ is injective. Consequently, we  consider $H_2^{m+1}(\overline{F})$ as a subgroup of $H_p^{m+1}( F).$
 \end{proof}

 \begin{cor}\label{T2c}
     By  Theorem \ref{imp},  $\alpha $ depends  on  $S$, $S'$. However, Corollary \ref{T1c} implies that  $\alpha'$ is independent of the choices of $S$ or $S'$.
 \end{cor}
 \begin{cor}\label{T3c}
     From Corollary \ref{T1c}, we get $ H_p^{m+1}(\overline{F})\cong \alpha'(H_p^{m+1}(\overline{F})) \subseteq H_p^{m+1}(F).$
 \end{cor}
\vspace{2mm}

{Now we define some residue maps which are used several times throughout this note}.

 \begin{defn}\label{defre}
          Let $F$ be a    field of characteristic $p$ and   {$v: F\xrightarrow{}\Z\cup\{\n\}$}  be a discrete valuation with   $\pi$ as a uniformizer. Let $\overline{F}_{\pi}$ and $ F_{\pi}$ be the residue field and the completion of $F$ respectively. Let $U_0$  denote a subgroup of $H_p^{m+1}( F_{\pi})$ generated by the elements  { $ a\frac{\dd x_1}{x_1}\wedge \ldots\wedge \frac{\dd x_m}{x_m}$}, where  $x_i\in F_{\pi}^*$ and $a\in F_{\pi}$ such that  $ v(a)\geq 0$.  By \cite[Theorem 2.5]{I},  the group $U_0$ is additively  generated by the following  forms $$ u \frac{\dd x_1}{x_1}\wedge \ldots\wedge \frac{\dd x_{m}}{x_m}, \hspace{2mm} u \frac{\dd x_1}{x_1}\wedge \ldots\wedge \frac{\dd x_{m-1}}{x_{m-1}}\wedge \frac{\dd \pi}{\pi}, $$ where  $u$, $ x_i\in F_{\pi}$ such that    $v(u)\geq 0$ and $v(x_i)=0$ for all $i$.\\
          \medskip

          \noindent{(i)}  {The residue maps $\zeta : U_0\xrightarrow{}  H_p^{m}(\overline{F})$ and $\chi_{\pi}: U_0\xrightarrow{} H_p^{m+1}( \overline{F})$ are defined  on the generators by}:  
     $$  \zeta \left( u \frac{\dd x_1}{x_1}\wedge \ldots\wedge \frac{\dd x_m}{x_m}+ u' \frac{\dd x'_1}{x'_1}\wedge \ldots\wedge \frac{\dd x'_{m-1}}{x'_{m-1}}\wedge \frac{\dd \pi}{\pi}\right)=  \overline{u'} \frac{\dd \overline{x'}_1}{\overline{x'}_1}\wedge\ldots\wedge \frac{\dd \overline{ x'}_{m-1}}{\overline{x'}_{m-1}}$$ and $$  \chi_{\pi}\left( u \frac{\dd x_1}{x_1}\wedge \ldots\wedge \frac{\dd x_m}{x_m}+ u' \frac{\dd x'_1}{x'_1}\wedge \ldots\wedge \frac{\dd x'_{m-1}}{x'_{m-1}}\wedge \frac{\dd \pi}{\pi}\right)= \overline{u} \frac{\dd \overline{x_1}}{\overline{x_1}}\wedge \ldots\wedge \frac{\overline{x_m}}{\overline{x_m}}.  $$  The maps $ \zeta$, $ \chi_{\pi}$ are well defined by \cite[Theorem 2.5]{Izhboldin}.   The residue map $\zeta$ is independent of the choice of $\pi$, but $\chi_{\pi}$ depends on the choice of a uniformizer.\\
     
    \noindent{(ii)} Let ${H_p^{m+1}}(F)'$ be  the subgroup of $H_p^{m+1}(F)$ generated by the elements {$ a\frac{\dd x_1}{ x_1}\wedge \ldots\wedge \frac{\dd x_m}{x_m} $}, where  $x_i\in F^*$ and $a\in F$ such that  $v(a)\geq 0$. Let $ i: H_p^{m+1}(F)'\xrightarrow{} H_p^{m+1}(F_{\pi})$ be the   restriction of the  homomorphism $j:H_p^{m+1}(F)\xrightarrow{}H_p^{m+1}(F_{\pi})$ on $H_p^{m+1}(F)'$, where $j$ is  induced by {scalar extension}.   Then $ {\rm Im}(i)\subset U_0$. Let $ \beta: H_p^{m+1}(F)'\xrightarrow{}  H_p^{m}(\overline{F}_{\pi})$     and $ \beta_{\pi}:  H_p^{m+1}(F)'\xrightarrow{} H_p^{m+1}(\overline{F}_{\pi})$  be the residue maps defined by:  $$  \beta:= \zeta\circ i, \hspace{2mm}  \beta_{\pi}= \chi_{\pi}\circ i.$$
\end{defn}
 The following lemma deals with the $p$-form, which is crucial to prove Lemma \ref{lepbas}.
\begin{lemma}\label{lan} Let $F$ be a field of characteristic $p$ such that $\{t_1,\ldots,t_n\}$ is $p$-independent over $F$ and $a\in F\setminus F^p$. If the $p$-form $\langle \langle a,t_1,\ldots,t_n\rangle\rangle$ is isotropic, then  we have $$ \langle\langle a, t_1,\ldots, t_n\rangle\rangle_{\rm an}\cong \langle\langle t_1,\ldots,t_n\rangle\rangle$$ where  the notation $\phi_{\rm an}$  denotes the anisotropic part\footnote[8]{See \cite[Definition 2.11]{Hoffmann1} for  explicit  definition of  the anisotropic part of a $p$-form.} of the  $p$-form $\phi$.
    
\end{lemma}

\begin{proof}
    Since $ \{t_1,\ldots, t_n\}$ is  $p$-independent over $K$, the $p$-form $\langle\langle t_1,\ldots,t_n\rangle\rangle$ is anisotropic over $K$.   By \cite[ Corollary 4.10]{Hoffmann1},  the dimension of the $p$-form $\langle\langle a,t_1,\ldots,t_n\rangle\rangle_{\rm an}$  is  $p^m$ for some $m\in \{0,\ldots,n+1\}$.  Since $\langle\langle a,t_1,\ldots,t_n\rangle\rangle$  is isotropic, $m<n+1$.  Therefore, {it follows from  \cite[Proposition 4.19]{Hoffmann1}  that }  \begin{equation}\label{ea1} \psi:=\langle\langle a, t_1,\ldots, t_n\rangle\rangle_{\rm an}\cong   a^i \langle \langle t_1,\ldots,t_n\rangle\rangle \hspace{2mm} \text{ for some $i=0,1,\ldots,p-1$}.\end{equation}
Thus $\dim \psi=p^n$. Hence, by \cite[Corollary 4.10]{Hoffmann1}, we have \begin{align*}
    \psi\cong \begin{cases}
        \langle\langle t_1,\ldots, t_n\rangle\rangle\quad \text{(case 1)}\\ \quad \text{or} \\
        \langle\langle a, t_{i_1},\ldots, t_{i_{n-1}}\rangle\rangle\quad \text{where $1\leq i_1<i_2<\ldots<i_n\leq n$ \quad (case 2).}
    \end{cases} 
\end{align*} In case 1, we are done.  In case $2$ by \cite[Proposition 4.6(i)]{Hoffmann1}, we get   \begin{equation}\label{ea2} \psi\cong a\psi\cong a^2\psi\cong \ldots \cong a^j\psi\quad \text{for any $j\in \N$.}  \end{equation}
By (\ref{ea1}) and (\ref{ea2}), we get  $$  \psi\cong a^{p-i} \psi\cong a^p\langle\langle t_1,\ldots,t_n\rangle\rangle\cong \langle\langle t_1,\ldots,t_n\rangle\rangle.$$
Thus the lemma is proved.
\end{proof}



\end{document}